\ifdefined\pdfminorversion\pdfminorversion=7\relax\fi
\documentclass[12pt]{amsart}
\usepackage[utf8]{inputenc}
\usepackage{amsmath,amssymb,amsthm}
\usepackage{url}
\usepackage{hyperref}
\usepackage{tikz-cd}
\usepackage{stmaryrd}

\usepackage{tikz} 
\usetikzlibrary{external} 
\usepackage{adjustbox}
\definecolor{color0}{HTML}{000000}
\definecolor{color1}{HTML}{000000}
\definecolor{bgColor}{HTML}{ffffff}
\usetikzlibrary {arrows.meta,bending} 

\theoremstyle{definition}
\newtheorem{definition}{Definition}[section]

\newtheorem{observation}[definition]{Observation}
\newtheorem{question}[definition]{Question}
\newtheorem{fact}[definition]{Fact}

\newtheorem{claim}[definition]{Claim}

\theoremstyle{plain}
\newtheorem*{theorem*}{Theorem}

\theoremstyle{plain}
\newtheorem{theorem}[definition]{Theorem}
\newtheorem{proposition}[definition]{Proposition}
\newtheorem{lemma}[definition]{Lemma}
\newtheorem{corollary}[definition]{Corollary}

\newtheorem{remark}[definition]{Remark}

\newcommand\last[2]{\eth_{#1,#2}}

\DeclareMathOperator{\cf}{cf}
\DeclareMathOperator{\cof}{cof}

\DeclareMathOperator{\Tr}{Tr}
\DeclareMathOperator{\tr}{tr}
\DeclareMathOperator{\im}{Im}
\DeclareMathOperator{\otp}{otp}
\DeclareMathOperator{\dom}{dom}

\DeclareMathOperator{\acc}{acc}

\usepackage{amsmath,amssymb,amsthm, xcolor, enumitem, comment, todonotes}

\title{Stepping up to higher dimensions}
\author{Stevo Todorcevic and Jing Zhang}
\date{}

\address[Todorcevic]{\hfill\break Department of Mathematics, University of Toronto,
Canada
\hfill\break stevo@math.toronto.edu}

\address[Zhang]{\hfill\break Department of Mathematics\\
University of North Texas \\
Texas, USA
\hfill\break Jing.Zhang2@unt.edu}

\thanks{
{\it 2010 \emph{MSC}.} 03E02, 03E10 \newline
{The research on this paper is partially supported by the
NSERC grant 455916%
}}

\hypersetup{hidelinks,pdftitle={Stepping up to higher dimensions},pdfauthor={Stevo Todorcevic and Jing Zhang}}
\begin{document}
\begin{abstract}
We design a coloring $c_3: [\omega_3]^4\to \omega$ such that  $c''_3 [\Gamma]^4=\omega$ for any uncountable $\Gamma\subseteq \omega_3.$ The method is based on far-reaching extensions of minimal walks and oscillation theory originally used by the first author in designing the corresponding lower-dimensional analogues $c_2: [\omega_2]^3\to \omega$ and $c_1: [\omega_1]^2\to \omega$ in  \cite{Todorcevic3D} and  \cite{TodorcevicWalks87}, respectively. 


\end{abstract}
\maketitle

\section{Introduction}

Designing structures that exhibit the high level of complexity typically necessary for determining the correct threshold measure for a given combinatorial problem, or merely proving that such structures exist, is a common theme in several areas of mathematics. What we mean by this is perhaps best explained by the example of determining the asymptotic growth of the $k$-dimensional Ramsey function, which in many ways also plays the role of the proper historical framework for the work in this paper.
In this context, for a reason very similar to Ramsey's original motivation in \cite{Ra30}, it is convenient to restrict (without losing generality) from the cubes $X^k$ of sets $X$ to their symmetric versions $[X]^k$ viewed as the collection of $k$-element subsets of  the set $X.$ 
For cardinals (finite or infinite) $\gamma,\beta, k,l$, the notation $\gamma\to (\beta)^k_l$ asserts that for any $f: [\gamma]^k \to l$, there is $B\subseteq \gamma$ with $|B|=\beta$ such that $f\restriction [B]^k$ is constant. For given $\beta, k, l$, let $(\beta)^k_l=\min \{\gamma: \gamma\to (\beta)^k_l\}.$ Erd\H{o}s \cite{erdosprob} has shown via a random coloring that $(n)^2_2$ grows exponentially, leaving two prominent problems open. The first asks if there is an explicit coloring showing that the 2-dimensional Ramsey function  grows exponentially, and the second asks for the corresponding  examples that would show the correct asymptotic behavior of the higher-dimensional Ramsey functions $(n)^k_l$.  Currently, we know essentially only one method for constructing higher-dimensional colorings relevant to the second problem. This is 
the well-known {\it stepping-up lemma} of Erd\H{o}s, Hajnal and Rado (see \cite{ErdosHajnalRado}),
asserting  that for any integers $k\geq 3$, $N$ and $n$, if $N\not\to (n)^k_l$, then $2^N \not\to (2n+k-4)^{k+1}_l$.
In a loose sense, this lemma is used to propagate a counterexample for a $k$-dimensional Ramsey statement to one for a $(k+1)$-dimensional statement, enlarging the source set appropriately.  
 This theorem reduces a well-known question by Erd\H{o}s, Hajnal and Rado about the asymptotic magnitude of the $n$-hypergraph Ramsey number to the 3-hypergraph one (see \cite{ImprovedBound} for more information). 
 
 The stepping-up method of Erd\H{o}s, Hajnal and Rado also works for stepping up counterexamples from dimension 2 to dimension 3 at the expense of doubling the number of colors. Thus, applying it to the 2-dimensional coloring of Erd\H{o}s establishing $\lfloor 2^{n/2}\rfloor\not\to (n)^2_2$, we get that the Ramsey function $(n)^3_4$ is indeed doubly exponential in $n$. 
It is for this reason that the Erd\H{o}s-Hajnal-Rado stepping-up lemma works well in the case of infinitely many colors and, in fact, this was the setting from which the method was imported to the finite case. For example, they used it to compute the exact values of the Ramsey numbers $(k+1)^k_\omega$  as ${\rm exp}_{k-1}(\omega)^+$ for every positive integer $k.$ In structural Ramsey theory, which has deep connections with topological dynamics (see \cite{KPT}) and where colors code relations in a given structure, the possibility of doubling colors when stepping up is not an option, so new stepping-up methods needed to be developed. This paper is a contribution in that direction.
 More precisely, in this paper, we consider stronger witnesses to the failure of Ramsey's theorem, originally proposed by Erd\H{o}s,  Hajnal  and Rado (see \cite{ErdosHajnalRado}, \cite{ErdosHajnalProblem}) and investigate various stepping-up methods to obtain such strong witnesses in higher dimensions. 
 
 Let us review the notation. 
For cardinals (finite or infinite) $\gamma,\beta, k,l$, the notation $\gamma\to [\beta]^k_l$ asserts that for any $f: [\gamma]^k \to l$, there is $B\subseteq \gamma$ with $|B|=\beta$ such that $f''[B]^k \neq l$. It is the negation $\gamma\not\to [\beta]^k_l$ that is usually more interesting and is frequently referred to as the \emph{strong coloring}, although the meaning may slightly vary in specific contexts.  The reason behind their interest is that they become canonical objects when we can prove the companion positive result $\gamma\to (\beta)^k_{m, l}$ for $m=l+1$ (or even better for all $m>l$), where  $\gamma\to (\beta)^k_{m, l}$  is the natural variation on the ordinary Erd\H{o}s-Rado symbol saying that for any $f: [\gamma]^k \to m$, there is $B\subseteq \gamma$ with $|B|=\beta$ such that $|f''[B]^k| \leq  l$.
Todorcevic \cite{TodorcevicWalks87} proved that in ZFC, $\omega_1\not\to [\omega_1]^2_{\omega_1}$, completing a long series of advances by Sierpinski \cite{Sierpinski}, Erd\H{o}s-Hajnal-Rado \cite{ErdosHajnalRado}, Laver \cite{LaverParition}, Galvin-Shelah \cite{GalvinShelah} and others. In doing so, he introduced two powerful new methods, walks on ordinals and oscillation theory.

Whether there are corresponding higher-dimensional strong colorings becomes a very natural question. An early result is due to Shore \cite{Shore}, who showed that in G\"{o}del's constructible universe $L$, for any $n\in \omega$, $\omega_n\not\to [\omega_1]^{n+1}_{\omega_1}$. His method involves stepping up using a Kurepa tree with many branches and a particular strong coloring on $[\omega_1]^2$ whose existence follows from CH but not ZFC alone. Todorcevic \cite{TodorcevicThesis} in his thesis achieved the same stepping-up result without using CH, at the cost of using a Kurepa tree that does not contain an Aronszajn subtree. The existence of such a tree does not follow from the existence of a Kurepa tree, by a result of Komj\'{a}th \cite{Komjath} (later improved in \cite{HLRandST}). Later, Todorcevic \cite{TodorcevicSpecial} found another argument to step up strong colorings using special Aronszajn trees that contain neither Aronszajn nor Cantor subtrees. As a result, he improved Shore's theorem by showing that if for all $n\in \omega$, $\square_{\omega_n}$ holds, then $\omega_n\not\to [\omega_1]^{n+1}_{\omega_1}$. However, a major caveat of these results is that they are not ZFC theorems. In fact, for $n\geq 2$, $\omega_n\not\to [\omega_1]^{n+1}_{\omega_1}$ contradicts the following version of Chang's conjecture: $(\omega_n, \omega_1)\twoheadrightarrow (\omega_1, \omega)$ which asserts that for any $f: [\omega_n]^{<\omega}\to \omega_1$, there is some $A\in [\omega_n]^{\omega_1}$ such that $|f''[A]^{<\omega}|\leq \aleph_0$. This version of Chang's conjecture is known to be consistent relative to the existence of large cardinals.

In the 1990s, Todorcevic \cite{Todorcevic3D} proved that in ZFC, 
\begin{enumerate}
\item $\omega_2\not\to[\omega_1]^3_{\omega_1}$ is equivalent to the failure of Chang's conjecture and 
\item $\omega_2\not\to [\omega_1]^3_{\omega}$.
\end{enumerate}
These results are the main motivation of our current work. Specifically, it is natural to ask if similar results hold for higher-dimensional colorings on larger cardinals (see \cite{partitionrelationsbook,higherwalks, BergfalkHigher} for explicit mentions). These questions are part of the emerging and exciting research program that investigates higher-dimensional combinatorics of the cardinals $\omega_k$ for positive integers $k$. Let us mention that the work in  \cite{TodorcevicWalks87}  and \cite{Todorcevic3D} inspired many follow-up investigations with connections to other fields, like homological algebra \cite{BergfalkHigher, BergfalkLambieHansonZhang}, group theory \cite{FeldmanRinot,COV2024,PR2025}, functional analysis \cite{ALAT2006, LAT2011, HRST2021} and many others.

In light of the previous results, the next natural question to attack is whether $\omega_3\not\to [\omega_1]^4_{\omega}$ in ZFC. Our main contribution is a positive resolution of this question. We give stepping-up arguments from various weak combinatorial assumptions at the level of $\omega_3$. More specifically, $\omega_3\not\to [\omega_1]^4_{\omega}$ holds if 
	\begin{enumerate}
	\item either there is an $\omega_2$-special $\omega_3$-Aronszajn tree, 
	\item or there is a weak Kurepa tree with $\geq \omega_3$ many branches,
	\item or the continuum is greater than or equal to $\omega_3$.
	\end{enumerate}
	
Compared to previous results in \cite{TodorcevicThesis, TodorcevicSpecial}, the properties of the trees that we use in the stepping-up argument are significantly weakened. In particular, we no longer require the tree to have no Cantor or Aronszajn subtrees. The stepping-up method when the continuum is large is an elaborate combination of an extension of the ordinary square-bracket operation of  \cite{TodorcevicWalks87} (see also \cite{walksbook}) to recover multiple points simultaneously and the oscillation of lower traces first discovered in \cite{MooreLspace}.

The reader may already have noticed that the number of colors is $\omega$ instead of $\omega_1$. This is due to the fact that suitable versions of Chang's conjecture imply there is no such strong coloring with $\omega_1$ many colors.

The main result of this paper is the following theorem in ZFC.

\begin{theorem}\label{theorem: main}
There exists $c: [\omega_3]^4\to \omega$ such that $c''[\Gamma]^4 = \omega$  for any uncountable $\Gamma\subseteq \omega_3$.
\end{theorem}

If we consider $\omega_1$ many colors instead, we are also able to establish an equivalence with a version of Chang's conjecture under a mild cardinal arithmetic hypothesis. 

\begin{theorem}\label{theorem: Chang}
If $n\geq 2$ and there exists a weak Kurepa tree with $\geq \omega_n$ many uncountable branches, then the following are equivalent: 
	\begin{enumerate}
	\item $(\omega_n, \omega_1)\not\twoheadrightarrow (\omega_1, \omega)$
	\item $\omega_n\not\to [\omega_1]^{n+1}_{\omega_1}$.
    \end{enumerate}	 
\end{theorem}
While we currently don't know if the hypothesis in Theorem \ref{theorem: Chang} could be removed, it is a consequence of a particular cardinal arithmetic configuration, as shown by Baumgartner \cite{BaumgartnerAlmostDisjoint}.

The paper is organized as follows. 
	\begin{enumerate}
	\item Sections \ref{section: prelim} and \ref{section: pruning} cover some preliminary facts useful for later sections.
	\item Section \ref{section: omega_2} contains a proof of $\omega_2\not\to [\omega_1]^3_\omega$ with new ingredients, which are crucially used in the next section.
	\item Section \ref{section: weaksquare} contains a proof of Theorem \ref{theorem: main} additionally conditioned on the existence of an $\omega_2$-special $\omega_3$-Aronszajn tree.
	\item Section \ref{section: weakKurepa} contains a proof of Theorem \ref{theorem: main} additionally conditioned on the existence of a weak Kurepa tree with $\geq \omega_3$ many uncountable branches.
	\item Section \ref{section: continuum} contains a proof of Theorem \ref{theorem: main} additionally conditioned on the continuum being greater than or equal to $\omega_3$. Combining this with the results from the previous sections, we deduce Theorem \ref{theorem: main}.
	\item Section \ref{section: cor} contains some remarks regarding the method of the proof. 
	\item Section \ref{section: WKTandCC} supplies a proof of Theorem \ref{theorem: Chang} as well as some extensions. 
	\item Section \ref{section: questions} finishes with some concluding remarks and open questions.
	\end{enumerate}
	
Sections \ref{section: omega_2}, \ref{section: weaksquare}, \ref{section: weakKurepa}, and \ref{section: continuum} build on one another in order of increasing complexity. 

\section{Preliminaries and notation}\label{section: prelim}
Let us collect the basic notation and facts concerning the method of \emph{walks on ordinals} developed by Todorcevic in \cite{TodorcevicWalks87}. Fix a regular cardinal $\kappa>\omega$ and a $C$-sequence $\vec{C}=\langle C_\alpha\subseteq \alpha: \alpha<\kappa\rangle$, where $C_\alpha$ is a closed unbounded subset of $\alpha$ when $\alpha$ is a limit ordinal and $C_\alpha= \{\alpha^-\}$ when $\alpha$ is a successor ordinal. For a set $A$ of ordinals, $\min(A\setminus\xi)$ denotes the least member of $A$ greater than or equal to $\xi$, whenever it exists.

\begin{definition}[Todorcevic \cite{TodorcevicWalks87}, {\cite[Chapter 10.3]{walksbook}}]\label{def: walks}
From $\vec{C}$, we derive maps $\Tr: [\kappa]^2\to {}^\omega \kappa$, $\rho_2: [\kappa]^2\to \omega$, $\tr: [\kappa]^2\to {}^{<\omega}\kappa$, $\lambda: [\kappa]^2\to \kappa$ and $\chi: [\kappa]^3\to \omega$ as follows. Let $(\alpha,\beta, \gamma)\in [\kappa]^3$ be arbitrary.

	\begin{itemize}
	\item Write $v_n=\Tr(\alpha,\beta)(n)$. Set $v_0=\beta$ and, for $n\in\omega$, put
    \[
      v_{n+1}=\begin{cases}
        \min(C_{v_n}\setminus\alpha), & v_n>\alpha,\\
        \alpha, & v_n=\alpha.
      \end{cases}
    \]
    \item $\rho_2(\alpha,\beta)=\min \{n\in \omega: \Tr(\alpha,\beta)(n)=\alpha\}$
	\item $\tr(\alpha,\beta)=\Tr(\alpha,\beta)\restriction \rho_2(\alpha,\beta)$
	\item $\lambda(\alpha,\beta)=\max \{ \sup (C_{\Tr(\alpha,\beta)(i)}\cap \alpha): i<\rho_2(\alpha,\beta)\}$.
	\item $\chi(\alpha,\beta,\gamma)= \max \{i\in \omega: \forall j\leq i, \Tr(\beta, \gamma)(j)=\Tr(\alpha,\gamma)(j)\}$.
	\end{itemize}
\end{definition}

We need two basic facts regarding the functions defined above. 

\begin{fact}
For $\alpha<\beta <\kappa$, if $i<\rho_2(\alpha,\beta)-1$, then $\sup (C_{\Tr(\alpha,\beta)(i)}\cap \alpha)<\alpha$. 
\end{fact}

\begin{fact}\label{fact: passing}
Whenever $\lambda(\beta,\gamma)<\alpha<\beta<\gamma<\kappa$, we have \[\tr(\alpha,\gamma)=\tr(\beta,\gamma)^\frown \tr(\alpha,\beta).\]
\end{fact}

Let us extend the definition of $\lambda$ as follows. 
\begin{definition}\label{def: lambdak}
Given $(\alpha,\beta)\in [\kappa]^2$ and $k\leq \rho_2(\alpha,\beta)$, let $$\lambda_k(\alpha,\beta)=\max\{\sup(C_{\Tr(\alpha,\beta)(i)}\cap\alpha):i<k\}.$$
\end{definition}

The following definition captures the first point in a walk (in decreasing order) whose club has the destination ordinal as a limit point, when such a point exists. 

\begin{definition}[{\cite[Definition 2.10]{transformation}}]
For every $(\alpha,\beta)\in[\kappa]^2$, we define an ordinal $\last{\alpha}{\beta}\in[\alpha,\beta]$ by
$$\last{\alpha}{\beta}:=\begin{cases}
\alpha,&\text{if }\lambda(\alpha,\beta)<\alpha;\\
\min(\im(\tr(\alpha,\beta))),&\text{otherwise}.
\end{cases}$$

\end{definition}

\begin{definition}[{\cite[Definition 8.14]{walksbook}}]
Let $\Gamma\subseteq \kappa$ and $\theta\in \kappa\cap \mathrm{cof}(>\omega)$ be given. We say $\langle C_\alpha: \alpha \in \Gamma\rangle$ 

\begin{itemize}
\item is \emph{stationary} at $\theta$ if $\bigcup_{\alpha\in \Gamma} \acc(C_\alpha) \cap \theta$ is a stationary subset of $\theta$;
\item is \emph{trivial} at $\theta$ if there exists a club $D\subseteq \theta$ such that for any $\alpha<\theta$, there exists some $\beta\in \Gamma$ such that $D\cap \alpha \subseteq C_\beta$.
\end{itemize}
\end{definition}

\begin{definition}[{\cite[Chapter 8.1]{walksbook}}]
Let $D,E$ be two closed sets of ordinals. We put $osc(D,E)=\left|\bigl(D-(\sup(D\cap E)+1)\bigr)/{\sim}\right|$, where for $a\leq b$ in $D-(\sup(D\cap E)+1)$, $a\sim b$ if $[a,b]\cap E=\emptyset$. In the colorings below, if the oscillation is infinite or the indicated difference is negative, we output $0$.
\end{definition}

\begin{proposition}[{\cite[Lemma 10.3.4]{walksbook}}]\label{prop: osc}
If there is some $\Gamma\subseteq \kappa$ with $\theta=\sup \Gamma$ such that $\cf(\theta)>\omega$ and $\langle C_\alpha: \alpha\in \Gamma\rangle$ is stationary and non-trivial at $\theta$, then for every positive integer $n$, there are $\alpha<\beta<\gamma\in \Gamma$ such that $C_\beta\cap C_\gamma\subseteq\alpha$ and $osc(C_\beta-\alpha, C_\gamma-\alpha)=n$.
\end{proposition}

\begin{definition}\label{def: delta}
Let $s,t$ be two functions whose domains are ordinals. If they differ on their common domain, let $\Delta(s,t)=\min\{\xi\in\dom(s)\cap\dom(t):s(\xi)\neq t(\xi)\}$.
\end{definition}

Given $A\subseteq \mathrm{Ord}$, we let $f_A: \otp(A) \to A$ denote the unique order-preserving enumeration. Let us record a useful special case of Definition \ref{def: delta}.

\begin{definition}\label{def: deltaclub}
For two sets $C,D \subseteq \mathrm{Ord}$, define $$\Delta(C,D)=\Delta(f_C, f_D) = \min \{\xi\in \otp(C)\cap \otp(D): f_C(\xi)\neq f_D(\xi)\}.$$
\end{definition}

Note that if $C$ and $D$ are closed sets, then $\Delta(C,D)$ is either zero or a successor ordinal whenever it is defined. The following relative of the $\Delta$-function will be crucially used in Section \ref{section: weakKurepa}.

\begin{definition}\label{def: deltaprime}
For $\bar{\beta}\leq \beta$ and two closed unbounded subsets $C\subseteq \bar{\beta}, D\subseteq \beta$, define $\Delta'(C,D)=\max\{\xi \leq   \max (\bar{\beta} \cap D): \otp(C\cap \xi) = \otp(D\cap \xi)\}$. If the associated club sets are clear from the context, we will sometimes write $\Delta'(\bar{\beta},\beta)$ in place of $\Delta'(C,D)$.
\end{definition}

\begin{definition}
We say $\langle C_\alpha: \alpha\in \Gamma\rangle$ is \emph{strongly non-trivial} at $\theta$ if there exist $\xi<\cf(\theta)$ and an unbounded set $A\subseteq \theta$ of order type $\cf(\theta)$ such that for any $\alpha\in A$, for all $\beta\in \Gamma-(\alpha+1)$, $\otp(C_\beta \cap \alpha)\leq \xi$.
\end{definition}

The reason why we consider such a variant is that it saves us one coordinate when we are trying to generate oscillation patterns.

\begin{proposition}\label{prop: deltaprime}
If there is some $\Gamma\subseteq \kappa$ with $\lambda=\sup \Gamma$ and $|\Gamma|=\cf(\lambda)$ such that $\cf(\lambda)>\omega$ and $\langle C_\alpha: \alpha\in \Gamma\rangle$ is stationary and strongly non-trivial at $\lambda$, then for every positive integer $n$, there are $\bar{\beta}<\beta\in \Gamma$ such that $$osc(C_\beta-\Delta'(C_{\bar{\beta}}, C_\beta), C_{\bar{\beta}}-\Delta'(C_{\bar{\beta}}, C_\beta))=n+1.$$
\end{proposition}

\begin{proof}
Fix a positive integer $n$.
Let $\nu=\cf(\lambda)$.
Let $\vec{M}=\langle M_\alpha: \alpha<\nu\rangle$ be an internally approachable sequence of elementary submodels of $H(\theta)$ of size $<\nu$ containing relevant objects with $M_\alpha\cap \nu \in \nu$ for each $\alpha<\nu$. Let $E=\{\sup(M_\alpha\cap\lambda): \alpha<\nu\}$. Then $E$ is a club in $\lambda$ of order type $\nu$. 

Since $\langle C_\alpha: \alpha\in \Gamma\rangle$ is stationary, we can find an elementary submodel $N\prec H((2^\theta)^+)$ containing relevant objects including $\vec{M}, \lambda$ and $\vec{C}$ such that 
	\begin{itemize}
	\item $|N|<\nu$,
	\item $N\cap \nu \in \nu$,
	\item $\sup(N\cap\lambda)\in \acc(C_\beta)$ for some $\beta\in \Gamma$.
	\end{itemize}
	Let $\lambda_N= \sup(N\cap\lambda)$.
By the fact that $\langle C_\alpha: \alpha\in \Gamma\rangle$ is strongly non-trivial, we know that $\otp(C_\beta \cap \lambda_N)\in N\cap \nu$. As $E\in N$ and $\otp(E\cap N)=N\cap \nu$, we know that for any $i<\lambda_N$, $ E\cap N - i \not \subseteq C_\beta \cap \lambda_N$.
By the order type calculation, there exist increasing sequences $\langle \lambda_i\in \lambda_N : i\in \omega\rangle$ and $\langle \xi_i\in N\cap \nu: i\in \omega\rangle$ such that for all $i\in \omega$,
	\begin{itemize} 
	\item  $\lambda_i\in E\cap N-C_{\beta}$,
	\item $(\lambda_i, \lambda_{i+1})\cap C_{\beta}\neq \emptyset$,
	\item $\otp(C_\beta \cap \lambda_0)=\xi_0$,
	\item $\otp((\lambda_i, \lambda_{i+1})\cap C_{\beta})=\xi_{i+1}$.
	\end{itemize}
For simplicity, let us denote the model witnessing $\lambda_i\in E$ by $M'_i$ for $i\in \omega$. We may also assume that $\otp(C_\beta \cap \lambda_N)\in M_0'$. Choose $\bar{\lambda}\in M'_0\cap\lambda$ with $\bar{\lambda}>\sup(C_\beta\cap\lambda_0)$.
Let $\exists^Q I$ abbreviate: for any $\xi<\lambda$, there exists an interval $I \subseteq \lambda$ with $\min I > \xi$. The standard reflection argument using the elementarity of $M_0',\cdots, M'_{n},N$ (see \cite[Lemma 8.1.2]{walksbook} for example) gives that $\exists^Q I_1 \exists^Q I_2 \cdots  \exists^Q I_n \psi(I_1,\cdots, I_n)$ where $\psi(I_1,\cdots, I_n)$ asserts that there is some $\bar{\beta}\in \Gamma$ these conditions.
	\begin{itemize}
	\item $\otp(C_{\bar{\beta}}\cap \bar{\lambda}) =\xi_0$,
	\item $\otp(C_{\bar{\beta}}\cap I_i)=\xi_i$ for $i=1,2, \cdots, n-1$,
	\item $C_{\bar{\beta}}\subseteq \bar{\lambda} \cup \bigcup_{i=1}^{n} I_i$ and $C_{\bar{\beta}}\cap I_n\neq\emptyset$.
	\end{itemize}
Observe that the parameters involved in the sentence $$\exists^Q I_1 \exists^Q I_2 \cdots  \exists^Q I_n \psi(I_1,\cdots, I_n)$$ are $\Gamma$, $\vec{C}$, $\bar{\lambda}$, $\{\xi_i: 0\leq i \leq n-1\}$, which are all in $M_0'$.
	
Let us briefly sketch how the reflection works. For each $1\leq i \leq n$, find some $\bar{\lambda}_i\in M'_i\cap \lambda$ such that $\max C_{\beta}\cap \lambda_i<\bar{\lambda}_i$. Consider the following intervals $K_i=[\lambda_{i-1}, \bar{\lambda}_i]$ for $1\leq i < n$ and $K_n=[\lambda_{n-1},\lambda)$. Note that $K_i\in M'_i$ for each $1\leq i\leq n$. Notice that $\psi(K_1,K_2, \cdots, K_n)$ is witnessed by $\beta$. By the elementarity of $M'_{n-1}$, we can conclude that $\exists^Q I_n \psi(K_1,\cdots, I_n)$. Repeating the same argument for the sentence $\exists^Q I_n \psi(K_1,\cdots, I_n)$ and $M'_{n-2}$, we get $\exists^Q I_{n-1}\exists^Q I_n \psi(K_1,\cdots, I_{n-1},I_n)$. After repeating this argument $n$ times, we arrive at the conclusion that $\exists^Q I_1 \exists^Q I_2 \cdots  \exists^Q I_n \psi(I_1,\cdots, I_n)$ holds.

Let $J_1 = (\lambda_0, \max C_{\beta}\cap \lambda_2+1)$, $J_2 = (\lambda_2, \max C_{\beta}\cap \lambda_4+1)$, $\cdots$, $J_n = (\lambda_{2(n-1)}, \max C_{\beta}\cap \lambda_{2n}+1)$. Reflecting the statements above in models $M_2', \cdots, M_{2n}'$, we get in $N$ intervals $I_1<\cdots<I_n$ and the witnessing $\bar{\beta}\in \Gamma$ satisfying the following:  $$J_1<I_1<\lambda_2 <J_2 < I_2<\lambda_4< J_3 < I_3 <\lambda_6 <\cdots < J_n<I_n< \lambda_{2n}$$ and $\otp(C_{\bar{\beta}}\cap \bar{\lambda})=\xi_0$.
The points of $C_\beta$ in each $J_j$ form one class, and the points above $I_n$ form one more, since $C_{\bar{\beta}}$ meets every $I_j$. Thus $osc(C_{\beta}-\lambda_0, C_{\bar{\beta}}-\lambda_0)=n+1$. 

	\begin{claim}\label{claim: computedeltaprime}
	$\Delta'(C_{\bar{\beta}}, C_{\beta})=\min(C_\beta\cap J_1)$.
	\end{claim}
	
	\begin{proof}[Proof of the claim]
	Note first that $\max \bar{\beta}\cap C_\beta = \max C_{\beta}\cap \lambda_{2n}= \max C_\beta \cap J_n$. Let $\eta=\min(C_\beta\cap J_1) < \max \bar{\beta}\cap C_\beta$. Clearly, 
		\begin{itemize}
		\item $\eta \leq \Delta'(C_{\bar{\beta}}, C_{\beta})$ since $\otp(C_{\bar{\beta}}\cap \eta)=\otp(C_{\bar{\beta}}\cap \bar{\lambda})=\xi_0 = \otp(C_\beta\cap \bar{\lambda})=\otp(C_\beta \cap \eta)$.
		\item For the other direction, we show that for any $\tau\in (\eta, \max \bar{\beta}\cap C_\beta]=(\eta, \max C_\beta\cap J_n]$, $\otp(C_{\bar{\beta}}\cap \tau)<\otp(C_\beta \cap \tau)$. Let us consider the following cases: there is $1\leq j\leq n$
				\begin{enumerate}
				\item $\tau \in (\min(C_\beta\cap J_j), \max J_j]$. We know that $\otp(C_{\bar{\beta}}\cap \tau)= \xi_0 + \cdots + \xi_{j-1} < \xi_0 + \cdots + \xi_{2j-2}+1 \leq \otp(C_\beta\cap \tau)$.
				\item $\tau \in (\max J_j, \min(C_\beta\cap J_{j+1})]$ where $j\leq n-1$. We know that $\otp(C_{\bar{\beta}}\cap \tau)\leq \xi_0 + \cdots + \xi_{j} < \xi_0 + \cdots + \xi_{2j-2}+  (\xi_{2j-1}+\xi_{2j})=\otp(C_\beta\cap \tau)$.
				\end{enumerate}
		\end{itemize}
	\end{proof}
Neither club meets $[\lambda_0,\min(C_\beta\cap J_1))$, so the claim gives $$osc(C_{\beta}-\Delta'(C_{\bar{\beta}}, C_{\beta}),C_{\bar{\beta}}-\Delta'(C_{\bar{\beta}}, C_{\beta}))=n+1.$$
\end{proof} 
Moore \cite{MooreLspace} proved the existence of a coloring $h: [\omega_1]^2\to \omega_1$ witnessing $\omega_1\not\to[\omega_1, \omega_1]^2_{\omega_1}$. Namely, for any $A,B\in [\omega_1]^{\omega_1}$, $h'' A\otimes B = h''\{(\alpha,\beta): \alpha<\beta, \alpha\in A, \beta\in B\}=\omega_1$. Let us supply a convenient strengthening of the conclusion.

\begin{lemma}\label{lemma: strengthenedMoore}
For any $h$ witnessing $\omega_1\not\to[\omega_1, \omega_1]^2_{\omega_1}$, $h$ satisfies the following strengthening: for any $\tau <\omega_1$ and any $A,B\in [\omega_1]^{\omega_1}$, there are $\alpha\in A$ and uncountably many $\beta\in B$ such that $h(\alpha,\beta)=\tau$.
\end{lemma}

\begin{proof}
Suppose $\tau, A, B$ witness the negation for the sake of contradiction: namely, for each $\alpha\in A$, there is $\beta_\alpha\in \omega_1$ such that for any $\xi\in B-\beta_\alpha$, $h(\alpha, \xi)\neq \tau$. Recursively build increasing sequences $\langle a'_\alpha\in A: \alpha<\omega_1\rangle$ and $\langle b'_\alpha\in B: \alpha<\omega_1\rangle$ 
such that for any $\alpha<\omega_1$,
		\begin{itemize}
		\item $a'_\alpha<b'_\alpha<a'_{\alpha+1}$ and 
		\item $\sup_{i\leq \alpha}\beta_{a'_i} +1<b'_\alpha$.
		\end{itemize}

Let $A'=\{a'_\alpha: \alpha<\omega_1\}$ and $B'=\{b'_\alpha: \alpha<\omega_1\}$.
 By the property of $h$ applied to $A'\otimes B'$, we can find $(a,b)\in A' \otimes B'$ such that $h(a,b)=\tau$. By the construction of $A'$ and $B'$, there are $\alpha\leq \beta$ such that $a=a'_\alpha$ and $b=b'_\beta$. As a result of the choice of $b'_\beta$, $h(a,b)\neq \tau$. This is a contradiction.
\end{proof}

We recall a construction of a well-known $\rho$-function at the level of $\omega_2$ and some of its properties from \cite{walksbook}.

\begin{definition}\label{def:walk-rank}
Fix a $C$-sequence $\langle C_\beta:\beta<\omega_2\rangle$ with $\otp(C_\beta)=\cf(\beta)$ for every $\beta<\omega_2$. Define $\rho:[\omega_2]^2\to\omega_1$ recursively. Put $\rho(\alpha,\alpha)=0$. For $\alpha<\beta<\omega_2$, let $s=\min(C_\beta\setminus\alpha)$ and set
\begin{equation}\label{eq:walk-rank}
\begin{aligned}
\rho(\alpha,\beta)=\sup\bigl(&\{\otp(C_\beta\cap\alpha),\rho(\alpha,s)\}\\
                            &{}\cup\{\rho(\xi,\alpha):\xi\in C_\beta\cap\alpha\}\bigr).
\end{aligned}
\end{equation}
\end{definition}

\begin{lemma}\label{lem:walk-rank}
For $\beta<\omega_2$ and $\nu<\omega_1$, the set $G_\nu(\beta)=\{\alpha<\beta:\rho(\alpha,\beta)\leq\nu\}$ is countable. For $\alpha<\beta<\gamma<\omega_2$,
\[
\begin{aligned}
\rho(\alpha,\gamma)&\leq\max\{\rho(\alpha,\beta),\rho(\beta,\gamma)\},\\
\rho(\alpha,\beta)&\leq\max\{\rho(\alpha,\gamma),\rho(\beta,\gamma)\}.
\end{aligned}
\]
If $y<\tau<\omega_2$ and $v>y$ is a node on the walk from $\tau$ to $y$, then
\[
 \rho(y,\tau)\geq\rho(y,v)\geq\otp(C_v\cap y).
\]
If $u<\tau<\omega_2$, $u$ is a limit ordinal and $v=\last{u}{\tau}>u$, then for all sufficiently large $y<u$,
\begin{equation}\label{eq:walk-tail}
   \rho(y,\tau)\geq\rho(y,v)\geq\otp(C_v\cap y).
\end{equation}
\end{lemma}

\begin{proof}
For countability and subadditivity see \cite[p.~271, (9.1.2) and Lemmas~9.1.1--9.1.2]{walksbook} and the proof of Lemma~7.2.2 there. Iterating \eqref{eq:walk-rank} along the walk to $v$ gives $\rho(y,\tau)\geq\rho(y,v)$, and the order type $\otp(C_v\cap y)$ occurs in the definition of $\rho(y,v)$. For the last assertion, take $y<u$ above $\sup(C_\beta\cap u)$ for every node $\beta$ before $v$ on the walk from $\tau$ to $u$. Each of these finitely many bounds is below $u$, so the walk from $\tau$ to $y$ also passes through $v$.
\end{proof}

\section{Monotone, increasing, and locally small functions and basic pruning}\label{section: pruning}

Fix a regular cardinal $\kappa\geq \omega_1$ and a function $f: [\kappa^+]^2 \to \kappa$ for the rest of the section. 

\begin{definition}
Given $\Gamma\in [\kappa^+]^{\omega_1}$, 
	\begin{enumerate}
	\item let $\lambda_{\Gamma}=\sup \{f(\alpha,\beta): \alpha<\beta\in \Gamma\}$;
	\item we say $f$ is \emph{monotone} on $\Gamma$ if for any $\alpha<\alpha'<\beta \in \Gamma$, $f(\alpha,\beta)\leq f(\alpha',\beta)$.
	\item we say that $f$ is \emph{locally small} on $\Gamma$ if for every $\alpha\in \Gamma$, $A_{\alpha}=_{def} \{f(\alpha,\beta): \beta\in \Gamma,\ \alpha<\beta\}$ is countable;
	\item we say that $f$ is \emph{increasing} on $\Gamma$ if there exists an unbounded $D\subseteq \sup \Gamma$ such that whenever $\alpha<\delta<\alpha'<\beta$ with $\alpha,\alpha',\beta\in \Gamma$ and $\delta\in D$, we have that $f(\alpha,\beta)<f(\alpha',\beta)$.
	\end{enumerate}
\end{definition}

Fix $\Gamma \in [\kappa^+]^{\omega_1}$ of order type $\omega_1$ and assume that $f$ is monotone, locally small and increasing on $\Gamma$. Note that for any $\Sigma\in [\Gamma]^{\omega_1}$, $f$ is monotone, locally small and increasing on $\Sigma$.

\begin{lemma}\label{lemma: least}
	There exists an uncountable set $\Sigma\subseteq \Gamma$ such that for each uncountable $\Gamma'\subseteq \Sigma$, $\lambda_{\Gamma'} = \lambda_{\Sigma}$.
	\end{lemma}
	\begin{proof}
	Take $\Sigma$ such that $\lambda_\Sigma$ is the least ordinal in $\{\lambda_{\Gamma'}: \Gamma'\in [\Gamma]^{\omega_1}\}$. Such a $\Sigma$ is as desired.
	\end{proof}

Shrinking to an uncountable subset if necessary, we may assume $\Gamma$ already satisfies the conclusion of Lemma \ref{lemma: least}. Furthermore, we can assume for each $\alpha\in \Gamma$, if $\tau\in A_\alpha$, then there are uncountably many $\beta\in \Gamma$ such that $f(\alpha,\beta)=\tau$.
Let $\lambda=\lambda_{\Gamma}$.

\begin{lemma}\label{lemma: limit}
No $\alpha,\delta\in \Gamma$ satisfy $\lambda=f(\alpha,\delta)$. In particular, $\lambda$ is a limit ordinal.
\end{lemma}

\begin{proof}
Suppose for the sake of contradiction there exist such $\alpha$ and $\delta$ in $\Gamma$.
Since there are uncountably many $\delta\in \Gamma$ with $\lambda=f(\alpha,\delta)$, we can pick some $\delta\in \Gamma$ large enough such that there is some $\beta\in \Gamma\cap(\alpha,\delta)$ with $D\cap (\alpha, \beta)\neq \emptyset$. But then the fact that $f$ is increasing implies that $f(\beta,\delta)>f(\alpha,\delta)=\lambda$, contradicting the definition of $\lambda$.
\end{proof}

\begin{lemma}\label{lemma: unbounded}
If $N\prec H(\theta)$ is a countable elementary submodel containing relevant objects including $\Gamma, f$ and $D$ where $\theta$ is a large enough regular cardinal, and $\delta\in\Gamma\setminus N$, then 
	\begin{enumerate}
	\item $\{f(\alpha,\delta): \alpha\in \Gamma\cap N\}$ is cofinal in $\lambda_N=_{def} \sup(N\cap\lambda)$, 
	\item $\lambda$ has cofinality $\omega_1$ and 
	\item there is an uncountable set $\Sigma\subseteq \Gamma$ such that $\bigcup\{A_\alpha: \alpha\in \Sigma\}$ has order type $\omega_1$.
	\end{enumerate}

\end{lemma}

\begin{proof}

\begin{enumerate}
\item Suppose for the sake of contradiction that $\{f(\alpha,\delta): \alpha\in \Gamma\cap N\}$ is bounded above by some $\tau\in \lambda\cap N$. The elementarity of $N$ implies that 
for any $\beta\in \Gamma$, there exists $\gamma\in \Gamma - (\beta+1)$ such that for all $\alpha\in \Gamma\cap \beta$, $f(\alpha,\gamma)< \tau$. We can then recursively build $\Gamma'\in [\Gamma]^{\omega_1}$ such that $\lambda_{\Gamma'}\leq \tau < \lambda_N\leq \lambda$. This contradicts the minimality of $\lambda=\lambda_\Gamma$.
\item Suppose for the sake of contradiction that $\cf(\lambda)= \omega$. Then we have $\sup \lambda
\cap N = \lambda$. Fix some $\delta\in \Gamma-N$ such that there are $\beta_0<\beta_1<\delta$ such that for $i<2$, $\beta_i\in \Gamma-N$ and $\emptyset\neq D\cap (\beta_0, \beta_1)$. We can find such ordinals since $\otp(\Gamma)=\omega_1$ and $D$ is unbounded in $\sup \Gamma$.
By (1), we know that $\{f(\alpha,\delta): \alpha\in N\cap\Gamma\}$ is cofinal in $\lambda_N$. In particular, since $f$ is monotone on $\Gamma$, we have $f(\beta_0, \delta)\geq \lambda_N = \lambda$. By the choice of $\beta_1$ and the fact that $f$ is increasing on $\Gamma$, we have that $f(\beta_1,\delta)> f(\beta_0, \delta)\geq \lambda$, which contradicts the definition of $\lambda$.
\item It suffices to show that for any $\alpha\in \Gamma$, there exists $\gamma\in \Gamma$ such that for any $\beta\in \Gamma\cap \alpha$, $A_\beta < A_\gamma$. Pick $\gamma\in \Gamma-N$ and fix any $\tau\in A_\gamma$. By the definition of $A_\gamma$, there is some $\delta\in \Gamma-N$ such that $\tau = f(\gamma,\delta)$. By (1) and the fact that $f$ is monotone on $\Gamma$, $\tau=f(\gamma,\delta)\geq \lambda_N$. For any $\alpha\in N\cap \Gamma$, since $A_\alpha\in N$ is countable, $A_\alpha\subseteq N\cap \lambda$. As a result, $A_\gamma > A_\alpha$ for any $\alpha\in \Gamma\cap N$. The conclusion follows from the elementarity of $N$.
\end{enumerate}
\end{proof}

\section{A 3-dimensional strong coloring on \texorpdfstring{$\omega_2$}{omega-2}: new ingredients}\label{section: omega_2}

The main difference between the proof presented here and that in \cite{Todorcevic3D} is that we remove the need for an arbitrary injection from $\omega_2$ to $2^{\omega_1}$, and use the $C$-sequence to generate a more canonical tree. We include a proof here because this idea will be crucial in Section \ref{section: weaksquare} where new results are obtained. The proof presented here is somewhat detailed since many elements will be recycled in the later sections.

Fix the following data for the coloring.

\begin{enumerate}
\item Fix a $C$-sequence $\langle C_\alpha\subseteq \alpha: \alpha<\omega_2\rangle$ with $\otp(C_\alpha)=\cf(\alpha)$ for all $\alpha<\omega_2$. Write $\Delta(\alpha,\beta)=\Delta(C_\alpha,C_\beta)$ as in Definitions \ref{def: delta} and \ref{def: deltaclub}.
\item Fix a coloring $h_2:[\omega_1]^2\to\omega$ witnessing $\omega_1\not\to[\omega_1]^2_\omega$ as in \cite{TodorcevicWalks87}.
\item Fix $d:\omega_1\to\omega$ with $d^{-1}\{n\}$ stationary for each $n\in\omega$.
\item For each $\beta<\omega_2$, fix an injection $l_\beta:\beta\to\omega_1$.
\item For each $\nu<\omega_1$, fix an injection $e_\nu:\nu\to\omega$. For $m\in\omega$, let $F_m(\nu)=\{\xi<\nu:e_\nu(\xi)\leq m\}$.
\end{enumerate}

Let us now define a function $c: [\omega_2]^3\to \omega$. Fix $\alpha<\beta<\gamma\in \omega_2$.
Let $k=\chi(\alpha,\beta,\gamma)$. Define $\gamma^*=\Tr(\alpha,\gamma)(k)$ and $\beta^*=\Tr(\alpha,\beta)(k)$. We decide the value $c(\alpha,\beta,\gamma)$ based on the following:

\begin{enumerate}

\item If $(C_{\beta^*}-\alpha)\cap (C_{\gamma^*}-\alpha)=\emptyset$, then we output $osc(C_{\beta^*}-\alpha, C_{\gamma^*}-\alpha)-1$. 
\item Otherwise. 
	\begin{enumerate}
	\item If $\lambda_k(\alpha,\gamma)\neq \lambda_k(\alpha,\beta)$ (recall Definition \ref{def: lambdak}), output $$h_2(l_\alpha(\lambda_k(\alpha,\gamma)), l_\alpha(\lambda_k(\alpha,\beta))).$$ \emph{Otherwise put $\sigma=\lambda_k(\alpha,\gamma)=\lambda_k(\alpha,\beta)$ and define $\alpha^*=\Tr(\sigma+1,\alpha)(k)$.}
	\item Else if $\xi=\Delta(C_{\alpha^*},C_{\beta^*})=\Delta(C_{\beta^*},C_{\gamma^*})=\Delta(C_{\alpha^*},C_{\gamma^*})$ and $C_{\alpha^*}(\xi)<C_{\beta^*}(\xi)<C_{\gamma^*}(\xi)$, output \[ h_2(l_{\alpha}(C_{\beta^*}(\xi)), l_{\alpha}(C_{\gamma^*}(\xi))). \]
	
\item Else if $\Delta(\alpha^*, \beta^*)=\min\{\Delta(\beta^*,\gamma^*),\Delta(\alpha^*, \beta^*), \Delta(\alpha^*,\gamma^*)\}<\max \{\Delta(\beta^*,\gamma^*), \Delta(\alpha^*,\gamma^*)\}$ and $\otp(C_{\beta^*}\cap \alpha)< \otp(C_{\gamma^*}\cap \alpha)$, 
	  then output $h_2(\otp(C_{\beta^*}\cap \alpha), \otp(C_{\gamma^*}\cap \alpha))$.
	 	\item Else if $\Delta(\alpha^*,\beta^*)<\Delta(\beta^*,\gamma^*)$ and $\otp(C_{\beta^*}\cap \alpha)\geq \otp(C_{\gamma^*}\cap \alpha)$, output $h_2(\Delta(\alpha^*,\beta^*),\Delta(\beta^*,\gamma^*))$.

	\end{enumerate}

\item Else if $\Delta(\alpha^*,\beta^*)>\Delta(\beta^*, \gamma^*)$, let $\tau=\otp(C_{\gamma^*}\cap \beta)$ and output $$d(\min\bigl(F_{e_\tau(\Delta(\alpha^*, \beta^*))}(\tau)\setminus\otp(C_{\beta^*}\cap \alpha)\bigr)).$$
\end{enumerate}

\begin{theorem}\label{theorem: omega2}
The coloring $c: [\omega_2]^3\to \omega$ defined above satisfies the following: for any $\Gamma\in [\omega_2]^{\omega_1}$ and any $n\in \omega$, there are $\alpha<\beta<\gamma\in \Gamma$ such that $c(\alpha,\beta,\gamma)=n$.
\end{theorem}

\begin{proof}
Fix $\Gamma\in [\omega_2]^{\omega_1}$ and $n\in \omega$. By passing to a subset, we may assume that $\otp(\Gamma)=\omega_1$. Let $\varepsilon=\sup \Gamma\in \cof(\omega_1)\cap \omega_2$.
Let $\vec{M}=\langle M_i: i<\omega_1\rangle$ be an $\in$-increasing sequence of elementary submodels of $H(\theta)$ containing all relevant objects and $D=\{\sup(M_i\cap\varepsilon): i<\omega_1\}$. Here $\theta$ is a large enough regular cardinal. Let $N\prec H((2^{\theta})^+)$ be another elementary submodel containing all relevant objects including $\vec{M}$. Let $\varepsilon_N=\sup(N\cap\varepsilon)$.

\textbf{(1):} If there exists some $\gamma\in \Gamma-N$ such that $\bar{\gamma}=\last{\varepsilon_N}{\gamma}=\Tr(\varepsilon_N, \gamma)(k)$ (for some $k\in \omega$) and it is not the case that a tail of $D\cap N$ is contained in $C_{\bar{\gamma}}$, then the proof of \cite[Lemma 10.3.4]{walksbook} applies to produce $\alpha<\beta \in N-(\lambda(\bar{\gamma},\gamma)+1)$ with $C_{\bar{\gamma}}\cap (\alpha,\beta) \neq \emptyset$ such that letting $\bar{\beta}=\Tr(\alpha,\beta)(k)$, we have that $(C_{\bar{\beta}}-\alpha) \cap (C_{\bar{\gamma}}-\alpha)=\emptyset$ and $osc(C_{\bar{\beta}}-\alpha, C_{\bar{\gamma}}-\alpha)=n+1$. Let us check that Case \textbf{(1)} applies to this triple. Since $\lambda(\bar{\gamma},\gamma)<\alpha<\beta<\bar{\gamma}<\gamma$, we have that $\tr(\alpha,\gamma)=\tr(\bar{\gamma},\gamma)^\frown \tr(\alpha,\bar{\gamma})$ and $\tr(\beta,\gamma)=\tr(\bar{\gamma},\gamma)^\frown \tr(\beta,\bar{\gamma})$. Since $C_{\bar{\gamma}}\cap (\alpha,\beta)\neq\emptyset$, we know that $\bar{\gamma}$ is the least point in the common part of $\tr(\alpha,\gamma)$ and $\tr(\beta,\gamma)$. Therefore $\chi(\alpha,\beta,\gamma)=k$. As a result, $\gamma^*=\Tr(\alpha,\gamma)(k)=\bar{\gamma}$ and $\beta^*=\Tr(\alpha,\beta)(k)=\bar{\beta}$. As $(C_{\bar{\beta}}-\alpha) \cap (C_{\bar{\gamma}}-\alpha)=\emptyset$, Case \textbf{(1)} of the definition applies and we have $c(\alpha,\beta,\gamma)=osc(C_{\bar{\beta}}-\alpha, C_{\bar{\gamma}}-\alpha)-1=n$.

\textbf{(2):} Assume the scenario in \textbf{(1)} does not happen. Fix some $\gamma\in \Gamma-N$ such that $\bar{\gamma}=\last{\varepsilon_N}{\gamma}=\Tr(\varepsilon_N, \gamma)(k)$. Removing some initial segment of $D$ if necessary, we may assume that $D\cap N \subseteq C_{\bar{\gamma}}$. By the elementarity of $N$, we know that there exist $\zeta<\varepsilon$ and a stationary set $S\subseteq \acc(D)$ such that for any $\nu\in 
S$, we can find some $\eta_\nu\in \Gamma$ such that, letting $\bar{\eta}_{\nu}=\Tr(\nu, \eta_{\nu})(k)$, we have that $D\cap \nu\subseteq C_{\bar{\eta}_\nu}$ and $\lambda_k(\nu,\eta_\nu)< \zeta$. For all sufficiently large $\nu\in S$, $\cf(\bar{\eta}_\nu)=\omega_1$.
By shrinking $\Gamma$ to an uncountable subset if necessary, we may assume $\Gamma=\{\eta_{\nu}: \nu \in \acc(S)\cap S-(\zeta+1)\}$ and for any $\nu'<\nu$, $\eta_{\nu'}<\nu\in \acc(D)$ and $\bar{\eta}_{\nu}\in \cof(\omega_1)$.
The following is the main point of the pruning above.

	\begin{claim}
	 For any $\{\eta_{\nu_0}, \eta_{\nu_1}, \eta_{\nu_2}\}\in [\Gamma]^3$, $\chi(\eta_{\nu_0}, \eta_{\nu_1}, \eta_{\nu_2})=k$. Furthermore, $\bar{\eta}_{\nu_i} = \Tr(\eta_{\nu_0}, \eta_{\nu_i})(k)$ for $i=1,2$.
	\end{claim}
	
	\begin{proof}[Proof of the Claim]
	First let us compare $\Tr(\eta_{\nu_0},\eta_{\nu_2})$ with $\Tr(\eta_{\nu_1},\eta_{\nu_2})$. Since $\lambda(\bar{\eta}_{\nu_2}, \eta_{\nu_2})<\zeta<\eta_{\nu_0}<\eta_{\nu_1}<\nu_2\leq \bar{\eta}_{\nu_2}<\eta_{\nu_2}$, by Fact \ref{fact: passing}, we have that for $i<2$, 
 $\tr(\eta_{\nu_i},\eta_{\nu_2})=\tr(\bar{\eta}_{\nu_2},\eta_{\nu_2})^\frown \tr(\eta_{\nu_i},\bar{\eta}_{\nu_2})$. From here we can conclude $\chi(\eta_{\nu_0}, \eta_{\nu_1}, \eta_{\nu_2}) \geq |\tr(\bar{\eta}_{\nu_2},\eta_{\nu_2})|=k$. To see that the value is exactly $k$, note that $C_{\bar{\eta}_{\nu_2}} \cap (\eta_{\nu_0},\eta_{\nu_1})\neq \emptyset$ since $\nu_1\in D\cap \nu_2 \subseteq C_{\bar{\eta}_{\nu_2}}$ and $\nu_1\in (\eta_{\nu_0},\eta_{\nu_1})$. 
		
	To see the ``furthermore'' part, first note that the argument above gives $\bar{\eta}_{\nu_2} = \Tr(\eta_{\nu_0}, \eta_{\nu_2})(k)$. Since $\lambda(\bar{\eta}_{\nu_1}, \eta_{\nu_1})\leq \zeta<\eta_{\nu_0}<\nu_1 \leq \bar{\eta}_{\nu_1}<\eta_{\nu_1}$, by Fact \ref{fact: passing}, we know that $\tr(\eta_{\nu_0}, \eta_{\nu_1})=\tr(\bar{\eta}_{\nu_1},\eta_{\nu_1})^\frown \tr(\eta_{\nu_0}, \bar{\eta}_{\nu_1})$ and $\tr(\nu_1, \eta_{\nu_1})=\tr(\bar{\eta}_{\nu_1}, \eta_{\nu_1})^\frown \tr(\nu_1, \bar{\eta}_{\nu_1})$. Since $\Tr(\nu_1,\eta_{\nu_1})(k)=\bar{\eta}_{\nu_1}$, we know that $|\tr(\bar{\eta}_{\nu_1}, \eta_{\nu_1})|=k$. It follows that $\Tr(\eta_{\nu_0}, \eta_{\nu_1})(k)=\bar{\eta}_{\nu_1}$.\end{proof}

	\begin{observation}\label{observation: not(1)}
	For any $\{\eta_{\nu_0}, \eta_{\nu_1}, \eta_{\nu_2}\}\in [\Gamma]^3$, $(C_{\bar{\eta}_{\nu_2}}- \eta_{\nu_0}) \cap (C_{\bar{\eta}_{\nu_1}}- \eta_{\nu_0}) \neq \emptyset$.
	\end{observation}
	\begin{proof}[Proof of the Observation]
	By our choice of $\Gamma$ and $D$, we have that $\nu_1\in \acc(D)$ and $D\cap \nu_1\subseteq C_{\bar{\eta}_{\nu_1}} \cap C_{\bar{\eta}_{\nu_2}}$. In particular, any $\xi\in D\cap \nu_1 - \eta_{\nu_0}$ is in $(C_{\bar{\eta}_{\nu_2}}- \eta_{\nu_0}) \cap (C_{\bar{\eta}_{\nu_1}}- \eta_{\nu_0})$.
	\end{proof}

As a result of Observation \ref{observation: not(1)}, for any $\bar{a}\in [\Gamma]^3$, the hypothesis of Case \textbf{(1)} does not hold.

\textbf{(2a):} If $A=\{\lambda(\bar{\eta}_\nu,\eta_\nu): \eta_\nu\in \Gamma\}$ is uncountable, letting $\alpha = \min \Gamma$, we know that $A$ is an uncountable subset of $\alpha$. Observe that for any $\eta_{\nu'}<\eta_{\nu}\in \Gamma$, $\lambda(\bar{\eta}_{\nu},\eta_\nu) = \lambda_k(\eta_{\nu'},\eta_\nu)$ since $\lambda(\bar{\eta}_{\nu},\eta_\nu) < \eta_{\nu'}<\bar{\eta}_\nu < \eta_\nu$ and Fact \ref{fact: passing} applies. As a result, $l_\alpha '' A$ is an uncountable subset of $\omega_1$. By the property of $h_2$, we can find $\eta_{\nu_1}<\eta_{\nu_2}\in \Gamma$ such that $$\lambda_k(\alpha,\eta_{\nu_1})=\lambda(\bar{\eta}_{\nu_1}, \eta_{\nu_1})\neq \lambda(\bar{\eta}_{\nu_2}, \eta_{\nu_2}) = \lambda_k(\alpha,\eta_{\nu_2})$$ and $h_2(l_\alpha(\lambda_k(\alpha,\eta_{\nu_1})), l_\alpha(\lambda_k(\alpha,\eta_{\nu_2})))=n$. By Case \textbf{(2a)}, we have that $c(\alpha,\eta_{\nu_1},\eta_{\nu_2})=n$.

\emph{Let us simplify some notation for the rest of the proof.} For each $\delta\in \Gamma$ with $\delta=\eta_{\nu}$ for some $\nu$, let $\delta^* = \bar{\eta}_\nu$. Note that this operation is definable from $S, \Gamma$ and the $C$-sequence. 
 
\textbf{(2b)} Otherwise, $A$ as defined above is countable. In such a case, by shrinking $\Gamma$ to an uncountable subset if necessary, we may assume that $A$ is a singleton, say $\{\sigma\}$. In such a case, we note that for each $\alpha\in \Gamma$, $\alpha^*=\Tr(\sigma+1, \alpha)(k)$. The reason is that we can apply Fact \ref{fact: passing} in the context where $|\tr(\alpha^*, \alpha)|=k$ and $\lambda(\alpha^*, \alpha)<\sigma+1 < \alpha^*<\alpha$.

Next, we consider the tree $T$ defined as follows: a function $t\in T$ iff there are uncountably many $\gamma\in \Gamma$ such that $t\sqsubseteq C_{\gamma^*}$, where we identify $C_{\gamma^*}$ with the function on $\otp(C_{\gamma^*})$ that enumerates $C_{\gamma^*}$ increasingly. The order on $T$ is end extension.

Suppose there is some $\xi\in \omega_1$ and an uncountable $\Gamma'\subseteq \Gamma$ such that for any $\beta<\gamma\in \Gamma'$, $C_{\beta^*}\restriction \xi = C_{\gamma^*}\restriction \xi$ and $C_{\beta^*}(\xi)<C_{\gamma^*}(\xi)$. Note that $\{C_{\beta^*}(\xi): \beta\in \Gamma'\}$ cannot be separated by points from $D$. The reason is that if there is some $\eta\in D\cap (C_{\beta^*}(\xi),C_{\gamma^*}(\xi))$, take $\delta\in \Gamma'$ such that $C_{\delta^*}(\xi)>C_{\gamma^*}(\xi)$. In particular, $\eta<C_{\delta^*}(\xi)$. The fact that $D\cap \gamma^*\subseteq C_{\delta^*}$ implies that $\eta\in C_{\delta^*}$ above $C_{\delta^*}\restriction \xi$, contradicting the definition of $C_{\delta^*}(\xi)$. In particular, we can find some $\alpha\in \Gamma'$ large enough such that for all $\gamma\in \Gamma'$, $C_{\gamma^*}(\xi)<\alpha$. Since $l_{\alpha}'' \{C_{\delta^*}(\xi): \delta\in \Gamma'\}\in [\omega_1]^{\aleph_1}$, we can find $\beta<\gamma\in \Gamma'- (\alpha+1)$ such that $h_2(l_\alpha(C_{\beta^*}(\xi)), l_\alpha(C_{\gamma^*}(\xi)))=n$. Note that 
	\begin{itemize}
	\item $\lambda_k(\alpha,\gamma)=\lambda_k(\alpha,\beta)$,
	\item  $\Delta(\alpha^*,\beta^*)=\Delta(\beta^*,\gamma^*)=\Delta(\alpha^*, \gamma^*)=\xi$ and
	\item  $C_{\alpha^*}(\xi)<C_{\beta^*}(\xi)<C_{\gamma^*}(\xi)$.
	\end{itemize}
Therefore, Case \textbf{(2b)} applies and we have that $c(\alpha,\beta,\gamma)=n$.

\textbf{(2c):} Suppose that for every $\xi\in \omega_1$ and every $t\in T_\xi$, $\{C_{\beta^*}(\xi): t=C_{\beta^*}\restriction \xi, \beta\in \Gamma\}$ is a countable set. In particular, $T$ is countably splitting. Suppose that there is some $\alpha\in \Gamma$ such that there is $\Gamma'\in [\Gamma]^{\omega_1}$ where $\beta<\gamma\in \Gamma'$ implies that $\otp(C_{\beta^*}\cap \alpha)<\otp(C_{\gamma^*}\cap \alpha)$.
	\begin{itemize}
	\item If $\{\Delta(\alpha^*, \gamma^*): \gamma\in \Gamma'\}$ is countable, then we can find $\beta<\gamma\in \Gamma'$ such that $h_2(\otp(C_{\beta^*}\cap \alpha), \otp(C_{\gamma^*}\cap \alpha))=n$ with $\Delta(\alpha^*,\beta^*)=\Delta(\alpha^*,\gamma^*)<\Delta(\beta^*, \gamma^*)$. 
	\item  If $\{\Delta(\alpha^*, \gamma^*): \gamma\in \Gamma'\}$ is uncountable, we can find $\beta<\gamma\in \Gamma'$ such that $h_2(\otp(C_{\beta^*}\cap \alpha), \otp(C_{\gamma^*}\cap \alpha))=n$ with $\Delta(\alpha^*,\beta^*)=\Delta(\beta^*,\gamma^*)<\Delta(\alpha^*,\gamma^*)$. 
	\end{itemize}
In either case, Case \textbf{(2c)} of the definition applies and we get $c(\alpha,\beta,\gamma)=n$.

\textbf{(2d):} Suppose that for all $\alpha\in \Gamma$, $A_\alpha=\{\otp(C_{\gamma^*}\cap \alpha): \gamma\in\Gamma,\ \alpha<\gamma\}$ is countable. Suppose additionally that $T$ contains an uncountable branch. Let $b$ be such an uncountable branch. We can find an uncountable $E\subseteq \Gamma$ such that for any $\beta<\gamma\in E$, $\Delta(b, C_{\beta^*})<\Delta(b, C_{\gamma^*})$. Shrinking $E$ to an uncountable subset if necessary, we may also assume that if $\alpha\in E$ and $\epsilon\in A_\alpha^E=_{\mathrm{def}} \{\otp(C_{\gamma^*}\cap \alpha): \gamma\in E,\ \alpha<\gamma\}$, then there are uncountably many $\gamma\in E$ such that $\epsilon=\otp(C_{\gamma^*}\cap \alpha)$.
Let $B=\{\Delta(b, C_{\delta^*}): \delta \in E\}\in [\omega_1]^{\omega_1}$.
Find $\xi<\zeta\in B$ such that $h_2(\xi,\zeta)=n$. Let $\alpha<\beta\in E$ such that $\Delta(b,C_{\alpha^*})=\xi$ and $\Delta(b,C_{\beta^*})=\zeta$. By the property of $A_\alpha^E$, we can find some $\gamma\in E-(\beta+1)$ such that $\otp(C_{\beta^*}\cap \alpha)=\otp(C_{\gamma^*}\cap \alpha)$. Therefore, we see that $\Delta(\alpha^*,\beta^*)<\Delta(\beta^*,\gamma^*)$ and $\otp(C_{\beta^*}\cap \alpha)\geq \otp(C_{\gamma^*}\cap \alpha)$. Case \textbf{(2d)} of the definition applies and we get $c(\alpha,\beta,\gamma)=h_2(\xi,\zeta)=n$.

\textbf{(3):} Finally, we deal with the remaining case. By shrinking $\Gamma$ to an uncountable subset if necessary, we may assume the tree $T$, defined in \textbf{(2b)}, is countably splitting and either contains a Cantor subtree or is an Aronszajn tree. To be more specific about $T$ containing a Cantor subtree, we mean there is some least $\nu<\omega_1$ such that $U_\nu=\{C_{\gamma^*}\restriction \nu: \gamma\in \Gamma\}$ is uncountable. We say that \emph{$T$ contains a Cantor subtree of height $\nu$}.
Note that such a $\nu$ is always a limit ordinal. Observe that in such a case, $T\restriction \nu$ is \emph{perfect}. Namely, for every $t\in T\restriction \nu$ and $\xi<\nu$, there exist $t_0,t_1\in T\restriction \nu$ of height $>\xi$ extending $t$ such that $t_0\perp t_1$.

Fix a countable elementary submodel $N\prec H(\theta)$, where $\theta$ is a large enough regular cardinal, containing all relevant objects such that $d(N\cap \omega_1)=n$. Let $$\bar{\nu}=\begin{cases}
\nu, & T\text{ contains a Cantor subtree of height }\nu \\
N\cap \omega_1, & \text{otherwise}
\end{cases}$$

Fix some $\beta<\gamma\in \Gamma$ with $\beta^*,\gamma^*\not\in N$ and $t_0,t_1\in T(\bar{\nu})-N$ such that $t_0\sqsubseteq C_{\beta^*}$, $t_1\sqsubseteq C_{\gamma^*}$ with $t_0\perp t_1$ and $\tau=\otp(C_{\gamma^*}\cap \beta) > N\cap \omega_1$. In particular, $\Delta(\beta^*,\gamma^*)\in N$. Let $m_0=e_\tau(N\cap \omega_1)$. 

There is some $t_0'\in N\cap T$ such that $\xi=\Delta(t_0', C_{\beta^*})>\Delta(t_0,t_1)$ with $e_\tau(\xi)=m\geq m_0$. To see this, 

	\begin{itemize}
	\item if $T$ contains a Cantor subtree of height $\nu$, then for any $\xi<\nu$, $C_{\beta^*}\restriction \xi\in N$, hence as $C_{\beta^*}\restriction \nu \in T- N$, there are unboundedly many $\xi<\nu$ for which there is $s\in T\restriction \nu$ with $\Delta(s, C_{\beta^*})=\xi$. As a result, such $t_0'$ can be found;
	\item if $T$ is an Aronszajn tree, then either for any $\eta\in N\cap \omega_1$, $C_{\beta^*}\restriction \eta\in T$ (hence in $N$) or there is a least $\eta \in N\cap \omega_1$ such that $C_{\beta^*}\restriction \eta \not \in T$. In the former situation, we can utilize the fact that $T$ does not have an uncountable branch to find such a desired $t_0'$. In the latter situation, we know that $C_{\beta^*}\restriction \eta \not \in N$ by the fact that $\beta^*\not \in N$ and the elementarity of $N$. Notice that $\eta$ has to be a limit ordinal since otherwise the situation would have been dealt with in \textbf{(2b)}. An argument similar to that in the previous case shows that if such $t_0'$ does not exist, then $C_{\beta^*}\restriction \eta$ is a member of $N$, which contradicts our assumption.
\end{itemize}

Let $G=F_{m}(\tau)\cap N \in [N\cap \omega_1]^{<\omega}$ and note that $N\cap\omega_1\in F_m(\tau)$ since $m\geq m_0$. Since $t_0'\in T$, the collection of $\alpha\in \Gamma$ such that $t_0'\sqsubseteq C_{\alpha^*}$ is uncountable. In particular, we can find some $\alpha\in N\cap \Gamma$ such that $t_0'\sqsubseteq C_{\alpha^*}$ and $\otp(C_{\beta^*}\cap \alpha)>\sup G$. The reason is that by the triviality assumption, $D \cap N\subseteq C_{\beta^*}$ and $\otp(D \cap N)=N\cap \omega_1$.
Observe that $\otp(C_{\beta^*}\cap \alpha)\in N\cap \omega_1$ since $\alpha\in N$ and thus the countable set $A_\alpha\subseteq N$. As a result, $\min\bigl(F_m(\tau)\setminus\otp(C_{\beta^*}\cap \alpha)\bigr)=N\cap \omega_1$. Since $\Delta(\alpha^*,\beta^*)>\Delta(\beta^*, \gamma^*)$, Case \textbf{(3)} of the definition applies and $c(\alpha,\beta,\gamma)=d(\min\bigl(F_m(\tau)\setminus\otp(C_{\beta^*}\cap \alpha)\bigr))=d(N\cap \omega_1)=n$.
\end{proof}

\section{Stepping up with a special Aronszajn tree}\label{section: weaksquare}

In this section, we prove the following theorem.
	\begin{theorem}\label{theorem: special}
	If there exists an $\omega_2$-special $\omega_3$-Aronszajn tree, then $\omega_3\not\to [\omega_1]^4_\omega$.
	\end{theorem}
By a theorem of Jensen (see \cite[Section 5]{CummingsNotes}), the existence of an $\omega_2$-special $\omega_3$-Aronszajn tree is equivalent to the existence of a $\square_{\omega_2}^*$-sequence. Fix the following data for the coloring.
\begin{enumerate}[label=\textbf{D\arabic*}]
\item Fix a $\square_{\omega_2}^*$-sequence $\langle C_\alpha:\alpha<\omega_3\rangle$. Each $C_\alpha$ is a club in $\alpha$ of order type $\cf(\alpha)$, and for each $\alpha<\omega_3$ the set $\{C_\beta\cap\alpha:\beta<\omega_3\}$ has size at most $\omega_2$.
For each $\alpha<\omega_3$, fix an enumeration $\langle C_{\alpha,i}:i<\omega_2\rangle$ of this set. Let $\operatorname{index}_\alpha(C_\beta)$, the \emph{$\alpha$-index} of $C_\beta\cap\alpha$, be the least $i$ such that $C_{\alpha,i}=C_\beta\cap\alpha$.
Here $\Delta(a,b)=\Delta(C_a,C_b)$ denotes the first disagreement of the increasing enumerations, as in Definition~\ref{def: deltaclub}.
%
\item \label{d1} Fix Todor\v{c}evi\'c's coloring $h_3:[\omega_2]^3\to\omega$ witnessing $\omega_2\not\to[\omega_1]^3_\omega$ \cite{Todorcevic3D}, and Moore's coloring $g:[\omega_1]^2\to\omega$ witnessing $\omega_1\not\to[\omega_1,\omega_1]^2_\omega$ \cite{MooreLspace}. Thus for any $A,B\in[\omega_1]^{\omega_1}$, $g''\{(\alpha,\beta):\alpha\in A,\ \beta\in B,\ \alpha<\beta\}=\omega$.
\item \label{d2} Fix $d:\omega_1\to\omega$ such that $d^{-1}\{n\}$ is stationary for each $n\in\omega$.
\item \label{d3} For each $\beta<\omega_3$, fix an injection $l_\beta:\beta\to\omega_2$. For each $\nu<\omega_1$, fix an injection $e_\nu:\nu\to\omega$. For $m\in\omega$, let $F_m(\nu)=\{\xi<\nu:e_\nu(\xi)\leq m\}$.
\item \label{d4} Use the $C$-sequence below $\omega_2$ for the function $\rho$ of Definition~\ref{def:walk-rank}.
\end{enumerate}

%

Let us now define a function $c: [\omega_3]^4\to \omega$. Fix $\alpha<\beta<\gamma<\delta\in \omega_3$ and let $r=\{\alpha,\beta,\gamma,\delta\}$.
Let $k=\chi(\beta,\gamma,\delta)$. Let $\delta^* = \Tr(\alpha,\delta)(k)$, $\gamma^*=\Tr(\alpha,\gamma)(k)$ and $\beta^*=\Tr(\alpha,\beta)(k)$.
We decide the value $c(\alpha,\beta,\gamma,\delta)$ based on the following cases (if the expression does not make sense at any point, just output $0$):
\begin{enumerate}
\item If $(C_{\Tr(\beta,\gamma)(k)}-\beta) \cap (C_{\Tr(\beta,\delta)(k)}-\beta) = \emptyset$, then output $$osc(C_{\Tr(\beta,\gamma)(k)}-\beta,C_{\Tr(\beta,\delta)(k)}-\beta)-1.$$ 
\item Otherwise. 
	\begin{enumerate}
	\item If $\{\lambda_k(\alpha,\gamma), \lambda_k(\alpha,\beta), \lambda_k(\alpha,\delta)\}$ is not a singleton, output $$h_3(l_\alpha(\lambda_k(\alpha,\gamma)), l_\alpha(\lambda_k(\alpha,\beta)), l_\alpha(\lambda_k(\alpha,\delta))).$$ Otherwise we define $\alpha^*=\Tr(\xi+1, \alpha)(k)$ where $\xi=\lambda_k(\alpha,\gamma)$.
	\item Else if $\{\Delta(a^*, b^*): \{a,b\}\in [r]^2\}=\{\xi\}$ and $C_{\alpha^*}(\xi)<C_{\beta^*}(\xi)<C_{\gamma^*}(\xi)<C_{\delta^*}(\xi)$, then we output $$h_3(l_{\alpha}(C_{\beta^*}(\xi)), l_{\alpha}(C_{\gamma^*}(\xi)),l_{\alpha}(C_{\delta^*}(\xi))).$$
	
\item Else if $\Delta(\alpha^*, \beta^*)=\min\{\Delta(\beta^*,\gamma^*),\Delta(\alpha^*, \beta^*), \Delta(\alpha^*,\gamma^*)\}<\Delta(\gamma^*, \delta^*)$, and $\operatorname{index}_\alpha(C_{\beta^*})< \operatorname{index}_\alpha(C_{\gamma^*})<\operatorname{index}_\alpha(C_{\delta^*})$, 
	  then output $$h_3(\operatorname{index}_\alpha(C_{\beta^*}), \operatorname{index}_\alpha(C_{\gamma^*}),\operatorname{index}_\alpha(C_{\delta^*})).$$
	 	\item Else if the three successive disagreements satisfy \[\Delta(\alpha^*,\beta^*)<\Delta(\beta^*,\gamma^*)<\Delta(\gamma^*,\delta^*)\] and $\operatorname{index}_\alpha(C_{\gamma^*})\geq \operatorname{index}_\alpha(C_{\delta^*})$, output \[h_3(\Delta(\alpha^*,\beta^*),\Delta(\beta^*,\gamma^*),\Delta(\gamma^*,\delta^*)).\]

	\end{enumerate}
	
\item Otherwise. Let $\tau=\otp(C_{\delta^*}\cap \gamma)$, $\tau'=\otp(C_{\beta^*}\cap \alpha)$ and $$w=\chi(\tau',\otp(C_{\delta^*}\cap \beta),\tau)$$ and let $\bar{\tau}=\Tr(\tau' ,\tau)(w)$. 
Let $\bar{\tau}'=\Tr(\Delta(\alpha^*, \gamma^*),\tau')(w)$.
	\begin{enumerate}
	\item If $\Delta(\alpha^*,\gamma^*)> \Delta(\alpha^*, \beta^*) =\Delta(\gamma^*, \delta^*)$, then output $$osc(C_{\bar{\tau}}-\Delta(\alpha^*, \gamma^*), C_{\bar{\tau}'}-\Delta(\alpha^*, \gamma^*))-1$$ if the value is non-negative, and output $0$ otherwise.
	\item Else if $\Delta(\gamma^*,\delta^*)< \Delta(\alpha^*,\beta^*)<\Delta(\beta^*,\delta^*)$, let $\xi=\Delta(\beta^*,\delta^*)$ and output $g(\rho(\xi, \tau'), \rho(\xi, \tau))$.
	\item Else if $\Delta(\alpha^*,\beta^*)> \max\{\Delta(\beta^*,\gamma^*),\Delta(\gamma^*,\delta^*)\}$, let $\nu=\otp(C_{\bar{\tau}}\cap \otp(C_{\delta^*}\cap \beta))$, $\xi=\Delta(\alpha^*, \beta^*)$, $m=e_\nu(\otp(C_{\bar{\tau}} \cap \xi))$ and output $d(\min(F_m(\nu)\setminus\otp(C_{\bar{\tau}}\cap \otp(C_{\delta^*}\cap \alpha))))$.
	\end{enumerate}
	
\end{enumerate}

\begin{figure}
    \centering 
   \includegraphics[scale=0.6]{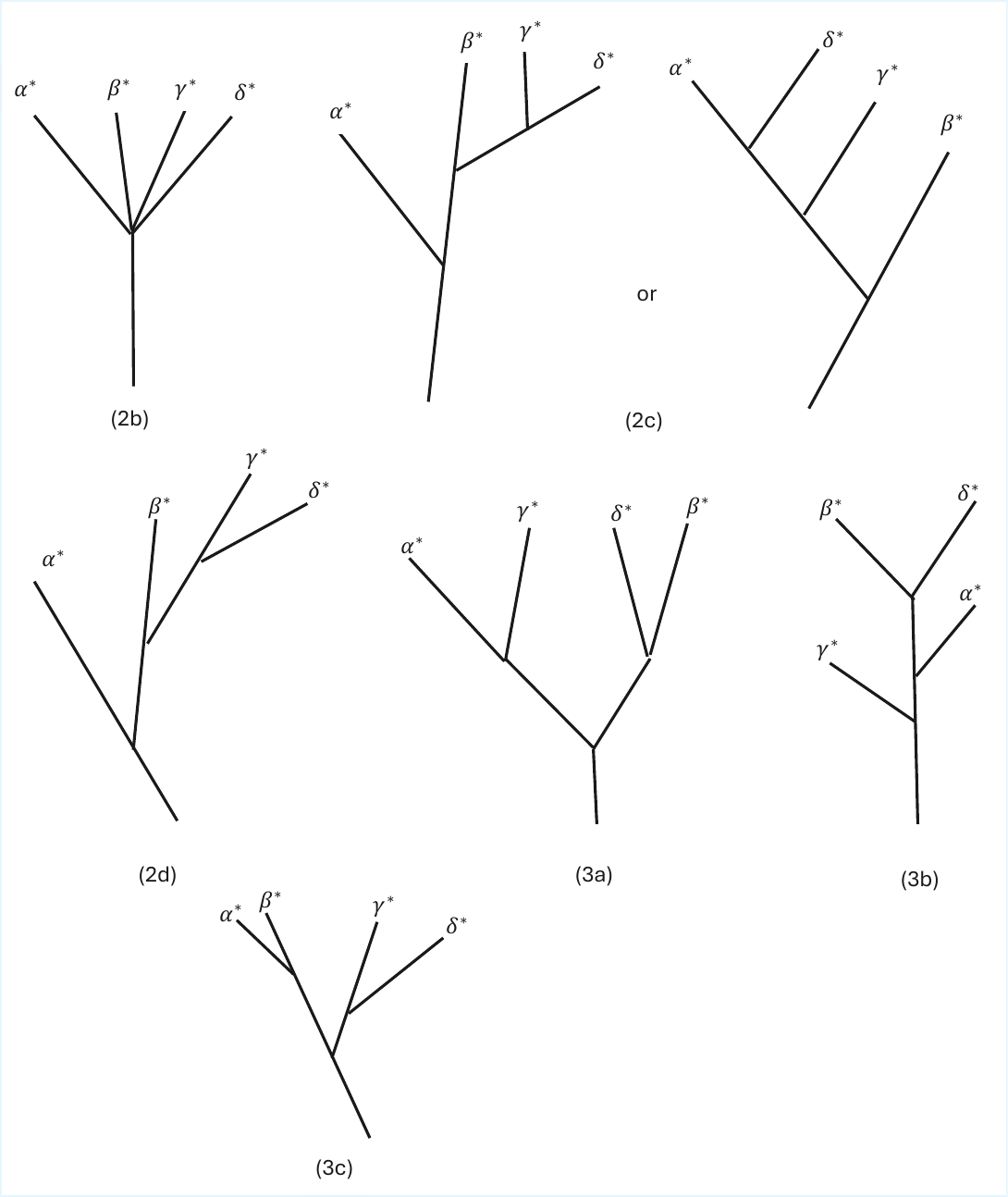}
    \caption{Tuple configurations based on a special Aronszajn tree} 
    \label{fig:special} 
\end{figure}

\begin{theorem}\label{theoerem: strongSpecial}
The coloring $c: [\omega_3]^4\to \omega$ defined above satisfies the following: for any $\Gamma\in [\omega_3]^{\omega_1}$ and any $n\in \omega$, there are $\alpha<\beta<\gamma<\delta\in \Gamma$ such that $c(\alpha,\beta,\gamma,\delta)=n$.
\end{theorem}

A proof of Theorem \ref{theoerem: strongSpecial} will occupy the rest of this section.

For cases \textbf{(1)} and \textbf{(2)} we follow the proof of Theorem~\ref{theorem: omega2}. We record the notation and the changes needed for \textbf{(2c)} and \textbf{(2d)}.

Let $\Gamma\in [\omega_3]^{\omega_1}$ and $n\in \omega$ be given. Let $\varepsilon=\sup \Gamma\in \cof(\omega_1)\cap \omega_3$.
Let $\vec{M}=\langle M_i: i<\omega_1\rangle$ be an $\in$-increasing sequence of elementary submodels of $H(\theta)$ containing all relevant objects and $D=\{\sup M_i\cap \varepsilon: i<\omega_1\}$. Here $\theta$ is a large enough regular cardinal. Let $N\prec H((2^{\theta})^+)$ be another elementary submodel containing all relevant objects including $\vec{M}$. Let $\varepsilon_N=\sup N\cap \varepsilon$.

\textbf{(1):} If there exists some $\delta\in \Gamma-N$ such that $\bar{\delta}=\last{\varepsilon_N}{\delta}=\Tr(\varepsilon_N, \delta)(k)$ (for some $k\in \omega$) and it is not the case that a tail of $D\cap N$ is contained in $C_{\bar{\delta}}$, then the argument for case \textbf{(1)} in the proof of Theorem \ref{theorem: omega2}, applied to $\Gamma\setminus\{\min\Gamma\}$, gives $\beta<\gamma<\delta$ for which the two tails in \textbf{(1)} are disjoint and have oscillation $n+1$ in the indicated order. Taking $\alpha=\min\Gamma$ gives $c(\alpha,\beta,\gamma,\delta)=n$.

Suppose the hypothesis in \textbf{(1)} does not hold. Then we prune exactly as in \textbf{(2)} in the proof of Theorem \ref{theorem: omega2}. For the sake of clarity, let us record the result of the pruning. 

	\begin{enumerate}[label=\textbf{B\arabic*}]
	\item \label{b1} There is $k\in \omega$ such that for any $\{\alpha,\beta,\gamma\}\in [\Gamma]^3$, $\chi(\alpha,\beta,\gamma)=k$.
	\item \label{b2} For any $\beta\in \Gamma$, there exists $\beta^*$ such that for any $\alpha\in \Gamma\cap \beta$, $\beta^*=\Tr(\alpha,\beta)(k)$.
	\item \label{b3} For any $\alpha<\beta\in \Gamma$, there is some $\nu\in \acc(D)\cap \acc(C_{\beta^*}) \cap \beta^*$ such that $\alpha<\nu$.
	\item \label{b4} For any $\beta<\gamma<\delta\in \Gamma$, $(C_{\gamma^*}-\beta)\cap (C_{\delta^*}-\beta)\neq \emptyset$.
	\end{enumerate}
	
For \textbf{(2a)}--\textbf{(2d)} we follow the corresponding arguments in the proof of Theorem~\ref{theorem: omega2}. In \textbf{(2c)} and \textbf{(2d)} we compare the $\alpha$-indices of the clubs in place of their order types below $\alpha$. This is where we use the $\square^*_{\omega_2}$-sequence.

From this point on, we may assume additionally the following: 

	\begin{enumerate}[label=\textbf{A\arabic*}]
	\item \label{a1} there is $\eta<\min \Gamma$ such that for all $\alpha<\beta \in \Gamma$, $\lambda_k(\alpha,\beta)=\eta$.
	\item \label{a2} for any $t$, the collection $\{C_{\alpha^*}(\dom(t)): t\sqsubseteq C_{\alpha^*}, \alpha\in \Gamma \}$ is countable.
	\item \label{a3} for each $\alpha\in \Gamma$, the collection $I_\alpha = \{\operatorname{index}_\alpha(C_{\delta^*}): \delta\in \Gamma\}$ is countable.
	\end{enumerate}
Consider the function $f$ mapping $(\alpha,\beta)\in [\omega_3]^2$ to $\otp(C_{\Tr(\alpha,\beta)(k)}\cap \alpha)\in \omega_2$. In particular, we see that $f$ is monotone, locally small and increasing on $\Gamma$, since for any $\alpha<\beta\in \Gamma$, $f(\alpha,\beta)=\otp(C_{\beta^*}\cap \alpha)$.
For each $\alpha\in\Gamma$, let $A_\alpha=\{f(\alpha,\delta):\delta\in\Gamma,\ \alpha<\delta\}$ be the row of order types. Put $W=\bigcup_{\alpha\in\Gamma} A_\alpha$ and $\lambda=\sup W$. By Lemmas~\ref{lemma: least}, \ref{lemma: limit}, and \ref{lemma: unbounded}, we may assume
	\begin{enumerate}[label=\textbf{P\arabic*}]
	\item \label{p1} for any $\Sigma\in [\Gamma]^{\omega_1}$, $\sup\{\otp(C_{\delta^*}\cap \alpha)\in \omega_2: \alpha<\delta\in \Sigma\}=\lambda$, 
	\item \label{p2} $\otp(W)=\omega_1=\cf(\lambda)$,
	\item \label{p3} for each $\alpha\in \Gamma$ and $\tau\in A_\alpha$, there are uncountably many $\delta\in \Gamma$ such that $\tau=\otp(C_{\delta^*}\cap \alpha)$.
	\item \label{p4} for any countable elementary submodel $N$ of $H((2^\theta)^+)$ containing relevant objects, if $\delta\in \Gamma$ satisfies that $\delta^* \not\in N$, $\{\otp(C_{\delta^*}\cap \alpha): \alpha\in \Gamma\cap N\}$ is unbounded in $\sup N\cap \lambda$.
	\end{enumerate}

Consider the tree $T$ consisting of elements $t$ such that there are uncountably many $\gamma\in \Gamma$ such that $t\sqsubseteq C_{\gamma^*}$. The order on $T$ is end extension. We may assume the following, since otherwise \textbf{(2b)} or \textbf{(2d)} applies, respectively. 
    \begin{enumerate}[label=\textbf{T\arabic*}]
	\item \label{T1} $T$ is countably splitting and 
	\item \label{T2} $T$ does not have an uncountable branch cofinal in $\lambda$.
	\end{enumerate}

Let us show that this tree $T$ is ``\emph{slim}''.

	\begin{claim}\label{claim: slim}
		\begin{itemize}
		\item For any $t\in T$, $\dom(t)<\lambda$.
		\item $T$ is normal. For any $t\in T$ and $\zeta<\lambda$, there exists an extension $t'$ of $t$ in $T$ with $\dom(t')\geq \zeta$. In particular,  $\{\dom(t): t\in T\}$ is cofinal in $\lambda$.
		\item For any $\xi<\lambda$, $\{t\in T: \dom(t)\in W\cap  (\xi+1)\}$ is countable.
		\end{itemize}
	\end{claim}
	
	\begin{proof}[Proof of the Claim]
	\begin{itemize}
	\item Let $t\in T$. By definition, there are uncountably many $\delta\in \Gamma$ such that $t\sqsubseteq C_{\delta^*}$. Let $\xi=\dom(t)$. Pick $\delta_0<\delta_1\in \Gamma$ with $t\sqsubseteq C_{\delta^*_0}, C_{\delta^*_1}$. Since $\mathrm{ran}(t)\subseteq \delta_0$, we have $\xi\leq\otp(C_{\delta_1^*}\cap\delta_0)<\lambda$ by Lemma~\ref{lemma: limit}. 
	
	\item Fix $t\in T$ and $\zeta<\lambda$. We may assume that $\zeta>\dom(t)$. We know that $B=\{\delta\in \Gamma: t\sqsubseteq C_{\delta^*}\}$ is uncountable by the definition of $T$. Let $N\prec H(\theta)$ be a countable elementary submodel containing all relevant objects including $\Gamma, t, \zeta$. Let $\delta\in B$ be such that $\delta^* \not \in N$. By \ref{p4}, $\{\otp(C_{\delta^*}\cap \alpha): \alpha\in N\cap \Gamma\}$ is cofinal in $N\cap \lambda$. In particular, we can find some $\alpha\in N\cap \Gamma$ with $\otp(C_{\delta^*}\cap \alpha)=\tau \geq \zeta$. Note that $\tau\in  N$. By the elementarity of $N$, there are uncountably many $\delta\in \Gamma$ such that $\otp(C_{\delta^*}\cap \alpha)=\tau$ and $t\sqsubseteq C_{\delta^*}$. Let us call this uncountable set $E\subseteq \Gamma$. Since $I_\alpha$ is countable, we can find some $i\in \omega_2\cap N$ such that there are uncountably many $\delta\in E$ such that $C_{\delta^*}\cap \alpha = C_{\alpha,i}$. Let $s: \otp(C_{\alpha,i})\to C_{\alpha,i}$ be the inverse of the transitive collapse of $C_{\alpha,i}$. It is now easy to see that $s\in T$, $t\sqsubseteq s$, and $\dom(s)=\tau\geq \zeta$.
 	\item It clearly suffices to show that for any $\xi < \lambda$, $\{t\in T: \dom(t)= \xi\}$ is countable, by the normality of $T$. Suppose for the sake of contradiction that there is some $\xi<\lambda$ such that $F=\{t\in T: \dom(t)= \xi\}$ is uncountable. Take some countable elementary submodel $N\prec H(\theta)$ containing all relevant objects including $\xi$ and $\Gamma$. Pick some $t\in F - N$ and some $\delta\in \Gamma$ such that $\delta^*\not\in N$ with $t\sqsubseteq C_{\delta^*}$. By \ref{p4}, we can find $\alpha\in \Gamma\cap N$ such that $\otp(C_{\delta^*}\cap \alpha)>\xi$. By \ref{a3}, $C_{\delta^*}\cap \alpha\in N$. But then $t=C_{\delta^*} \restriction \xi=(C_{\delta^*}\cap \alpha)\restriction \xi \in N$, which is a contradiction.
	\end{itemize}
	\end{proof}

As a result, $\{\Delta(\gamma^*,\delta^*): \gamma,\delta\in \Gamma\}\cap \lambda$ is cofinal in $\lambda$. We use this in the rest of the proof.

Let $\chi$ be a large enough regular cardinal and fix a continuous $\in$-increasing sequence $\vec{M}=\langle M_i:i<\omega_1\rangle$ of countable elementary submodels of $H(\chi)$ containing all relevant objects. Let $E=\{\sup(M_i\cap\lambda):i<\omega_1\}$.

\subsection{The downstairs oscillation case (3a)}\label{osc} Suppose there is a countable elementary submodel $N\prec H((2^\chi)^+)$ containing all relevant objects including $\vec{M}$, with $\lambda_N=\sup(N\cap\lambda)$, and 
there is some $\tau\in \lambda-\lambda_N$ such that there are $\gamma, \delta\in \Gamma$ with $\tau = \otp(C_{\delta^*}\cap \gamma)$ and $\Delta(\gamma^*, \delta^*)\in N\cap \lambda$, and the following holds: 
$E\cap N\cap \lambda \not \subseteq^* \acc(C_{\last{\lambda_N}{\tau}})\cap \lambda_N$. Recall ``$\subseteq^*$'' here means containment modulo a bounded subset of $\lambda_N$.

Let $\bar{\tau}=\last{\lambda_N}{\tau}$ and $w\in \omega$ be such that $\bar{\tau}=\Tr(\lambda_N, \tau)(w)$. Let $\xi= \lambda(\bar{\tau}, \tau)<\lambda_N$. For the interval construction, fix a positive integer $n$. Take an increasing sequence $\langle\delta_i:i<n\rangle$ in $E\cap N\cap\lambda_N$ and put $\delta_n=\lambda$ so that for each $i<n$,
	\begin{itemize}
	\item $\delta_i\not\in \acc(C_{\bar{\tau}})$,
	\item $C_{\bar{\tau}}\cap (\delta_i, \delta_{i+1})\neq \emptyset$, and 
	\item $\delta_0>\max\{\xi, \Delta(\gamma^*, \delta^*)\}$.
	\end{itemize}

For each $i<n$, let $M'_i$ be the model from $\vec M$ with $\sup(M'_i\cap\lambda)=\delta_i$, choosing $\delta_0$ sufficiently large that $\Delta(\gamma^*,\delta^*)\in M'_0$.
For $i+1<n$, choose $\bar\delta_i\in M'_{i+1}\cap\lambda$ above $\sup(C_{\bar\tau}\cap\delta_{i+1})$ and put $J_i=[\delta_i,\bar\delta_i]$. Put $J_{n-1}=[\delta_{n-1},\lambda)$.

Fix $\zeta\in \delta_0\cap M_0'-(\max \{\xi, \max C_{\bar{\tau}}\cap \delta_0, \Delta(\gamma^*, \delta^*)\}+1)$. For each $t\in T$, let $\mathcal{F}_t$ consist of the collection of intervals $I_1<I_2<\cdots<I_n$ such that the following statement $\psi(t, I_1,\cdots, I_n)$ holds: there are $\alpha,\beta\in \Gamma$ such that $t\sqsubseteq C_{\alpha^*}$, $\Delta(\alpha^*, \beta^*)=\Delta(\gamma^*,\delta^*)<\zeta$ and letting $\tau=\otp(C_{\beta^*}\cap \alpha)$ and $\tau'=\Tr(\zeta, \tau)(w)$, we have that $C_{\tau'}\subseteq \zeta \cup \bigcup_{1\leq j \leq n} I_j$ and $C_{\tau'}\cap I_j\neq \emptyset$ for each $1\leq j\leq n$. 
 We say $\mathcal{F}_t$ is \emph{large} if $$\exists^Q I_1 \exists^Q I_2 \cdots \exists^Q I_n \ (I_1, I_2, \cdots, I_n)\in \mathcal{F}_t,$$ where $\exists^Q I$ means that for any $\xi\in \lambda$, there is an interval $I$ with $\min I > \xi$.
Note that for any $t\in M_0'$, $\mathcal{F}_t\in M_0'$.

\begin{claim}\label{claim: goodt}
There exists $t\in M_0'$ such that $\Delta(t, C_{\gamma^*})=\xi'>\zeta$ and $\mathcal{F}_t$ is large.
\end{claim}
\begin{proof}[Proof of the Claim]
Suppose otherwise. First we check that for any $t\in M_0'\cap T$ with $t\sqsubseteq C_{\gamma^*}$, $\mathcal{F}_t$ is large. The intervals $J_0,\ldots,J_{n-1}$ cover $C_{\bar\tau}\setminus\zeta$, and each meets $C_{\bar\tau}$. Thus $\gamma,\delta$ witness $\psi(t,J_0,\ldots,J_{n-1})$.

Starting with $j=n-1$ and proceeding down to $0$, replace $J_j$ by an unboundedly quantified interval. Indeed, $t,\zeta$ and the earlier intervals $J_i$ ($i<j$) belong to $M'_j$, and $\min J_j=\delta_j=\sup(M'_j\cap\lambda)$. For every bound in $M'_j\cap\lambda$, the interval $J_j$ witnesses the statement with the remaining unbounded quantifiers. By elementarity, the statement holds for every bound below $\lambda$. At $j=0$ this gives the largeness of $\mathcal F_t$.

Note that for any $\xi \in M_0'\cap \lambda$, $C_{\gamma^*}\restriction \xi\in M_0'\cap T$. To see this, by \ref{p4}, we can find some $\alpha\in M_0'\cap \Gamma$ such that $\otp(C_{\gamma^*}\cap \alpha)\geq \xi$. Since $C_{\gamma^*}\cap \alpha\in M_0'$ by \ref{a3}, $C_{\gamma^*}\restriction \xi=(C_{\gamma^*}\cap \alpha)\restriction \xi\in M_0'$. The elementarity of $M'_0$ ensures that $C_{\gamma^*}\restriction \xi \in T$.
Let $t_0 = C_{\gamma^*}\restriction(\zeta+1) \in M_0'$. Then $\{t\in T: t_0\sqsubseteq t, \mathcal{F}_t \text{ is large}\}$ defines a cofinal branch on $T$ (since every element in this set is compatible with $C_{\gamma^*}$), which contradicts \ref{T2}.
\end{proof}

Fix $t\in M'_0$ as in Claim~\ref{claim: goodt} and put $\xi'=\Delta(t,C_{\gamma^*})$. By the largeness of $\mathcal F_t$, choose $I_1$ in $M'_0$ with $\min I_1>\dom(t)$. For $0<j<n$, choose $I_{j+1}$ in $M'_j$ with $\min I_{j+1}>\bar\delta_{j-1}$, retaining the remaining unbounded quantifiers at each step. Take the witnesses $\alpha,\beta$ in $M'_{n-1}\cap\Gamma$. Then 
\begin{itemize}
\item $t\sqsubseteq C_{\alpha^*}$, $\Delta(\alpha^*, \beta^*)=\Delta(\gamma^*,\delta^*)<\zeta$. Put $\tilde{\tau}=\otp(C_{\beta^*}\cap \alpha)$ and $\tilde{\tau}'=\Tr(\zeta, \tilde{\tau})(w)$.

\item $\xi'<I_1<J_0<I_2<J_1<\cdots<I_n<J_{n-1}$, where for $n=1$ this reads $\xi'<I_1<J_0$. Indeed, each $I_{j+1}$ for $j+1<n$ is bounded in $\lambda$ and belongs to $M'_j$, so $\sup I_{j+1}<\delta_j$. We may truncate $I_n$ at $\tilde\tau\in M'_{n-1}\cap\lambda$, since $C_{\tilde\tau'}\subseteq\tilde\tau$.
\end{itemize}

Since $t\sqsubseteq C_{\alpha^*}$, we have $\Delta(C_{\alpha^*},C_{\gamma^*})=\xi'$. Above $\xi'$ each $J_i$ contributes one class of $C_{\bar\tau}$. The separating interval $I_{i+2}$ meets $C_{\tilde\tau'}$ whenever $i+1<n$. Hence $osc(C_{\bar\tau}-\xi',C_{\tilde\tau'}-\xi')=n$. By \ref{p3}, we may assume that we have chosen $\beta\in N\cap \Gamma$ large enough to satisfy $C_{\bar{\tau}}\cap (\tilde{\tau}, \otp(C_{\delta^*}\cap \beta))\neq \emptyset$ by \ref{p4}.
Finally, let us check that the tuple $(\alpha,\beta,\gamma,\delta)$ satisfies the hypothesis of \textbf{(3a)}. 

\begin{itemize}
\item By the choice of $\beta$ and the fact that $\lambda_w(\lambda_N, \tau)<\tilde{\tau}$, $$w=\chi(\tilde{\tau},\otp(C_{\delta^*}\cap \beta),\tau).$$
\item Since $\lambda_w(\zeta, \tilde{\tau})<\xi'$, we know $\Tr(\xi', \tilde{\tau})(w)=\Tr(\zeta, \tilde{\tau})(w)=\tilde{\tau}'$.
\item $\Delta(\alpha^*,\gamma^*)=\xi'>\Delta(\alpha^*,\beta^*)= \Delta(\gamma^*,\delta^*)$.
\end{itemize}
Hence \textbf{(3a)} applies and
\[
c(\alpha,\beta,\gamma,\delta)=osc(C_{\bar\tau}-\xi', C_{\tilde\tau'}-\xi')-1=n-1.
\]
As $n$ ranges over the positive integers, $n-1$ ranges over $\omega$.

\subsection{The unbounded case (3b)}\label{unbounded}

Suppose there are $t_0,t_1\in T$ with $\xi=\Delta(t_0,t_1)$ and $t=t_0\wedge t_1$ such that there are $\xi_0 < \xi_1 <  \xi$ such that 
	\begin{itemize}
	\item $\{\rho(\xi,\tau): \exists\,\alpha,\beta\in \Gamma, \tau= \otp(C_{\beta^*}\cap \alpha), t_0\sqsubseteq C_{\beta^*}, \Delta(t, C_{\alpha^*})=\xi_1\}$ and
	\item  $\{\rho(\xi,\tau): \exists\,\gamma,\delta\in \Gamma, \tau= \otp(C_{\delta^*}\cap \gamma) , t_1\sqsubseteq C_{\delta^*}, \Delta(t, C_{\gamma^*})=\xi_0\}$
	\end{itemize}
 are both uncountable.

Since $g: [\omega_1]^2\to \omega$ witnesses $\omega_1\not\to [\omega_1, \omega_1]^2_\omega$,  we can pick $\alpha<\beta<\gamma<\delta \in \Gamma$ such that $\Delta(\alpha^*, \beta^*)=\xi_1$ and $t_0\sqsubseteq C_{\beta^*}$ and $\Delta(\gamma^*, \delta^*)=\xi_0$ and $t_1\sqsubseteq C_{\delta^*}$ such that $g(\rho(\xi, \otp(C_{\beta^*}\cap \alpha)),\rho(\xi, \otp(C_{\delta^*}\cap \gamma)))=n$.
Note that $\Delta(\gamma^*,\delta^*)=\xi_0< \Delta(\alpha^*,\beta^*)=\xi_1<\Delta(\beta^*,\delta^*)=\xi$. Case \textbf{(3b)} applies and we get $c(\alpha,\beta,\gamma,\delta)=n$.
\subsection{The square-bracket operation case (3c)}
We now assume the hypotheses of cases \textbf{(3a)} and \textbf{(3b)} in Subsections \ref{osc} and \ref{unbounded} do not hold.
Recall that we fixed a sequence of countable elementary submodels $\vec{M}$ of $H(\chi)$ before Subsection \ref{osc}.
	\begin{claim}\label{claim: full}
	There exists a club $F\subseteq \omega_1$ such that for any $i\in F$ and any $\tau\not\in M_i\cap \lambda$ for which there are $\gamma<\delta\in \Gamma$ satisfying $\otp(C_{\delta^*}\cap \gamma)=\tau$, we have $E\cap M_i\cap \lambda\subseteq^* \acc(C_{\bar{\tau}})$ where $\bar{\tau}=\last{\lambda_{M_i}}{\tau}$ (in particular, $\otp(C_{\bar{\tau}}\cap \lambda_{M_i})\geq M_i\cap \omega_1$).
	\end{claim}
	\begin{proof}[Proof of the Claim]
 Suppose otherwise. Then we can take a countable elementary submodel $N\prec H((2^\chi)^+)$ containing relevant objects such that $N\cap \omega_1= i\in F^c$. Fix some $\tau \not\in N$ such that there are $\gamma,\delta\in \Gamma$ with $\tau=\otp(C_{\delta^*}\cap \gamma)$. Notice that $N\cap \omega_1 = M_i\cap \omega_1$ and $N\cap \lambda=M_i\cap \lambda$. The latter follows from the fact that $|\lambda|=\omega_1$. In particular, $\bar{\tau}=\last{\lambda_{N}}{\tau}= \last{\lambda_{M_i}}{\tau}$. Hence, we get that $E\cap N\cap \lambda = E\cap M_i\cap \lambda\not \subseteq^* \acc(C_{\bar{\tau}})$. As a result, the hypothesis of \textbf{(3a)} in Subsection \ref{osc} holds, contrary to our assumption.
	\end{proof}

	\begin{definition}
	Let $i\in \omega_1$ be given. We say $(\gamma,\delta)\in [\Gamma]^2$ is \emph{good} for $M_i \in \vec{M}$ if 
		\begin{itemize}
		\item $\gamma^*, \delta^*\not\in M_i$,
		\item $\Delta(\gamma^*, \delta^*)\in M_i$,
		\item letting $\tau=\otp(C_{\delta^*}\cap \gamma)$ and $\bar{\tau}=\last{\lambda_{M_i}}{\tau}$, we have that $\otp(C_{\bar{\tau}}\cap \lambda_{M_i})=M_i\cap \omega_1$.
		\end{itemize}
	\end{definition}

	\begin{claim}\label{claim: stat}
For any countable elementary submodel $N\prec H((2^\chi)^+)$ containing all relevant objects, there exists $(\gamma,\delta)\in [\Gamma]^2$ that is good for $M_i$ for unboundedly many $i\in N\cap \omega_1$ with $d(M_i\cap \omega_1)=n$. 
	\end{claim}
	
	\begin{proof}[Proof of the Claim]
	Suppose for the sake of contradiction that the conclusion is not true. Let us pick some $\alpha,\gamma,\delta\in \Gamma$ such that 
	\begin{itemize}
	\item $\alpha^*, \gamma^*, \delta^*\not\in N$,
	\item $\Delta(C_{\alpha^*}, C_{\delta^*})<\Delta(C_{\gamma^*}, C_{\delta^*}) \in N$, 
	\item for a tail of $i\in N\cap \omega_1$, if $d(M_i\cap \omega_1)=n$, then neither $(\gamma,\delta)$ nor $(\alpha,\delta)$ is good for $M_i$.
	\end{itemize}
Fix $s_0, s_1,t\in T\cap N$ such that 
	\begin{itemize}
	\item $t\sqsubseteq C_{\delta^*}$
	\item $\Delta(s_0,t)=\Delta(C_{\alpha^*}, C_{\delta^*})$ and $s_0\sqsubseteq C_{\alpha^*}$, and 
	\item $\Delta(s_1,t)=\Delta(C_{\gamma^*}, C_{\delta^*})$ and $s_1\sqsubseteq C_{\gamma^*}$.
	\end{itemize}
These choices are possible by \ref{p4} and \ref{a3}.

Let $\tau=\otp(C_{\delta^*}\cap \gamma), \tau' = \otp(C_{\delta^*}\cap \alpha)$ and $\bar{\tau}=\last{\lambda_N}{\tau}, \bar{\tau}' = \last{\lambda_N}{\tau'}$. Then since $E\cap\lambda_N\subseteq^* \acc(C_{\bar{\tau}})\cap\acc(C_{\bar{\tau}'})$ we know that for a tail of $i\in N\cap \omega_1$, it is the case that $\bar{\tau}=\last{\lambda_{M_i}}{\tau}$ and $\bar{\tau}'=\last{\lambda_{M_i}}{\tau'}$. Fix one such $i$ such that $d(M_i\cap \omega_1)=n$ and $s_0,s_1, t\in M_i$.

By the hypothesis that neither $(\alpha,\delta)$ nor $(\gamma,\delta)$ is good for $M_i$, we know that there exists $\bar{\xi}\in M_i\cap \lambda$ such that for any $\xi\in \lambda\cap M_i-(\bar{\xi}+1)$, $\rho(\xi, \tau),\rho(\xi, \tau')\geq M_i\cap \omega_1$. Indeed, the tail of $E$ below $\lambda_{M_i}$ gives order type at least $M_i\cap\omega_1$ in both clubs. Failure of goodness makes both order types strictly larger. Now apply Lemma~\ref{lem:walk-rank} to $\tau$ and $\tau'$.

We claim that there exist $t_0,t_1\in M_i$ such that \[\begin{aligned}\xi&=\Delta(t_0,t_1)=\min\{\Delta(t_0,C_{\delta^*}),\Delta(t_1,C_{\delta^*})\}\\&>\max\{\bar{\xi}, \lambda(\bar{\tau},\tau), \lambda(\bar{\tau}',\tau'), \Delta(s_0,t), \Delta(s_1, t)\},\end{aligned}\] and that, for each $j<2$, the following set is uncountable: $\{\rho(\xi, \tau): \exists \delta,\gamma\in \Gamma, s_j\sqsubseteq C_{\gamma^*}, t_j\sqsubseteq C_{\delta^*}, \tau = \otp(C_{\delta^*}\cap \gamma)\}$. To see this, for each $s\in T$, consider the following property $\varphi(s)$: for any $\xi\in (\bar{\xi},\lambda)\cap \operatorname{ht}(s)$, the following sets are uncountable: $\{\rho(\xi, \tau): \exists \delta,\alpha\in \Gamma, s_0\sqsubseteq C_{\alpha^*}, s\sqsubseteq C_{\delta^*}, \tau = \otp(C_{\delta^*}\cap \alpha)\}$ and $\{\rho(\xi, \tau): \exists \delta,\gamma\in \Gamma, s_1\sqsubseteq C_{\gamma^*}, s\sqsubseteq C_{\delta^*}, \tau = \otp(C_{\delta^*}\cap \gamma)\}$.

Notice that for any proper initial segment $t'\sqsubseteq C_{\delta^*}\restriction \lambda_{M_i}$ of height $>\bar{\xi}$ with $t'\in M_i$, we have that $\varphi(t')$ is true. To see this, fix $\zeta\in (\bar{\xi}, \operatorname{ht}(t'))\cap M_i$ and $\mu<M_i\cap\omega_1$. By the choice of $\alpha, \gamma$, $\delta$ and $\bar{\xi}$, we know that 
		\begin{itemize}
		\item $t'\sqsubseteq C_{\delta^*}$
		\item $s_0\sqsubseteq C_{\alpha^*}$
		\item $s_1\sqsubseteq C_{\gamma^*}$ and 
		\item $\rho(\zeta, \tau), \rho(\zeta, \tau')>\mu$ where $\tau = \otp(C_{\delta^*}\cap \gamma)$ and $\tau' = \otp(C_{\delta^*}\cap \alpha)$.
		\end{itemize}
Applying the elementarity of $M_i$, we get that $\varphi(t')$ holds.

Since $C_{\delta^*}\restriction \lambda_{M_i}\not\in M_i$ and all its proper initial segments are in $T\cap M_i$ and $T$ does not contain an uncountable branch by \ref{T2}, we must be able to find $t_0, t_1\in M_i\cap T$, extending $t$, not compatible with $C_{\delta^*}\restriction \lambda_{M_i}$ with $\max\{\bar{\xi}, \lambda(\bar{\tau},\tau), \lambda(\bar{\tau}',\tau'), \Delta(s_0,t), \Delta(s_1, t)\}<\Delta(t_0,C_{\delta^*})<\Delta(t_1,C_{\delta^*})$ such that $\varphi(t_0)$ and $\varphi(t_1)$ both hold. Note that $\Delta(t_0,t_1)=\Delta(t_0,C_{\delta^*})$. Letting $\xi_0=\Delta(s_0, t)$ and $\xi_1=\Delta(s_1,t)$, we know that $\xi=\Delta(t_0,t_1)>\max\{\xi_0,\xi_1\}$ and \begin{itemize}
	\item $\{\rho(\xi,\tau): \exists\,\alpha,\beta\in \Gamma, \tau= \otp(C_{\beta^*}\cap \alpha), t_0\sqsubseteq C_{\beta^*}, \Delta(t_0\wedge t_1, C_{\alpha^*})=\xi_1\}$ and
	\item  $\{\rho(\xi,\tau): \exists\,\gamma,\delta\in \Gamma, \tau= \otp(C_{\delta^*}\cap \gamma) , t_1\sqsubseteq C_{\delta^*}, \Delta(t_0\wedge t_1, C_{\gamma^*})=\xi_0\}$
	\end{itemize}
 are both uncountable. However, this contradicts our assumption that the hypothesis of \textbf{(3b)} in Subsection \ref{unbounded} does not hold.
	\end{proof}

Let us find some countable $N \prec H((2^\chi)^+)$ containing relevant objects. By Claim \ref{claim: stat}, we may fix $(\gamma,\delta)$ that is good for $M_i$ for unboundedly many $i\in N\cap \omega_1$ with $d(M_i\cap \omega_1)=n$. Let $\tau=\otp(C_{\delta^*}\cap \gamma)$ and $\bar{\tau}=\last{\lambda_N}{\tau}$. Find $i\in N\cap \omega_1$ so large that $d(M_i\cap \omega_1)=n$, $\Delta(\gamma^*, \delta^*)\in M_i$, $\lambda(\bar{\tau}, \tau) \in M_i$ and $\lambda_{M_i}\in \acc(C_{\bar{\tau}})$. In particular, $\bar{\tau}=\last{\lambda_{M_i}}{\tau}$ and $\otp(C_{\bar{\tau}}\cap \lambda_{M_i})=M_i\cap \omega_1$.

Fix $\beta\in \Gamma\cap N$ with $\beta^*\not\in M_i$ such that $\Delta(\beta^*,\gamma^*)\in M_i$ and $$(\lambda_{M_i}, \otp(C_{\delta^*}\cap \beta))\cap C_{\bar{\tau}} \neq \emptyset.$$ Such a $\beta$ exists by \ref{p4}.

Let $\nu=\otp(C_{\bar{\tau}}\cap \otp(C_{\delta^*}\cap \beta))\geq\otp((C_{\bar{\tau}}\cap \lambda_{M_i}) \cup \{\lambda_{M_i}\})>M_i\cap \omega_1$ by Claim \ref{claim: stat} and \ref{p4}. Put $m_0=e_{\nu}(M_i\cap \omega_1)$. Since $\beta^* \not \in M_i$ and $C_{\beta^*}\restriction \xi\in M_i$ for all $\xi\in M_i\cap \lambda$, we can find $s\in T\cap M_i$ such that $\Delta(s,C_{\beta^*})>\max\{\Delta(\beta^*,\gamma^*), \Delta(\gamma^*, \delta^*)\}$ and such that, if $\eta=\otp(C_{\bar{\tau}}\cap\Delta(s,C_{\beta^*}))$, then $e_\nu(\eta)=m\geq m_0$. Finally, since there are uncountably many $\alpha\in \Gamma$ with $s\sqsubseteq C_{\alpha^*}$, we can find $\alpha\in\Gamma\cap M_i$ such that 
	\begin{itemize}
	\item $s\sqsubseteq C_{\alpha^*}$,
	\item $\otp(C_{\bar{\tau}}\cap \otp(C_{\delta^*}\cap \alpha))\in (\max (F_m(\nu)\cap M_i), M_i\cap \omega_1)$, and 
	\item $\otp(C_{\beta^*}\cap \alpha)> \lambda(\bar{\tau}, \tau)$
	\end{itemize}
by \ref{p4}. We have $\otp(C_{\bar{\tau}}\cap \otp(C_{\delta^*}\cap \alpha))\in M_i$ because $\otp(C_{\bar{\tau}}\cap \lambda_{M_i})=M_i\cap \omega_1$ and $\otp(C_{\delta^*}\cap \alpha)\in M_i\cap \lambda$.

Let us verify that $c(\alpha,\beta,\gamma,\delta)$ equals $d(M_i\cap \omega_1)=n$. First of all, note $\Delta(\alpha^*,\beta^*)=\Delta(s,C_{\beta^*})>\max\{\Delta(\beta^*,\gamma^*), \Delta(\gamma^*, \delta^*)\}$. As a result, the definition enters Case \textbf{(3c)}.

Let $\tau'=\otp(C_{\beta^*}\cap \alpha)$. Let $w=\chi(\tau',\otp(C_{\delta^*}\cap \beta),\tau)$.
Since $\lambda(\bar{\tau},\tau)<\tau'<\otp(C_{\delta^*}\cap \beta)<\bar{\tau}<\tau$ and $(\lambda_{M_i}, \otp(C_{\delta^*}\cap \beta))\cap C_{\bar{\tau}} \neq \emptyset$, by Fact \ref{fact: passing}, we have that $\Tr(\tau',\tau)(w)=\bar{\tau}$.

By the definition of $c$ in case \textbf{(3c)}, $$c(\alpha,\beta,\gamma,\delta)=d(\min(F_m(\nu)\setminus\otp(C_{\bar{\tau}}\cap \otp(C_{\delta^*}\cap \alpha))))= d(M_i\cap \omega_1)=n.$$

This concludes the proof of Theorem~\ref{theoerem: strongSpecial}.

\section{Stepping up with a weak Kurepa tree}\label{section: weakKurepa}

Recall that $T\subseteq 2^{<\omega_1}$ is a \emph{weak Kurepa tree} if the height and size of $T$ are both $\omega_1$ and $T$ has at least $\omega_2$ uncountable branches.

\begin{theorem}\label{theorem: wkt}
If there exists a weak Kurepa tree with at least $\omega_3$ uncountable branches, then $\omega_3\not\to [\omega_1]^4_\omega$.
\end{theorem}

Fix the following data for the coloring.
\begin{enumerate}
\item Fix a $C$-sequence $\langle C_\alpha:\alpha<\omega_3\rangle$ such that $C_\alpha$ is a club in $\alpha$ of order type $\cf(\alpha)$. Use this sequence for the walks and, below $\omega_2$, for the rank $\rho$ of Definition~\ref{def:walk-rank}.
\item Fix a weak Kurepa tree $T\subseteq 2^{<\omega_1}$, distinct cofinal branches $\langle b_\alpha:\alpha<\omega_3\rangle$ through $T$, and a bijection $\operatorname{index}_T:T\to\omega_1$. In this section $\Delta(\alpha,\beta)=\Delta(b_\alpha,b_\beta)$ for $\alpha<\beta<\omega_3$. Thus $\Delta(\alpha,\beta)$ is the least $\xi<\omega_1$ such that $b_\alpha(\xi)\neq b_\beta(\xi)$.
\item Fix a coloring $h_3:[\omega_2]^3\to\omega$ witnessing $\omega_2\not\to[\omega_1]^3_\omega$ \cite{Todorcevic3D}.
\item Fix $d:\omega_1\to\omega$ such that $d^{-1}\{n\}$ is stationary for each $n\in\omega$.
\item For each $\beta<\omega_3$, fix an injection $l_\beta:\beta\to\omega_2$. For each $\nu<\omega_1$, fix an injection $e_\nu:\nu\to\omega$. For $m\in\omega$, let $F_m(\nu)=\{\xi<\nu:e_\nu(\xi)\leq m\}$.
\end{enumerate}

Let us now define a function $c: [\omega_3]^4\to \omega$. Fix $\alpha<\beta<\gamma<\delta\in \omega_3$ and let $r=\{\alpha,\beta,\gamma,\delta\}$.
Let $k=\chi(\beta,\gamma,\delta)$. Let $\delta^* = \Tr(\alpha,\delta)(k)$, $\gamma^*=\Tr(\alpha,\gamma)(k)$ and $\beta^*=\Tr(\alpha,\beta)(k)$.
We decide the value $c(\alpha,\beta,\gamma,\delta)$ based on the following cases.

\begin{figure}
    \centering 
   \includegraphics[scale=0.6]{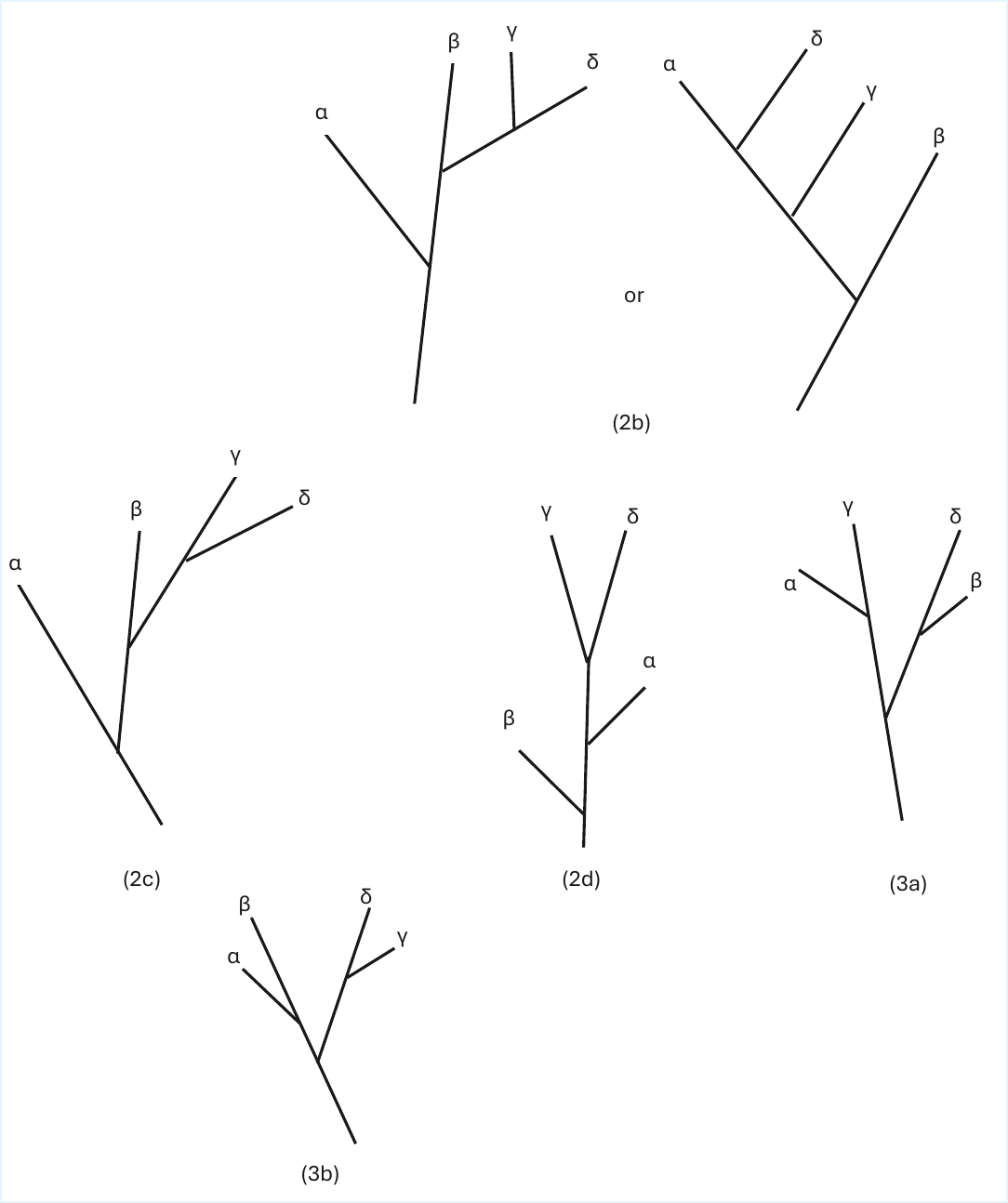}
    \caption{Tuple configurations based on a weak Kurepa tree} 
    \label{fig:kurepa} 
\end{figure}	  
\begin{enumerate}
\item If $(C_{\Tr(\beta,\gamma)(k)}-\beta) \cap (C_{\Tr(\beta,\delta)(k)}-\beta) = \emptyset$, then output $$osc(C_{\Tr(\beta,\gamma)(k)}-\beta,C_{\Tr(\beta,\delta)(k)}-\beta)-1.$$ 
\item Otherwise. 
	\begin{enumerate}
	\item If $\lambda_k(\alpha,\gamma)$, $\lambda_k(\alpha,\beta)$ and $\lambda_k(\alpha,\delta)$ are pairwise distinct, output $$h_3(l_\alpha(\lambda_k(\alpha,\gamma)), l_\alpha(\lambda_k(\alpha,\beta)), l_\alpha(\lambda_k(\alpha,\delta))).$$

\item Else if $\Delta(\alpha, \beta)=\min\{\Delta(\beta,\gamma),\Delta(\alpha, \beta), \Delta(\alpha,\gamma)\}<\Delta(\gamma,\delta)$ and $\otp(C_{\beta^*}\cap \alpha)< \otp(C_{\gamma^*}\cap \alpha)<\otp(C_{\delta^*}\cap \alpha)$, 
	  then output $$h_3(\otp(C_{\beta^*}\cap \alpha), \otp(C_{\gamma^*}\cap \alpha),\otp(C_{\delta^*}\cap \alpha)).$$
	  
	  \item Else if $\Delta(\alpha,\beta)<\Delta(\beta,\gamma)<\Delta(\gamma,\delta)$ and $\otp(C_{\gamma^*}\cap \alpha)\geq \otp(C_{\delta^*}\cap \alpha)$, output $h_3(\Delta(\alpha,\beta),\Delta(\beta,\gamma),\Delta(\gamma,\delta))$.
	  
    \item Else if $\max\{\rho(\otp(C_a\cap\alpha),\otp(C_b\cap\beta)):a,b\in\{\gamma^*,\delta^*\}\}<\Delta(\gamma,\delta)$ and $\Delta(\beta,\gamma)<\Delta(\alpha,\gamma)<\Delta(\gamma,\delta)$, let $\xi_0<\xi_1<\xi_2$ enumerate $\{\Delta(x,y):x,y\in r,\ x<y\}$. Put $m=e_{\xi_2}(\xi_0)$ and output $d(\min(F_m(\xi_2)\setminus\xi_1))$.

	\end{enumerate}
	
\item Otherwise. Let $\tau=\otp(C_{\delta^*}\cap \gamma)$, $\tau'=\otp(C_{\beta^*}\cap \alpha)$ and $$w=\chi(\tau',\otp(C_{\delta^*}\cap \beta),\tau)$$ and let $\bar{\tau}=\Tr(\tau' ,\tau)(w)$. 
Let $\bar{\tau}'=\Tr(\sup C_{\bar{\tau}}\cap \tau',\tau')(w)$. 
	\begin{enumerate}
	\item If $\min\{\Delta(\alpha,\gamma), \Delta(\beta,\delta)\}> \Delta(\alpha, \beta) =\Delta(\gamma, \delta)$, then output $$osc(C_{\bar{\tau}}-\Delta'(\bar{\tau}', \bar{\tau}), C_{\bar{\tau}'}-\Delta'(\bar{\tau}', \bar{\tau}))-2$$ if the value is non-negative, and output $0$ otherwise. Recall from Definition~\ref{def: deltaprime} that $\Delta'(a,b)=\max \{\xi\leq \sup C_b \cap a: \otp(C_a\cap \xi)=\otp(C_b\cap \xi)\}$.
	\item Else if $\Delta(\beta,\delta)< \Delta(\alpha,\beta)<\Delta(\gamma,\delta)$, output $$d(\min(F_m(i)\setminus \operatorname{index}_T(b_\gamma\restriction \nu))),$$ where $\nu=\rho(\otp(C_{\beta^*}\cap \alpha),\otp(C_{\delta^*}\cap \beta))$, $i=\operatorname{index}_T(b_\delta\restriction \nu)$ and $m=e_i(\Delta(\gamma,\delta))$.

	\end{enumerate}
	
\end{enumerate}

\begin{theorem}\label{theoerem: strongwKurepa}
The coloring $c: [\omega_3]^4\to \omega$ defined above satisfies the following: for any $\Gamma\in [\omega_3]^{\omega_1}$ and any $n\in \omega$, there are $\alpha<\beta<\gamma<\delta\in \Gamma$ such that $c(\alpha,\beta,\gamma,\delta)=n$.
\end{theorem}

The proof of Theorem \ref{theoerem: strongwKurepa} occupies the rest of this section.

Let $\Gamma\in [\omega_3]^{\omega_1}$ and $n\in \omega$ be given and let $\varepsilon=\sup \Gamma$. 
Cases \textbf{(1)}, \textbf{(2a)}, \textbf{(2b)}, and \textbf{(2c)} follow the arguments for \textbf{(1)}, \textbf{(2a)}, \textbf{(2c)}, and \textbf{(2d)}, respectively, in the proof of Theorem~\ref{theoerem: strongSpecial}. Since $T$ splits in a binary fashion, there is no analogue of \textbf{(2b)} of that theorem.

Let $\theta$ be a large enough regular cardinal. Let $\vec{M}=\langle M_i:i<\omega_1\rangle$ be an $\in$-increasing sequence of elementary submodels of $H(\theta)$ containing all relevant objects, and let $D=\{\sup(M_i\cap\varepsilon):i<\omega_1\}$.
We use the notation of Theorem~\ref{theoerem: strongSpecial}. After pruning $\Gamma$ as in that proof, we have the following properties.

\begin{enumerate}
\item Properties \ref{b1}, \ref{b2}, \ref{b3}, \ref{b4}.
\item Properties \ref{p1}, \ref{p2}, \ref{p3}, \ref{p4}.
\end{enumerate}

\begin{enumerate}[label=\textbf{A'\arabic*}]
	\item \label{a'1} there is $\eta<\min \Gamma$ such that for all $\alpha<\beta \in \Gamma$, $\lambda_k(\alpha,\beta)=\eta$.
    \item \label{a'2} for each $\alpha\in\Gamma$, the row $A_\alpha=\{\otp(C_{\delta^*}\cap\alpha):\delta\in\Gamma,\ \alpha<\delta\}$ is countable. For each $\xi\in A_\alpha$, there are uncountably many $\delta\in\Gamma$ above $\alpha$ such that $\xi=\otp(C_{\delta^*}\cap\alpha)$.
	\end{enumerate}

Consider the subtree $S$ of $T$: $t\in S$ iff there are uncountably many $\alpha\in \Gamma$ such that $t\sqsubseteq b_\alpha$. If $S$ contains an uncountable branch, then we enter case \textbf{(2c)}, arguing exactly as in Theorems \ref{theoerem: strongSpecial} and \ref{theorem: omega2}. We may assume that $S$ does not have an uncountable branch.

\textbf{(2d):} Suppose that $S$ is an Aronszajn tree. Let $N\prec H(\theta)$ be a countable elementary submodel containing all relevant objects such that $d(N\cap\omega_1)=n$. Since $S$ has height $\omega_1$, we can find $\gamma,\delta\in\Gamma$ such that $\Delta(\gamma,\delta)=\nu>N\cap\omega_1$ and $b_\gamma\restriction\nu\in S$. Put $m_0=e_\nu(N\cap\omega_1)$.
Note that $b_\gamma\restriction(N\cap\omega_1)\subseteq N$. By elementarity, choose $s_0\in S\cap N$ such that $\xi=\Delta(s_0,b_\gamma)$ satisfies $m=e_\nu(\xi)\geq m_0$. Again by elementarity, choose $s_1\in S\cap N$ such that $\xi'=\Delta(s_1,b_\gamma)>\xi$ and $\xi'$ is above every member of $F_m(\nu)\cap N$.
Finally, pick $\alpha<\beta$ in $\Gamma\cap N$ such that $s_1\sqsubseteq b_\alpha$ and $s_0\sqsubseteq b_\beta$. The values of $\rho$ on $A_\alpha\times A_\beta$ lie in $N\cap\omega_1<\nu$, and $\xi=\Delta(\beta,\gamma)<\xi'=\Delta(\alpha,\gamma)<\nu=\Delta(\gamma,\delta)$. Thus \textbf{(2d)} applies and
\[
c(\alpha,\beta,\gamma,\delta)=d(\min(F_m(\nu)\setminus\xi'))=d(N\cap\omega_1)=n.
\]

\textbf{(3):} Suppose $S$ is not an Aronszajn tree. Shrinking $\Gamma$ to an uncountable subset if necessary, we may assume that $S$ is a countable normal perfect subtree of $T$ of height $\eta\in \acc(\omega_1)\cap \omega_1$ such that
		\begin{enumerate}[label=\textbf{Cantor\arabic*}]
	\item \label{cantor1} for each $t$ that is a cofinal branch through $S$, there is at most one $\alpha\in \Gamma$ such that $t\sqsubseteq b_\alpha$, and 
	\item \label{cantor2} for each $\alpha\in \Gamma$, $b_\alpha\restriction \eta$ is a cofinal branch through $S$.
	\end{enumerate}

To see this, each $\Sigma\in [\Gamma]^{\omega_1}$ is naturally associated with a subtree $S_\Sigma$ of $T$, namely, $t\in S_\Sigma$ iff there are uncountably many $\alpha\in \Sigma$ such that $t\sqsubseteq b_\alpha$. Notice that $S_{\Sigma_0}\subseteq S_{\Sigma_1}$ whenever $\Sigma_0\subseteq \Sigma_1$. Take $\Sigma\in [\Gamma]^{\omega_1}$ such that $\operatorname{ht}(S_\Sigma)$ is least, in the sense that for any $\Sigma'\in [\Sigma]^{\omega_1}$, $\operatorname{ht}(S_{\Sigma'})=\operatorname{ht}(S_{\Sigma})$. We claim that such a $\Sigma$ is as desired. Let $S=S_{\Sigma}$ and let $\eta=\operatorname{ht}(S)$. Note that $\eta$ is a limit ordinal by the definition. 

	\begin{claim}
	For all but countably many $\alpha\in \Sigma$, $b_\alpha\restriction \eta$ is a branch through $S$.
	\end{claim}
	
	\begin{proof}[Proof of the Claim]
	Suppose toward a contradiction that there is $\Sigma_0\in [\Sigma]^{\omega_1}$ such that for each $\alpha\in \Sigma_0$, there is some least $\nu_\alpha<\eta$ such that $b_\alpha\restriction \nu_\alpha\not\in S$. By the pigeonhole principle, we may find $\Sigma_1\in [\Sigma_0]^{\omega_1}$ and $\nu<\eta$ such that $\nu_\alpha=\nu$ for all $\alpha\in \Sigma_1$. By the minimality of $\eta$, we know that $\operatorname{ht}(S_{\Sigma_1})=\eta$. Take some $t\in S_{\Sigma_1}$ of height $\nu$. By the definition, there are uncountably many $\alpha\in \Sigma_1$ such that $t\sqsubseteq b_\alpha$. But since $\alpha\in \Sigma_1$ and $\nu_\alpha=\nu$, by the choice of $\Sigma_0\supseteq \Sigma_1$, we have that $t=b_\alpha\restriction \nu \not \in S\supseteq S_{\Sigma_1}$, which is a contradiction.
	\end{proof}
With the claim above, it is easy to see that
	$S$ is 	
		\begin{itemize}
		\item \emph{normal}: for any $t\in S$ and any $\operatorname{ht}(t)\leq\xi<\eta$, there is an extension $t'\in S$ of $t$ such that $\operatorname{ht}(t')=\xi$.
		\item \emph{perfect}: for any $t\in S$, there are incompatible $t_0,t_1\in S$ extending $t$.
		\end{itemize}

For the rest of the proof, we assume $\Gamma$ and its derived tree $S=S_{\Gamma}$ satisfy properties \ref{cantor1} and \ref{cantor2}.

\subsection{The downstairs oscillation case (3a)}

Before we embark on the proof, let us mention that the overall proof structure is similar to that of Proposition \ref{prop: deltaprime}. The added difficulty comes from the fact that in order to determine an ordinal in $\omega_2$ we need to use two ordinals in $\omega_3$. As a result, the reflection arguments become a bit more delicate. 

Recall that $W=\bigcup_{\alpha\in\Gamma}A_\alpha$, that $\lambda=\sup W$, and that \ref{p1}--\ref{p4} hold.
Let $\vec{D}=\langle M_\alpha: \alpha<\omega_1\rangle$ be an internally approachable sequence of countable elementary submodels of $H(\chi)$ containing the relevant objects, where $\chi$ is a sufficiently large regular cardinal. Let $E=\{\sup \lambda\cap M_\alpha: \alpha<\omega_1\}$. Then $E$ is a club in $\lambda$.
Suppose there is an elementary submodel $N\prec H((2^\chi)^+)$ containing the relevant objects and a pair $\gamma<\delta$ in $\Gamma$ with the following property. Put $\lambda_N=\sup(N\cap\lambda)$, $\tau=\otp(C_{\delta^*}\cap\gamma)\notin N$, and $\bar{\tau}=\last{\lambda_N}{\tau}$. Then $\otp(C_{\bar{\tau}}\cap\lambda_N)<N\cap\omega_1$. This holds in particular when $\rho(\lambda_N,\tau)<N\cap\omega_1$. Indeed, when $\bar\tau>\lambda_N$, Lemma~\ref{lem:walk-rank} gives $\otp(C_{\bar\tau}\cap\lambda_N)\leq\rho(\lambda_N,\tau)$, and when $\bar\tau=\lambda_N$, the order type is $\omega$. Let $l\in\omega$ be such that $\bar{\tau}=\Tr(\lambda_N,\tau)(l)$.

By the order type bound, we can choose a sequence $\langle \lambda_i\in N: i\in \omega\rangle$ cofinal in $\lambda_{N}$ and a sequence of nonzero ordinals $\langle \xi_i\in N\cap \omega_1 : i\in \omega\rangle$ such that
	\begin{itemize}
	\item $\lambda_0 > \lambda_l(\lambda_N, \tau)$,
	\item  $\lambda_i\in E- C_{\bar{\tau}}$,
	\item $(\lambda_i, \lambda_{i+1})\cap C_{\bar{\tau}}\neq \emptyset$,
	\item $\otp(C_{\bar{\tau}}\cap \lambda_0)=\xi_0$,
	\item $\otp((\lambda_i, \lambda_{i+1})\cap C_{\bar{\tau}})=\xi_{i+1}$.
	\end{itemize}
Let $\langle M'_i: i\in \omega\rangle$ be models witnessing that $\lambda_i\in E$ for $i\in \omega$. We may also assume that $\otp(C_{\bar{\tau}}\cap \lambda_N)\in M_0'$ (as a result, $\xi_i\in M'_0$ for all $i\in \omega$). Let $\bar{\lambda}\in M_0'$ be such that $$\bar{\lambda}> \max\{\sup C_{\bar{\tau}}\cap \lambda_0, \lambda_l(\lambda_N, \tau) \}.$$

Let $\exists^Q I$ mean that for any $\xi<\lambda$, there exists an interval $I\subseteq \lambda$ with $\min I > \xi$. Let $n$ be a positive integer, the number of intervals in the reflection. Let $s,t\in S$ be such that $s\sqsubseteq b_\gamma, t\sqsubseteq b_\delta$ and $\Delta(s,t)=\Delta(\gamma,\delta)$. Note that $S\in N$ is countable, so $S\subseteq N$. The standard reflection argument (see Proposition \ref{prop: deltaprime}) using the elementarity of $M_0', \cdots, M_n', N$ gives the following statement $\psi$: $\exists^Q I_1 \exists^Q I_2 \cdots \exists^Q I_n$, and there exist $\epsilon\in I_n$ and $\tau'=\otp(C_{\beta^*}\cap \alpha)$ for some $\alpha,\beta\in \Gamma$ with $\bar{\tau}'=\Tr(\epsilon, \tau')(l)$ satisfying the following conditions.
	\begin{itemize}
	\item $s\sqsubseteq b_\alpha$ and $t\sqsubseteq b_\beta$,
	\item $\bar{\lambda} > \lambda_l(\epsilon, \tau')$,
	\item $\otp(C_{\bar{\tau}'}\cap \bar{\lambda}) =\xi_0$,
	\item $\otp(C_{\bar{\tau}'}\cap I_i)=\xi_i$ for $i=1,2,\cdots, n-1$,
	\item $C_{\bar{\tau}'}\subseteq \bar{\lambda} \cup \bigcup_{i=1}^n I_i$ and $C_{\bar{\tau}'}\cap I_n\neq\emptyset$.
	\end{itemize}

Note that the parameters that appear in $\psi$ are all in $M'_0$. Hence, the elementarity of $M'_0$ ensures that $M_0'\models \psi$.
For $1\leq j\leq n$, let $J_j=(\lambda_{2j-2},\sup(C_{\bar{\tau}}\cap\lambda_{2j})+1)$. Reflecting $\psi$ in $M_2', M_4', \cdots, M'_{2n}$, we obtain in $N$ intervals $I_1<I_2<\cdots<I_n$ and the witnessing $\epsilon\in I_n$ and $\alpha,\beta\in \Gamma\cap N$ with the following properties. 
	\begin{itemize}
	\item  $J_1<I_1<\lambda_2 <J_2 < I_2<\lambda_4< \cdots<J_n < I_n <\lambda_{2n}$ and 
	\item $\tau', \bar{\tau}' \in \lambda_{2n}-\min I_n$.
	\end{itemize}

As in Proposition~\ref{prop: deltaprime}, $osc(C_{\bar{\tau}}-\lambda_0, C_{\bar{\tau}'}-\lambda_0)=n+1$, including the class above $I_n$. Since there are uncountably many $\beta\in \Gamma$ satisfying $\otp(C_{\beta^*}\cap \alpha)=\tau'$, we may pick such a $\beta\in N\cap \Gamma$ large enough that $C_{\bar{\tau}}\cap (\tau', \otp(C_{\delta^*}\cap \beta))\neq \emptyset$ by \ref{p4}.

It remains to verify that we can recover all the data from the tuple $(\alpha,\beta, \gamma,\delta)$. Indeed, if $w=\chi(\tau',\otp(C_{\delta^*}\cap \beta),\tau)$, then by Fact \ref{fact: passing}, the definition of the function $\chi$ and the choice of $\beta$, we get that $\bar{\tau}=\Tr(\tau', \tau)(w)$ and $w=l$. Letting $\xi'=\sup(C_{\bar{\tau}}\cap \tau')$, we note that 

$$\tau'>\epsilon>\xi'>\lambda_{2n-1} >\lambda_0> \lambda_l(\epsilon,\tau').$$

 As a result of Fact \ref{fact: passing}, $\Tr(\xi', \tau')(w)=\Tr(\epsilon, \tau')(l)=\bar{\tau}'$. The computation of $\Delta'(\bar{\tau}', \bar{\tau})$ is a key point of the argument.
 
\begin{claim}
$\Delta'(\bar{\tau}', \bar{\tau})=\min(C_{\bar{\tau}}\cap J_1)$.
\end{claim} 
 \begin{proof}
  The proof is the same as that of Claim \ref{claim: computedeltaprime} in Proposition \ref{prop: deltaprime}.
 \end{proof}

Finally, let us check that the hypothesis of \textbf{(3a)} is satisfied. Since $s\sqsubseteq b_\alpha, b_\gamma$ and $t\sqsubseteq b_\beta, b_\delta$, we have that $\Delta(\alpha,\gamma)\geq |s|$ and $\Delta(\beta,\delta)\geq |t|$. As $\Delta(s,t)=\Delta(\gamma,\delta)$, we have that $$\Delta(\alpha,\beta)=\Delta(\gamma,\delta)=\Delta(s,t)<\min \{|s|,|t|\}\leq \min \{\Delta(\alpha,\gamma),\Delta(\beta,\delta)\}.$$

Hence \textbf{(3a)} applies and
\[
\begin{aligned}
c(\alpha,\beta,\gamma,\delta)
 &=osc(C_{\bar{\tau}}-\min(C_{\bar{\tau}}\cap J_1),C_{\bar{\tau}'}-\min(C_{\bar{\tau}}\cap J_1))-2\\
 &=osc(C_{\bar{\tau}}-\lambda_0,C_{\bar{\tau}'}-\lambda_0)-2=(n+1)-2=n-1.
\end{aligned}
\]
As $n$ ranges over the positive integers, $n-1$ ranges over $\omega$.
\begin{figure}
    \centering 
   \includegraphics[scale=1]{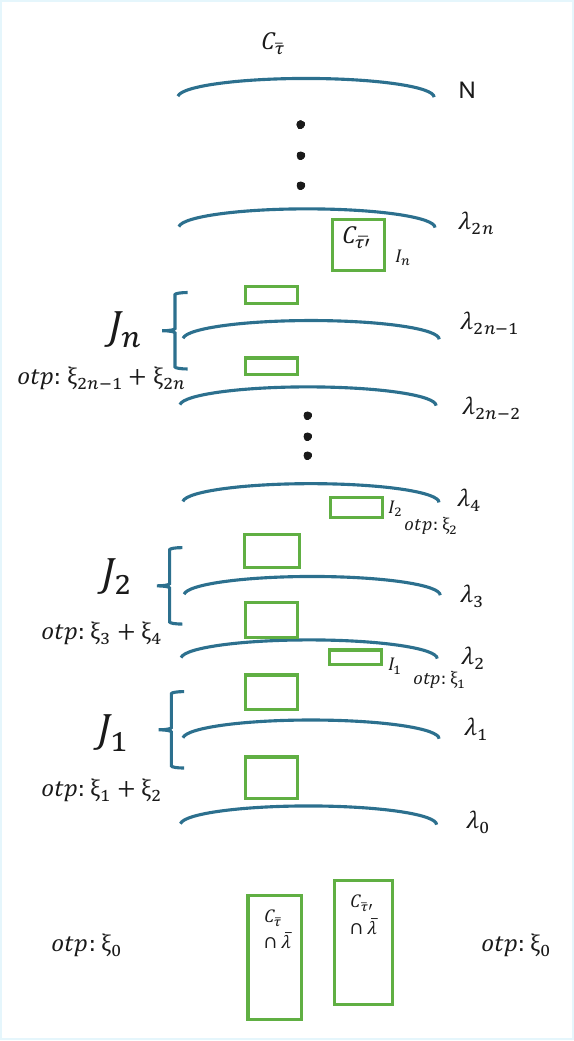}
    \caption{A demonstration of the downstairs oscillation (3a)} 
    \label{fig:oscillation} 
\end{figure}

\subsection{The unbounded case (3b)}\label{subsection: unbounded(3b)}

Suppose for any elementary submodel $N\prec H((2^\chi)^+)$ containing relevant objects and for any $\gamma,\delta\in \Gamma$ with $\tau= \otp(C_{\delta^*}\cap \gamma)\not\in N$, it is the case that if $\bar{\tau} = \last{\lambda_N}{\tau}$, then $\otp(C_{\bar{\tau}}\cap \lambda_N)\geq N\cap \omega_1$.

Fix a countable elementary submodel $N_0\prec H((2^\chi)^+)$ containing the relevant objects. Take a countable $X\prec H((2^\chi)^+)$ with $N_0\in X$ and choose $\beta\in\Gamma$ above $\sup(\Gamma\cap X)$. By \ref{b3}, \ref{p4} and monotonicity, every $\tau\in A_\beta$ is at least $\lambda_X>\lambda_{N_0}$. Since $S$ and $A_\beta$ are countable, by \ref{a'2} we may fix $\tau$ and $\xi<\eta$ such that there are uncountably many $\delta\in\Gamma$ above $\beta$ with $\otp(C_{\delta^*}\cap\beta)=\tau$ and $\Delta(\beta,\delta)=\xi$.
We call such $\delta\in \Gamma$ \emph{good}. The hypothesis gives $\otp(C_{\last{\lambda_{N_0}}{\tau}}\cap\lambda_{N_0})\geq N_0\cap\omega_1$. Recall that $\eta$ is the height of $S$.

Next we find $\alpha\in N_0\cap \Gamma$ large enough with $\Delta(\alpha,\beta)>\xi$ such that if $\tau_0 = \otp(C_{\beta^*}\cap \alpha)$, then $\rho(\tau_0, \tau)=\nu>\eta$. Indeed, for all sufficiently large $y<\lambda_{N_0}$ the walk from $\tau$ to $y$ passes through $\bar\tau=\last{\lambda_{N_0}}{\tau}$. Thus Lemma~\ref{lem:walk-rank} gives $\rho(y,\tau)\geq\otp(C_{\bar\tau}\cap y)$, which eventually exceeds $\eta<N_0\cap\omega_1$. We can take $y=\otp(C_{\beta^*}\cap\alpha)$ by \ref{p4} and monotonicity, restricting $\alpha$ to the cone agreeing with $b_\beta$ through level $\xi$.

Find another countable elementary submodel $N\prec H((2^\chi)^+)$ containing $N_0$ and $\beta$ such that $d(N\cap \omega_1)=n$.
We then find some good $\delta\in \Gamma-N$ such that $i=\operatorname{index}_T(b_\delta\restriction \nu)>N\cap \omega_1$. Put $m_0=e_i(N\cap \omega_1)$. Find $s\in S$ such that $\Delta(b_\delta,s)>\Delta(\alpha,\beta)$ and $m=e_i(\Delta(b_\delta,s))\geq m_0$. The reason why such $s$ exists is that if it did not, since $\{b_\delta\restriction \xi: \xi<\eta\} \subseteq N \cap S$, we would have been able to define $b_\delta\restriction \eta$ in $N$, contrary to the fact that $\delta\not\in N$. 

Since $s\in S\cap N$, there are uncountably many $\gamma\in\Gamma$ with $s\sqsubseteq b_\gamma$. By elementarity choose such a $\gamma\in\Gamma\cap N$ above $\beta$ with $\operatorname{index}_T(b_\gamma\restriction\nu)>\sup(F_m(i)\cap N)$.
We have $\alpha<\beta<\gamma<\delta$. Since $\nu,\gamma\in N$, $\operatorname{index}_T(b_\gamma\restriction\nu)\in N\cap\omega_1$. Finally, $\Delta(\beta,\delta)<\Delta(\alpha,\beta)<\Delta(\gamma,\delta)$ and $\rho(\otp(C_{\beta^*}\cap\alpha),\otp(C_{\delta^*}\cap\beta))=\rho(\tau_0,\tau)=\nu$. Thus \textbf{(3b)} applies and
\[
c(\alpha,\beta,\gamma,\delta)=d(\min(F_m(i)\setminus \operatorname{index}_T(b_\gamma\restriction\nu)))=d(N\cap\omega_1)=n.
\]

\section{Stepping up with a large continuum}\label{section: continuum}

The main result of this section is the following. 
\begin{theorem}\label{theorem: continuum}
If $\omega_3\leq 2^\omega$, then $\omega_3\not\to [\omega_1]^4_{\omega}$ holds.
\end{theorem}

Fix the following data for the coloring.
\begin{enumerate}
\item Fix a $C$-sequence $\langle C_\alpha:\alpha<\omega_3\rangle$ such that $C_\alpha$ is a club in $\alpha$ of order type $\cf(\alpha)$. Use this sequence for the walks and, below $\omega_2$, to define the function $\rho$ of Definition~\ref{def:walk-rank}.
\item Fix distinct reals $\langle r_\alpha:\alpha<\omega_3\rangle$ in $2^\omega$. In this section $\Delta(\alpha,\beta)=\Delta(r_\alpha,r_\beta)$ for $\alpha<\beta<\omega_3$. Thus $\Delta(\alpha,\beta)$ is the least $m\in\omega$ such that $r_\alpha(m)\neq r_\beta(m)$.
\item Fix a coloring $h_3:[\omega_2]^3\to\omega$ witnessing $\omega_2\not\to[\omega_1]^3_\omega$ and a coloring $g:[\omega_1]^2\to\omega$ witnessing $\omega_1\not\to[\omega_1,\omega_1]^2_\omega$, as in Section~\ref{section: weaksquare}.
\item Fix $d:\omega_1\to\omega$ such that $d^{-1}\{n\}$ is stationary for each $n\in\omega$.
\item For each $\beta<\omega_3$, fix an injection $l_\beta:\beta\to\omega_2$. For each $\nu<\omega_1$, fix an injection $e_\nu:\nu\to\omega$. For $m\in\omega$, let $F_m(\nu)=\{\xi<\nu:e_\nu(\xi)\leq m\}$.
\item For $\beta<\omega_2$ and $\nu<\omega_1$, let $G_\nu(\beta)=\{\alpha<\beta:\rho(\alpha,\beta)\leq\nu\}$. This set is countable by Lemma~\ref{lem:walk-rank}. Fix an increasing sequence of finite sets $\langle K_m^{\beta,\nu}:m\in\omega\rangle$ with union $G_\nu(\beta)$.
\end{enumerate}

Let us now define a function $c: [\omega_3]^4\to \omega$. Fix $\alpha<\beta<\gamma<\delta\in \omega_3$.
Let $k=\chi(\beta,\gamma,\delta)$. Let $\delta^* = \Tr(\alpha,\delta)(k)$, $\gamma^*=\Tr(\alpha,\gamma)(k)$ and $\beta^*=\Tr(\alpha,\beta)(k)$.
Use the first applicable clause below to define $c(\alpha,\beta,\gamma,\delta)$. If the value is undefined, output $0$.

\begin{enumerate}
\item If $(C_{\Tr(\beta,\gamma)(k)}-\beta) \cap (C_{\Tr(\beta,\delta)(k)}-\beta) = \emptyset$, then output $$osc(C_{\Tr(\beta,\gamma)(k)}-\beta,C_{\Tr(\beta,\delta)(k)}-\beta)-1.$$ 
\item Otherwise. 
	\begin{enumerate}
	\item If $\lambda_k(\alpha,\gamma)$, $\lambda_k(\alpha,\beta)$ and $\lambda_k(\alpha,\delta)$ are pairwise distinct, output $$h_3(l_\alpha(\lambda_k(\alpha,\gamma)), l_\alpha(\lambda_k(\alpha,\beta)), l_\alpha(\lambda_k(\alpha,\delta))).$$

\item Else if $\Delta(\alpha, \beta)=\min\{\Delta(\beta,\gamma),\Delta(\alpha, \beta), \Delta(\alpha,\gamma)\}<\Delta(\gamma,\delta)$ and $\otp(C_{\beta^*}\cap \alpha)< \otp(C_{\gamma^*}\cap \alpha)<\otp(C_{\delta^*}\cap \alpha)$, 
	  then output $$h_3(\otp(C_{\beta^*}\cap \alpha), \otp(C_{\gamma^*}\cap \alpha),\otp(C_{\delta^*}\cap \alpha)).$$
	  
	  \item Else if $\Delta(\alpha,\beta)<\Delta(\beta,\gamma)<\Delta(\gamma,\delta)$ and $\otp(C_{\gamma^*}\cap \alpha)\geq \otp(C_{\delta^*}\cap \alpha)$, output $h_3(\Delta(\alpha,\beta),\Delta(\beta,\gamma),\Delta(\gamma,\delta))$.

	\end{enumerate}
	
\item Otherwise. Let $\tau=\otp(C_{\delta^*}\cap \gamma)$, $\tau'=\otp(C_{\beta^*}\cap \alpha)$ and $$w=\chi(\tau',\otp(C_{\delta^*}\cap \beta),\tau)$$ and let $\bar{\tau}=\Tr(\tau' ,\tau)(w)$. 
Let $\bar{\tau}'=\Tr(\sup C_{\bar{\tau}}\cap \tau',\tau')(w)$. 
	\begin{enumerate}
	\item If $\min\{\Delta(\alpha,\gamma), \Delta(\beta,\delta)\}> \Delta(\alpha, \beta) =\Delta(\gamma, \delta)$, then output $$osc(C_{\bar{\tau}}-\Delta'(\bar{\tau}', \bar{\tau}), C_{\bar{\tau}'}-\Delta'(\bar{\tau}', \bar{\tau}))-2$$ if the value is non-negative and output $0$ otherwise. Recall from Definition~\ref{def: deltaprime} that $\Delta'(a,b)=\max \{\xi\leq \sup C_b \cap a: \otp(C_a\cap \xi)=\otp(C_b\cap \xi)\}$.
    \item Else if $\Delta(\gamma,\delta)<\Delta(\beta,\gamma)<\Delta(\alpha,\beta)$, put
    \[
    \begin{aligned}
    m&=\Delta(\alpha,\beta), & \tau_0&=\otp(C_{\gamma^*}\cap\beta),\\
    \nu&=\rho(\tau_0,\tau), & u&=\min(K_m^{\tau_0,\nu}\setminus\tau').
    \end{aligned}
    \]
    Let $\xi=\min(F_m(\nu)\setminus\rho(\tau',u))$ and output $g(\xi,\nu)$.

    \item Else if $\Delta(\gamma,\delta)<\Delta(\beta,\delta)<\Delta(\alpha,\delta)$, put $m=\Delta(\alpha,\delta)$ and $\nu=\rho(\otp(C_{\gamma^*}\cap\beta),\tau)$. Output
    \[
      d\bigl(\min(F_m(\nu)\setminus\rho(\otp(C_{\delta^*}\cap\alpha),\tau))\bigr).
    \]
	\end{enumerate}
	
\end{enumerate}

\begin{theorem}\label{theoerem: continuumVerify}
The coloring $c: [\omega_3]^4\to \omega$ defined above satisfies the following: for any $\Gamma\in [\omega_3]^{\omega_1}$ and any $n\in \omega$, there are $\alpha<\beta<\gamma<\delta\in \Gamma$ such that $c(\alpha,\beta,\gamma,\delta)=n$. 
\end{theorem}

The proof of Theorem~\ref{theoerem: continuumVerify} occupies the rest of this section.

Fix $\Gamma\in[\omega_3]^{\omega_1}$ and $n\in\omega$. We may assume that $\otp(\Gamma)=\omega_1$. Put $\varepsilon=\sup\Gamma$. For cases \textbf{(1)}, \textbf{(2a)} and \textbf{(2b)} we follow the proof of Theorem~\ref{theoerem: strongwKurepa}. 
We may assume there are a club $D\subseteq\varepsilon$, a $k\in\omega$ and a function $b\mapsto b^*$ on $\Gamma$ such that the following hold. 
\begin{enumerate}
\item \emph{The upper clubs are trivial at $\varepsilon$.} The pruning gives $D\cap a\subseteq C_{b^*}$ whenever $a<b$ are in $\Gamma$. Thus for every $\xi<\varepsilon$ some member of $\{C_{b^*}:b\in\Gamma\}$ contains $D\cap\xi$. For each $a<b$ in $\Gamma$, \ref{b3} gives a $\nu\in\acc(D)\cap\acc(C_{b^*})$ with $a<\nu<b^*$. By \ref{b4}, $(C_{b^*}\setminus a)\cap(C_{c^*}\setminus a)\neq\emptyset$ for $a<b<c$ in $\Gamma$.
\item \emph{The walk parameters are fixed.} The function $\chi$ is constantly $k$ on $[\Gamma]^3$. There is an $\eta<\min\Gamma$ such that for all $a<b$ in $\Gamma$,
\[
 \Tr(a,b)(k)=b^*,\qquad \lambda_k(a,b)=\eta
 \qquad\text{(\ref{b1}, \ref{b2}, \ref{a'1}).}
\]
\end{enumerate}

Let $S$ consist of the $t\in 2^{<\omega}$ extended by uncountably many $r_a$ with $a\in\Gamma$. We discard the countably many $a$ for which some finite initial segment of $r_a$ lies outside $S$. For $a<b$ in $\Gamma$ put
\[
 f(a,b)=\otp(C_{b^*}\cap a),\qquad A_a=\{f(a,b):b\in\Gamma,\ a<b\}.
\]
Let $W=\bigcup_{a\in\Gamma}A_a$ and $\lambda=\sup W$. For sufficiently large regular $\theta$, fix a continuous $\in$-increasing chain $\vec M=\langle M_i:i<\omega_1\rangle$ of countable elementary submodels of $H(\theta)$ containing $\Gamma$, the map $b\mapsto b^*$, $D,S,f$ and the fixed coloring data. The club in $\lambda$ used in \textbf{(3a)} is $\{\lambda_{M_i}:i<\omega_1\}$, where $\lambda_{M_i}=\sup(M_i\cap\lambda)$.

We are now in the following situation.
\begin{enumerate}

\item \emph{The rows are countable.} For every $a\in\Gamma$, $A_a$ is countable and each member of $A_a$ occurs as $f(a,b)$ for uncountably many $b>a$ in $\Gamma$ (\ref{a'2}, \ref{p3}). The function $f$ is monotone, increasing and locally small on $\Gamma$ in the sense of Section~\ref{section: pruning}.
\item \emph{The order types are cofinal in $\lambda$.} We have $\otp(W)=\cf(\lambda)=\omega_1$ and $\sup f''[\Sigma]^2=\lambda$ for every uncountable $\Sigma\subseteq\Gamma$ (\ref{p1}, \ref{p2}). For every countable $N\prec H((2^\theta)^+)$ containing $\theta,\Gamma,D,f$ and the other fixed parameters, and every $b\in\Gamma$ with $b^*\notin N$, the values $f(a,b)$ for $a\in\Gamma\cap N$, $a<b$, are cofinal in $\lambda_N$ (\ref{p4}).
\item \emph{The reals give distinct branches through $S$.} The tree $S$ is countable, normal and perfect of height $\omega$. This rules out the uncountable-branch and Aronszajn-tree alternatives of Section~\ref{section: weakKurepa}.
\item \emph{The hypothesis of \textbf{(3a)} fails.} For every countable $N\prec H((2^\theta)^+)$ containing the fixed data and $\vec M$, and every $\gamma<\delta$ in $\Gamma$ with $\tau=f(\gamma,\delta)>\lambda_N$, put $\bar\tau=\last{\lambda_N}{\tau}$. Then
\[
 \otp(C_{\bar\tau}\cap\lambda_N)\geq N\cap\omega_1 .
\]
\end{enumerate}
It remains to verify \textbf{(3b)} and \textbf{(3c)}.

\subsection{The unbounded case (3b)}\label{subsection: continuumunbounded(3b)}

Suppose there is a stationary $Z\subseteq\omega_1$ such that for every $i\in Z$ there is a $\zeta<\sup\Gamma$ with the following property. Whenever $\zeta<\gamma<\delta$ are in $\Gamma$ and $\tau=f(\gamma,\delta)\notin M_i$, we have
\[
  \sup\{\rho(y,\tau)+1:y\in W\cap M_i\}>M_i\cap\omega_1.
\]
In particular, some $y\in W\cap M_i$ satisfies $\rho(y,\tau)\geq M_i\cap\omega_1$. The next observation lets us use any sufficiently large $y$ in its place.

\begin{claim}\label{claim: rank-tail}
If $y_0<y_1$ belong to $W\cap M_i$, $y_1<\tau$, and $\rho(y_0,\tau)\geq M_i\cap\omega_1$, then $\rho(y_1,\tau)\geq\rho(y_0,\tau)$.
\end{claim}
\begin{proof}
By subadditivity, \[\rho(y_0,\tau)\leq\max\{\rho(y_0,y_1),\rho(y_1,\tau)\}.\] By elementarity, $\rho(y_0,y_1)<M_i\cap\omega_1$.
\end{proof}

\begin{claim}
For each $i\in Z$ there are $\beta_i<\gamma_i<\delta_i$ in $\Gamma$ with $\Delta(\gamma_i,\delta_i)<\Delta(\beta_i,\gamma_i)$. Put
\[
  v_i=f(\beta_i,\gamma_i),\qquad \tau_i=f(\gamma_i,\delta_i),
  \qquad t_i=\rho(v_i,\tau_i).
\]
Then $t_i\geq M_{i'}\cap\omega_1$, where $i'=\min(Z\setminus(i+1))$, and $\lambda_{M_s}\in G_{t_i}(v_i)$ for every $s\leq i$.
\end{claim}
\begin{proof}
Choose sufficiently large $\gamma_i<\delta_i$ above $\Gamma\cap M_{i'}$, with $\tau_i>\lambda_{M_{i'}}$. This is possible by \ref{p1} applied to a tail of $\Gamma$. Take $y\in W\cap M_{i'}$ with $\rho(y,\tau_i)\geq M_{i'}\cap\omega_1$. Choose $\beta_i\in\Gamma\cap M_{i'}$ with $v_i>\max\{y,\lambda_{M_i}\}$ and $\Delta(\beta_i,\gamma_i)>\Delta(\gamma_i,\delta_i)$. Indeed, the finite initial segment of $r_{\gamma_i}$ of length $\Delta(\gamma_i,\delta_i)+1$ determines an uncountable subset of $\Gamma$ in $M_{i'}$. Its intersection with $M_{i'}$ is cofinal in $\Gamma\cap M_{i'}$, so we can use \ref{p4} and monotonicity. Since $A_{\beta_i}$ is countable and belongs to $M_{i'}$, $v_i\in W\cap M_{i'}$. Claim~\ref{claim: rank-tail} now gives $t_i\geq M_{i'}\cap\omega_1$. For $s\leq i$ both $\lambda_{M_s}$ and $v_i$ belong to $M_{i'}$, and $\lambda_{M_s}\leq \lambda_{M_i}<v_i$. Hence $\rho(\lambda_{M_s},v_i)<M_{i'}\cap\omega_1\leq t_i$, as desired.
\end{proof}

For each $i\in Z\cap\acc(Z)$, consider
\[
  d_i=\sup\{\rho(y,\lambda_{M_i}):y\in W\cap M_i\}.
\]

\begin{remark}
The consideration of $d_i$ is inspired by the proof of the same statement when the $C$-sequence on $\omega_2$ is a $\square_{\omega_1}$-sequence. In that case one can use $\sup\{\otp(C_{\lambda_{M_i}}\cap y):y\in W\cap M_i\}=\otp(C_{\lambda_{M_i}})$. For our order-minimal $C$-sequence, $\otp(C_{\lambda_{M_i}})=\omega$. Here the $\rho$-function gives a suitable replacement, as the next claim shows.
\end{remark}

We would like to recover $\lambda_{M_i}$ and $d_i$ simultaneously with a suitably chosen quadruple. If $\rho(\cdot, \lambda_{M_i})\restriction W\cap M_i$ is nondecreasing, then the task is easy. However, in general, this is not true. It turns out that the following weaker property about $\rho(\cdot, \lambda_{M_i})\restriction W\cap M_i$ suffices for our purpose.

\begin{claim}\label{claim: dominate}
We have $d_i\geq M_i\cap\omega_1$. Moreover, if $\Gamma'\in M_i$ is an uncountable subset of $\Gamma$, $b\in\Gamma\setminus M_i$, $\xi\in M_i\cap\lambda$, and $q<d_i$, there is an $a\in\Gamma'\cap M_i$ such that
\[
   f(a,b)>\xi\quad\text{and}\quad \rho(f(a,b),\lambda_{M_i})>q .
\]
\end{claim}
\begin{proof}
First take $q<M_i\cap\omega_1$. By continuity choose $l\in Z\cap i$ such that $q,\Gamma',\xi\in M_l$. By \ref{p4} and monotonicity choose a sufficiently large $a_1\in\Gamma'\cap M_i$ with $x_1=f(a_1,b)>\lambda_{M_l}$. We have $x_1\in W\cap M_i$, since $A_{a_1}\in M_i$ is countable. Since $l\in Z$ and $a_1$ was chosen sufficiently large, we can find $y\in W\cap M_l$ with $\rho(y,x_1)\geq M_l\cap\omega_1$. By \ref{p4}, choose $a_0\in\Gamma'\cap M_l$ with $x_0=f(a_0,b)>\max\{y,\xi\}$. Then $x_0<\lambda_{M_l}<x_1<\lambda_{M_i}$. By Claim~\ref{claim: rank-tail} and subadditivity,
\[
 M_l\cap\omega_1\leq\rho(x_0,x_1)\leq\max\{\rho(x_0,\lambda_{M_i}),\rho(x_1,\lambda_{M_i})\}.
\]
Since $q<M_l\cap\omega_1$, one of $a_0,a_1$ works. This also proves $d_i\geq M_i\cap\omega_1$.

If $M_i\cap\omega_1\leq q<d_i$, choose $y\in W\cap M_i$ with $\rho(y,\lambda_{M_i})>q$. By \ref{p4} and monotonicity take $a\in\Gamma'\cap M_i$ with $f(a,b)>\max\{y,\xi\}$. Apply Claim~\ref{claim: rank-tail} with $\tau=\lambda_{M_i}$.
\end{proof}

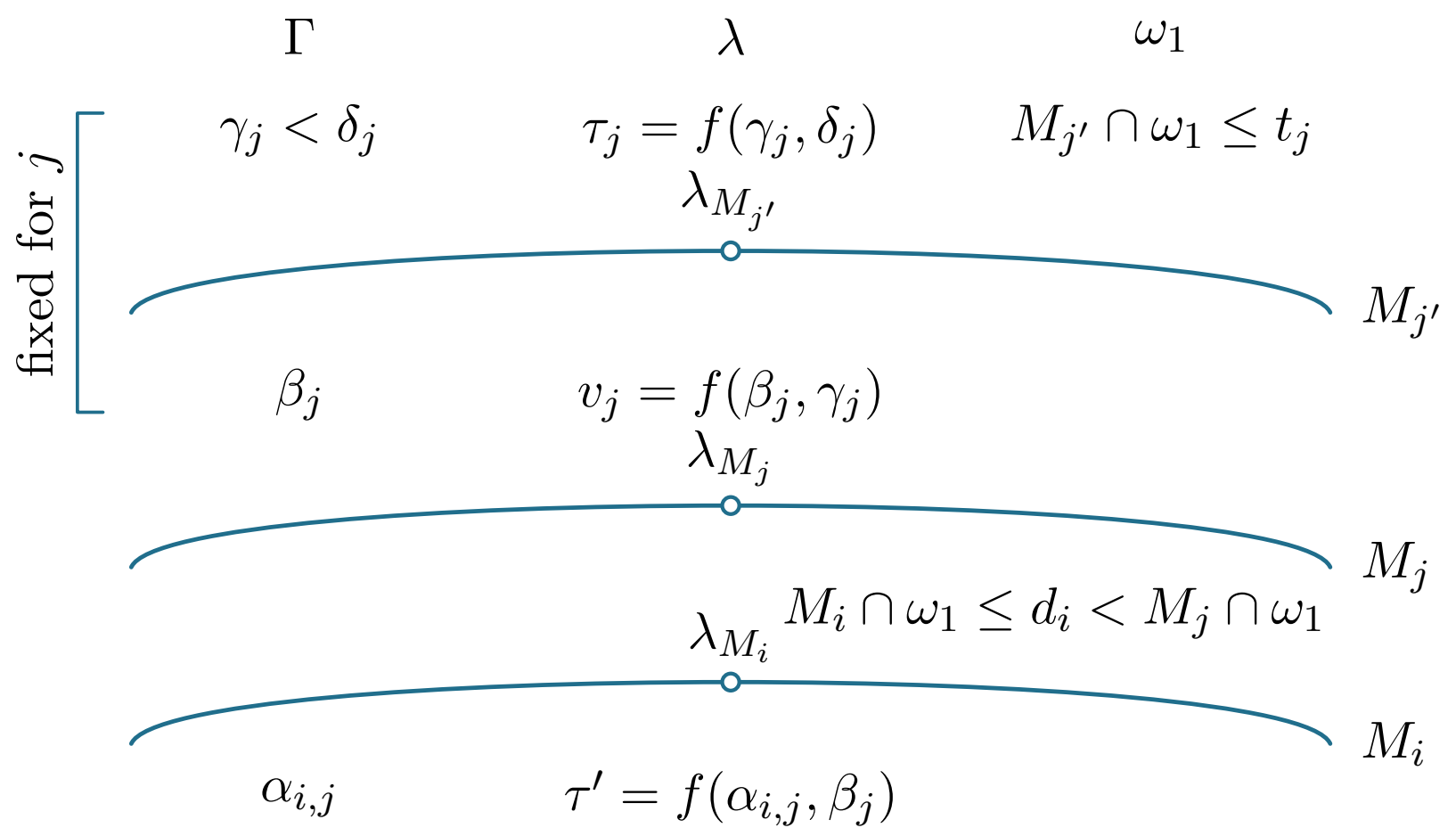
\begin{figure}[htbp]
  \centering
  \tikzexternaldisable
\begin{tikzpicture}[x=1cm,y=1cm,line cap=round,line join=round,
    font=\small,every node/.style={inner sep=2pt}]
  \definecolor{modelblue}{RGB}{32,110,140}

  \foreach \base/\index/\key in {.40/i/i,1.71/j/j,3.60/{j'}/next} {
    \draw[modelblue,line width=.9pt]
      (.40,\base) .. controls (.67,{\base+.61}) and (9.03,{\base+.61})
      .. (9.30,\base) coordinate[pos=.5] (cut-\key);
    \draw[modelblue,fill=white,line width=.8pt] (cut-\key) circle (1.8pt);
    \node[anchor=west] at (9.45,\base) {$M_{\index}$};
  }
  \node[above=2pt] at (cut-i) {$\lambda_{M_i}$};
  \node[above=2pt] at (cut-j) {$\lambda_{M_j}$};
  \node[above=2pt] at (cut-next) {$\lambda_{M_{j'}}$};

  \node at (1.65,5.65) {$\Gamma$};
  \node at (4.85,5.65) {$\lambda$};
  \node at (8.05,5.65) {$\omega_1$};

  \node at (1.65,4.95) {$\gamma_j<\delta_j$};
  \node at (1.65,2.99) {$\beta_j$};
  \node at (1.65,.00) {$\alpha_{i,j}$};
  \draw[modelblue,line width=.7pt]
    (.19,5.08) -- (0,5.08) -- (0,2.86) -- (.19,2.86);
  \node[rotate=90] at (-.28,3.97) {fixed for $j$};

  \node at (4.85,4.95) {$\tau_j=f(\gamma_j,\delta_j)$};
  \node at (4.85,2.99) {$v_j=f(\beta_j,\gamma_j)$};
  \node at (4.85,.00) {$\tau'=f(\alpha_{i,j},\beta_j)$};

  \node at (8.05,4.95) {$M_{j'}\cap\omega_1\leq t_j$};
  \node at (7.25,1.37) {$M_i\cap\omega_1\leq d_i<M_j\cap\omega_1$};
\end{tikzpicture}
  \caption{The models in \textbf{(3b)}, with $i<j$ and $j'=\min\{k\in Z:k>j\}$. The circles mark the cuts $\lambda_{M_k}=\sup(M_k\cap\lambda)$. The same $\beta_j,\gamma_j,\delta_j$ work for every $i\in Z\cap\acc(Z)$ below $j$. Choosing $\alpha_{i,j}\in M_i\cap\Gamma$ gives $c(\alpha_{i,j},\beta_j,\gamma_j,\delta_j)=g(d_i,t_j)$.}
  \label{fig:continuum-fixed-triple}
\end{figure}

We now choose $i<j$ in $Z\cap\acc(Z)$ with $g(d_i,t_j)=n$. Since $\lambda_{M_i}\in G_{t_j}(v_j)$ and $d_i<t_j$, choose $m_0\in\omega$ with
\[
 \lambda_{M_i}\in K^{v_j,t_j}_{m_0},\qquad d_i\in F_{m_0}(t_j).
\]
Choose $\sigma\in S$ with $\Delta(\sigma,r_{\beta_j})=m\geq \max\{m_0,\Delta(\beta_j,\gamma_j)\}+1$. Put
\[
 \xi_0=\sup(K_m^{v_j,t_j}\cap \lambda_{M_i}),\qquad
 q=\sup(F_m(t_j)\cap d_i).
\]
These bounds are strictly below $\lambda_{M_i}$ and $d_i$, respectively, by finiteness. Apply Claim~\ref{claim: dominate} to $b=\beta_j$ and the uncountable set $\{a\in\Gamma:\sigma\sqsubseteq r_a\}\in M_i$. We obtain $\alpha\in M_i\cap\Gamma$ such that $\sigma\sqsubseteq r_\alpha$ and, with $\tau'=f(\alpha,\beta_j)$,
\[
  \xi_0<\tau'<\lambda_{M_i},\qquad q<\rho(\tau',\lambda_{M_i})\leq d_i.
\]
Here $\tau'\in W\cap M_i$ by countability of $A_\alpha$. Thus
\[
 \Delta(\gamma_j,\delta_j)<\Delta(\beta_j,\gamma_j)
       <\Delta(\alpha,\beta_j)=m .
\]
We thus enter case \textbf{(3b)}. The upper rank is $\rho(f(\beta_j,\gamma_j),f(\gamma_j,\delta_j))=t_j$. Using the two bounds just obtained, we have
\[
 \min(K_m^{v_j,t_j}\setminus\tau')=\lambda_{M_i},\qquad
 \min(F_m(t_j)\setminus\rho(\tau',\lambda_{M_i}))=d_i .
\]
Consequently $c(\alpha,\beta_j,\gamma_j,\delta_j)=g(d_i,t_j)=n$.

\subsection{Square-bracket operation (3c)}\label{subsection: continuum(3c)}

We first record a consequence of excluding \textbf{(3a)}. For a countable elementary submodel $X\prec H((2^\theta)^+)$ containing the relevant objects, put $\lambda_X=\sup(X\cap\lambda)$.

\begin{claim}\label{claim: big}
Let $\tau\in W$ with $\lambda_X<\tau$. For every $\xi<X\cap\omega_1$, $\rho(y,\tau)>\xi$ for all sufficiently large $y<\lambda_X$. In particular, if $b\in\Gamma-X$ and $\Gamma'\in X$ is an uncountable subset of $\Gamma$, then there is $a\in\Gamma'\cap X$ such that $\rho(f(a,b),\tau)>\xi$.
\end{claim}

\begin{proof}
Write $\bar\tau=\last{\lambda_X}{\tau}$. Since $\tau\in W$, it has the form $f(\gamma,\delta)$, so the failure of \textbf{(3a)} gives $\otp(C_{\bar\tau}\cap\lambda_X)\geq X\cap\omega_1$. Now $\cf(\lambda_X)=\omega$ and $\otp(C_{\lambda_X})=\omega<X\cap\omega_1$. Hence $\bar\tau>\lambda_X$.

On the walk from $\tau$ to $\lambda_X$, each node $v$ before $\bar\tau$ satisfies $\sup(C_v\cap\lambda_X)<\lambda_X$. Choose $y<\lambda_X$ above all these finitely many bounds. At each such $v$ we have $\min(C_v\setminus y)=\min(C_v\setminus\lambda_X)$, so the walk from $\tau$ to $y$ also passes through $\bar\tau$. By Lemma~\ref{lem:walk-rank},
\[
\rho(y,\tau)\geq\rho(y,\bar\tau)\geq\otp(C_{\bar\tau}\cap y).
\]
As $y$ tends to $\lambda_X$, the order types on the right have supremum $\otp(C_{\bar\tau}\cap\lambda_X)\geq X\cap\omega_1$. This gives the first assertion. For the second, \ref{b3} puts $b^*$ above every member of $\Gamma\cap X$, hence outside $X$ since $\Gamma$ is cofinal in $\varepsilon$. By \ref{p4}, the values $f(a,b)$ for $a\in\Gamma\cap X$ are cofinal in $\lambda_X$. The same holds for $a\in\Gamma'\cap X$ by monotonicity, since $\Gamma'\cap X$ is cofinal in $\Gamma\cap X$.
\end{proof}

We now use the failure of \textbf{(3b)}. For club many $i<\omega_1$, every tail of $\Gamma$ contains $\gamma<\delta$ such that $\tau=f(\gamma,\delta)\notin M_i$ and
\[
\sup_{y\in W\cap M_i}\bigl(\rho(y,\tau)+1\bigr)\leq M_i\cap\omega_1.
\]
We may choose $i$ in this club and a countable $N\prec H((2^\theta)^+)$ containing all the relevant objects so that $N\cap\omega_1=M_i\cap\omega_1=i$ and $d(N\cap\omega_1)=n$. In particular $\lambda_N=\lambda_{M_i}$.

Choose a countable elementary submodel $N_0\prec H((2^\theta)^+)$ containing the relevant objects with $N\in N_0$. Take a large enough witnessing pair $\gamma<\delta$ on the tail above $\sup(\Gamma\cap N_0)$ and put $\tau=f(\gamma,\delta)$. Since $\delta^*\notin N_0$, \ref{p4} and monotonicity give
\[
\tau > \sup\{f(a,\delta):a\in\Gamma\cap N_0\}=\lambda_{N_0}>\lambda_N.
\]
Claim~\ref{claim: big} and the cofinality of $W\cap N$ in $\lambda_N$ now give
\[
\sup_{y\in W\cap N}\bigl(\rho(y,\tau)+1\bigr)=\sup_{y\in W\cap M_i}\bigl(\rho(y,\tau)+1\bigr)=M_i\cap \omega_1=N\cap\omega_1.
\]

Put $q=\Delta(\gamma,\delta)+1$. The set of $\beta\in\Gamma$ for which $r_\delta\restriction q\sqsubseteq r_\beta$ is uncountable. Indeed, it belongs to $N_0$ and contains $\delta\notin N_0$. Removing its members at or below $\sup(\Gamma\cap N)$ leaves an uncountable set in $N_0$. Apply Claim~\ref{claim: big} in $N_0$ to this set, with $b=\gamma$ and $\xi=N\cap\omega_1<N_0\cap\omega_1$. We obtain $\beta\in\Gamma\cap N_0$ above $\sup(\Gamma\cap N)$ such that
\[
\Delta(\gamma,\delta)<\Delta(\beta,\delta)
\quad\hbox{and}\quad
\nu=\rho(f(\beta,\gamma),\tau)>N\cap\omega_1.
\]

Put $m_0=e_\nu(N\cap\omega_1)$. Every finite initial segment of $r_\delta$ lies in $S$, by the same countability argument. Since $S$ is perfect, choose $\sigma\in S$ incompatible with $r_\delta$ so that their first disagreement $m=\Delta(\sigma,r_\delta)$ satisfies $m\geq m_0$ and $m>\Delta(\beta,\delta)$. The set $\Gamma'=\{a\in\Gamma:\sigma\sqsubseteq r_a\}$ is uncountable and belongs to $N$. Apply Claim~\ref{claim: big} again, now in $N$ with $b=\delta$, to choose $\alpha\in\Gamma'\cap N$ with
\[
\sup(F_m(\nu)\cap(N\cap\omega_1))<\rho(f(\alpha,\delta),\tau)<N\cap\omega_1.
\]
For the upper bound, $f(\alpha,\delta)\in A_\alpha\subseteq W\cap N$ since $A_\alpha$ is countable and belongs to $N$, so it follows from the hypothesis.

We have $\alpha<\beta<\gamma<\delta$ and $\Delta(\gamma,\delta)<\Delta(\beta,\delta)<\Delta(\alpha,\delta)=m$. Since $m\geq m_0$, $N\cap\omega_1\in F_m(\nu)$. Entering case \textbf{(3c)}, we get
\[
c(\alpha,\beta,\gamma,\delta)
=d\bigl(\min(F_m(\nu)\setminus\rho(f(\alpha,\delta),\tau))\bigr)
=d(N\cap\omega_1)=n.
\]

\begin{proof}[Proof of Theorem~\ref{theorem: main}]
If $2^{\omega_1}=\omega_2$, then a classical result of Specker gives an $\omega_2$-special $\omega_3$-Aronszajn tree, and Theorem~\ref{theorem: special} applies. If $2^{\omega_1}\geq\omega_3$ and $2^\omega<\omega_3$, then a theorem of Baumgartner \cite[Theorem 3.5]{BaumgartnerAlmostDisjoint} gives a weak Kurepa tree with at least $\omega_3$ uncountable branches, and Theorem~\ref{theorem: wkt} applies. Finally, if $\omega_3\leq 2^\omega$, apply Theorem~\ref{theorem: continuum}.
\end{proof}

\section{Some remarks}\label{section: cor}

%

It is unknown whether $\omega_2\not\to [\omega_1]^2_{\omega}$ is equivalent to the failure of $\mathrm{CH}$ (\cite[Question 10.3.9]{walksbook}). See \cite{shelahPositive} for some recent progress in the forcing direction. A positive answer would give a short proof of $\omega_3\not\to [\omega_1]^4_{\omega}$ by considering separately the cases of $\mathrm{CH}$ and $\neg\mathrm{CH}$.

\begin{corollary}\label{cor: omega_2to3}
If $\omega_2\not\to [\omega_1]^2_{\omega}$, then $\omega_3 \not\to [\omega_1]^4_\omega$.
\end{corollary}

\begin{proof}
The idea of the proof already appears in the previous sections, so we give a shorter argument here.

Fix $g: [\omega_2]^2\to \omega$ witnessing $\omega_2\not\to [\omega_1]^2_{\omega}$ and $h_2: [\omega_1]^2\to \omega$ witnessing $\omega_1\not\to [\omega_1]^2_\omega$.
We define $c: [\omega_3]^4\to \omega$ as follows.
Fix $\alpha<\beta<\gamma<\delta<\omega_3$.
Let $k=\chi(\beta,\gamma,\delta)$. Let $\delta^* = \Tr(\alpha,\delta)(k)$, $\gamma^*=\Tr(\alpha,\gamma)(k)$ and $\beta^*=\Tr(\alpha,\beta)(k)$.
We define $c(\alpha,\beta,\gamma,\delta)$ by the following cases.
	\begin{enumerate}
     \item If $(C_{\Tr(\beta,\gamma)(k)}-\beta) \cap (C_{\Tr(\beta,\delta)(k)}-\beta) = \emptyset$, then output $$osc(C_{\Tr(\beta,\gamma)(k)}-\beta,C_{\Tr(\beta,\delta)(k)}-\beta)-1.$$ 
     \item Else if $\beta^*, \gamma^*, \delta^*\in \cof(\omega_1)$, output $h_2(\otp(C_{\beta^*}\cap \alpha), \otp(C_{\delta^*}\cap \gamma))$.
     \item Else if $\beta^*, \gamma^*, \delta^*\in \cof(\omega_2)$, output $g(\otp(C_{\beta^*}\cap \alpha), \otp(C_{\delta^*}\cap \gamma))$.
     \item Else output $0$.
	\end{enumerate}
Fix $\Gamma\in [\omega_3]^{\omega_1}$. For case \textbf{(1)}, apply the argument for case \textbf{(1)} of Theorem \ref{theorem: omega2} to find $\beta<\gamma<\delta$ in $\Gamma\setminus\{\min\Gamma\}$ and take $\alpha=\min\Gamma$. Hence we may assume that there are a club $D\subseteq\sup\Gamma$ of order type $\omega_1$ and $k\in\omega$ such that
\begin{itemize}
\item for any $\nu\in\acc(D)$, there is some $\eta_\nu\in\Gamma$ such that, letting $\bar{\eta}_\nu=\Tr(\nu,\eta_\nu)(k)$, we have $D\cap\nu\subseteq C_{\bar{\eta}_\nu}$,
\item for any $\nu_0<\nu_1<\nu_2\in \acc(D)$, $\chi(\eta_{\nu_0}, \eta_{\nu_1}, \eta_{\nu_2})=k$, and 
\item for any $\nu_0<\nu_1\in \acc(D)$, $\Tr(\eta_{\nu_0}, \eta_{\nu_1})(k)=\bar{\eta}_{\nu_1}$.
\end{itemize}

Suppose there is an uncountable $A\subseteq\acc(D)$ such that $\bar{\eta}_\nu\in\cof(\omega_1)$ for every $\nu\in A$. Then $\{\otp(C_{\bar{\eta}_{\nu_1}}\cap\eta_{\nu_0}):\nu_0<\nu_1\in A\}$ is an uncountable subset of $\omega_1$. We can therefore step up $h_2$ using case \textbf{(2)} of the definition.

Otherwise, by passing to a tail of $D$, we may assume that $\bar{\eta}_\nu\in\cof(\omega_2)$ for every $\nu\in\acc(D)$. Then $\{\otp(C_{\bar{\eta}_{\nu_1}}\cap\eta_{\nu_0}):\nu_0<\nu_1\in\acc(D)\}$ is an uncountable subset of $\omega_2$, and we can step up $g$ using case \textbf{(3)} of the definition.

\end{proof}

The following corollary of the proof of Proposition \ref{prop: deltaprime} appeared in \cite{TodorcevicZhangOsc}.

\begin{corollary}
For any cardinal $\kappa\geq \omega_1$, $\kappa\not\to [\mathrm{stat}(\omega_1)]^2_{\omega}$. Namely, there exists a coloring $c: [\kappa]^2\to \omega$ such that for any $A\subseteq \kappa$ of order type $\omega_1$ that is stationary in $\sup A$, $c''[A]^2 =\omega$.
\end{corollary}

In contrast, the following is consistent relative to the existence of a huge cardinal (\cite[Lemma 2.20]{GartiZhang}): for any $c: [\omega_2]^{<\omega_1}\to \omega_1$, there exists $A\subseteq \omega_2$ of order type $\omega_1$ stationary in $\sup A$ such that $c''[A]^{<\omega}$ is countable.

\section{Weak Kurepa trees and versions of Chang's conjecture}\label{section: WKTandCC}

The technique employed in Section \ref{section: weakKurepa} is suitable for connecting the existence of strong colorings of various dimensions to the failure of versions of Chang's conjecture. The unconditional result for dimension 3 was obtained by Todorcevic \cite{Todorcevic3D}, who showed in ZFC that Chang's conjecture fails if and only if $\omega_2\not\to [\omega_1]^3_{\omega_1}$.

\begin{theorem}\label{theorem: WKTdimensionN}
Suppose $n\geq 3$ and there exists a weak Kurepa tree $T$ with at least $\kappa$ branches. Then (2) implies (1) for the following statements.
	\begin{enumerate}
	\item $\kappa\to (\omega_1)^n_{\omega_1, \omega}$, namely, for any $f: [\kappa]^{n}\to \omega_1$, there is $A\in [\kappa]^{\omega_1}$ such that $f''[A]^{n}$ is countable,
	\item $\kappa \to [\omega_1]^{n+1}_{\omega_1}$. 
	\end{enumerate}
\end{theorem}

\begin{proof}
We largely follow the notation from the proof of Theorem \ref{theorem: wkt}. We assume that (1) fails and show that (2) fails. Fix the following data.

	\begin{enumerate}
	\item We may assume $T\subseteq 2^{<{\omega_1}}$, ordered by end extension.
	\item Let $\langle b_\alpha:\alpha<\kappa\rangle$ be a sequence of distinct uncountable branches through $T$. As before, write $\Delta(\alpha,\beta)=\Delta(b_\alpha,b_\beta)$.
	\item Let $g: [\kappa]^n\to \omega_1$ witness the failure of $\kappa\to (\omega_1)^n_{\omega_1, \omega}$. Namely, for any $A\in [\kappa]^{\omega_1}$, $g''[A]^n$ is uncountable.
	\item Fix $h: [\omega_1]^2\to \omega_1$ witnessing $\omega_1\not\to [\omega_1, \omega_1]^2_{\omega_1}$. Namely, for any $A,B\in [\omega_1]^{\omega_1}$, $h'' \{(\alpha,\beta): \alpha<\beta, \alpha\in A, \beta\in B\}=\omega_1$.
	\item Fix $d: \omega_1\to \omega_1$ such that each fiber of $d$ is stationary in $\omega_1$.
	\item Let $\operatorname{index}_T: T\to \omega_1$ be a bijection.
	\item For each $\nu<\omega_1$, fix an injection $e_\nu:\nu\to\omega$. For $m\in\omega$, let $F_m(\nu)=\{\xi<\nu:e_\nu(\xi)\leq m\}$.
	\end{enumerate}

Let us define $c:[\kappa]^{n+1}\to\omega_1$. Given $a=\{a_0<\cdots<a_n\}$, let $a'=\{a_0,\ldots,a_{n-2}\}$ and $\delta=\Delta(a_{n-1},a_n)$. Put
\[x=g(a'\cup\{a_{n-1}\}),\qquad y=g(a'\cup\{a_n\}).\]
    \begin{enumerate}
	\item If $\max\{\delta,x\}<y$, output $h(x,y)$.
	\item Else if $\max\{x,y\}<\delta$ and $\Delta(a_{n-3},a_n)<\Delta(a_{n-2},a_n)<\delta$, output $h(\Delta(a_{n-3},a_n),\Delta(a_{n-2},a_n))$.
	\item Else if $\max\{x,y\}<\delta$ and $\Delta(a_{n-2},a_n)\leq\Delta(a_{n-3},a_n)<\delta$, let $m=e_\delta(\Delta(a_{n-2},a_n))$ and output
	\[
	d\bigl(\min(F_m(\delta)\setminus \Delta(a_{n-3},a_n))\bigr).
	\]
	\item Else if $\delta<x=y$, let $\nu=x=y$, $\xi_0=\operatorname{index}_T(b_{a_{n-1}}\restriction\nu)$ and $\xi_1=\operatorname{index}_T(b_{a_n}\restriction\nu)$. For $m=e_{\xi_1}(\delta)$, output
	\[
	d\bigl(\min(F_m(\xi_1)\setminus \xi_0)\bigr).
	\]
	\item If no case applies or the indicated expression is undefined, output $0$.
	\end{enumerate}
	
Let us verify that the coloring works. Fix $\Gamma\in[\kappa]^{\omega_1}$ and $\tau<\omega_1$. We may take $\Gamma$ of order type $\omega_1$, since it suffices to find witnesses in an uncountable subset. We consider the following cases.

	\textbf{(1):} Suppose there is some $b\in[\Gamma]^{n-1}$ such that $E=\{g(b\cup\{\delta\}):\delta\in\Gamma,\ \delta>\max b\}$ is uncountable. By shrinking the part of $\Gamma$ above $\max b$, we may assume that $\delta\mapsto g(b\cup\{\delta\})$ is strictly increasing there. For each $s\in T$, let $E_s=\{g(b\cup\{\delta\}):\delta\in\Gamma,\ \delta>\max b,\ s\sqsubseteq b_\delta\}$.
		\begin{itemize}
		\item If there are $s,t \in T$ with $s\perp t$ such that $E_s$ and $E_t$ are both uncountable, then by the property of $h$, we can find $\delta_0<\delta_1\in \Gamma- (\max b+1)$ such that $s\sqsubseteq b_{\delta_0}$, $t\sqsubseteq b_{\delta_1}$, $\Delta(s,t)<g(b\cup \{\delta_0\}) < g(b\cup \{\delta_1\})$ and $h(g(b\cup \{\delta_0\}), g(b\cup \{\delta_1\}))=\tau$. Then by case \textbf{(1)} of the definition, we have that $c(b\cup \{\delta_0, \delta_1\})=\tau$.
		\item Otherwise we have an uncountable branch $b^*$ through $T$ such that if $E_s$ is uncountable, then $s\sqsubseteq b^*$. We may find $\Sigma\in [\Gamma]^{\omega_1}$ such that for each $\alpha\in \Sigma$, there is some $\xi_\alpha\in \omega_1$ such that $\Delta(b_\alpha,b^*)=\xi_\alpha$ and $\langle \xi_\alpha: \alpha\in \Sigma\rangle$ is strictly increasing. By Lemma \ref{lemma: strengthenedMoore}, we can find $\alpha\in \Sigma$ such that there are uncountably many $\beta\in \Sigma$ such that $h(g(b\cup \{\alpha\}), g(b\cup \{\beta\}))=\tau$. Choose such a $\beta>\alpha$ with $g(b\cup \{\beta\}) > \xi_\alpha$. Since $\max\{\Delta(\alpha,\beta), g(b\cup \{\alpha\})\}<g(b\cup \{\beta\})$, case \textbf{(1)} applies and we get $c(b\cup \{\alpha,\beta\})=h(g(b\cup \{\alpha\}), g(b\cup \{\beta\}))=\tau$.
		\end{itemize}

Suppose now that $A_b=\{g(b\cup\{\delta\}):\delta\in\Gamma,\ \delta>\max b\}$ is countable for every $b\in[\Gamma]^{n-1}$. By shrinking $\Gamma$ to an uncountable subset, we may assume that each $i\in A_b$ is realized by uncountably many $\delta\in\Gamma$ above $\max b$. Let $S$ be the subtree of $T$ consisting of those $t$ for which $\{\delta\in\Gamma:t\sqsubseteq b_\delta\}$ is uncountable. The remaining cases depend on $S$.

	\textbf{(2):} If $S$ has an uncountable branch, then, shrinking $\Gamma$ if necessary, we may assume that there is a strictly increasing sequence $\langle\xi_\alpha:\alpha\in\Gamma\rangle$ such that $\Delta(\alpha,\beta)=\xi_\alpha$ whenever $\alpha<\beta$ are in $\Gamma$. By the property of $h$, find $\alpha<\beta$ in $\Gamma$ large enough that $h(\xi_\alpha,\xi_\beta)=\tau$ and choose $a'\in[\Gamma]^{n-1}$ whose last two members are $\alpha,\beta$. Since $A_{a'}$ is countable, find $\beta<a_{n-1}<a_n$ in $\Gamma$ such that $\xi_{a_{n-1}}=\Delta(a_{n-1},a_n)>\sup A_{a'}+1$ and $g(a'\cup\{a_{n-1}\})=g(a'\cup\{a_n\})$. By case \textbf{(2)} of the definition, $c(a'\cup\{a_{n-1},a_n\})=h(\xi_\alpha,\xi_\beta)=\tau$.

	\textbf{(3):} Suppose now that $S$ is an Aronszajn tree. Let $\theta$ be a sufficiently large regular cardinal and choose a countable elementary submodel $N\prec H(\theta)$ containing the relevant objects such that $d(N\cap\omega_1)=\tau$. Since $S$ is normal and perfect of height $\omega_1$, we can find $a_{n-1}<a_n$ in $\Gamma-N$ and $t\in T_{N\cap\omega_1}\cap S$ such that $t\sqsubseteq b_{a_n}$ and $\delta=\Delta(a_{n-1},a_n)>N\cap\omega_1$. Let $m_0=e_\delta(N\cap\omega_1)$. Since $t\notin N$, we can find $s_0\in S\cap N$ such that $\Delta(s_0,t)=\xi_0$ and $m=e_\delta(\xi_0)\geq m_0$. Repeating the same procedure, find $s_1\in S\cap N$ such that $\Delta(s_1,t)=\xi_1>\max\{\sup(F_m(\delta)\cap N),\xi_0\}+1$.
	
	Choose $a_{n-3}<a_{n-2}$ in $\Gamma\cap N$ such that $s_1\sqsubseteq b_{a_{n-3}}$ and $s_0\sqsubseteq b_{a_{n-2}}$, with room for $n-3$ earlier members of $\Gamma\cap N$. Complete these two points to $a'\in[\Gamma\cap N]^{n-1}$ by choosing those earlier members. Since $A_{a'}$ is a countable set in $N$, we have
	\[
	\max\{g(a'\cup\{a_{n-1}\}),g(a'\cup\{a_n\})\}<N\cap\omega_1<\delta.
	\]
	Moreover,
	\[
	\Delta(a_{n-2},a_n)=\xi_0<\Delta(a_{n-3},a_n)=\xi_1<N\cap\omega_1<\delta.
	\]
Thus case \textbf{(3)} of the definition applies, and
\[
c(a'\cup\{a_{n-1},a_n\})=d\bigl(\min(F_m(\delta)\setminus \xi_1)\bigr)=d(N\cap\omega_1)=\tau.
\]

		\textbf{(4):} If $S$ is not an Aronszajn tree, by shrinking $\Gamma$ to an uncountable subset, we may assume that $S$ satisfies \ref{cantor1} and \ref{cantor2}. Let $\eta=\operatorname{ht}(S)\in \acc(\omega_1)\cap \omega_1$.

Let $N\prec H(\theta)$ be a countable elementary submodel containing the relevant objects such that $d(N\cap\omega_1)=\tau$. Since $g''[\Gamma]^n$ is unbounded in $\omega_1$, by elementarity we can find $a'\in[\Gamma]^{n-1}\cap N$ and $\eta<\nu<N\cap\omega_1$ such that $g(a'\cup\{\delta\})=\nu$ for uncountably many $\delta\in\Gamma$ above $\max a'$. Among these, choose $a_n\in\Gamma-N$ such that $\xi_1=\operatorname{index}_T(b_{a_n}\restriction\nu)>N\cap\omega_1$ and $b_{a_n}\restriction\eta\notin N$. We say that $b_{a_n}\restriction\eta$ is \emph{fresh over $N$} because all its proper restrictions belong to $S\subseteq N$.
Let $m_0=e_{\xi_1}(N\cap\omega_1)$. Find $s\in S\cap N$ such that $m=e_{\xi_1}(\Delta(s,b_{a_n}))\geq m_0$ and there are uncountably many $\delta\in\Gamma$ with $s\sqsubseteq b_\delta$ and $g(a'\cup\{\delta\})=\nu$. Indeed, $F_{m_0}(\xi_1)\cap\eta$ is finite and hence bounded in $\eta$. If there were no such $s$, we could define $b_{a_n}\restriction\eta$ in $N$, contradicting the fact that it is fresh over $N$.

Finally, by the property of $s$, we can find $a_{n-1}\in\Gamma\cap N$ above $\max a'$ such that $s\sqsubseteq b_{a_{n-1}}$, $\xi_0=\operatorname{index}_T(b_{a_{n-1}}\restriction\nu)>\sup(F_m(\xi_1)\cap N)$ and $g(a'\cup\{a_{n-1}\})=\nu$.
Note that $\xi_0\in N$ since $a_{n-1}$ and $\nu$ both belong to $N$.

As $g(a'\cup \{a_{n-1}\})=g(a'\cup \{a_{n}\})>\eta>\Delta(a_{n-1}, a_n)$, we enter case \textbf{(4)} of the definition. 
Thus
\begin{align*}
c(a'\cup\{a_{n-1},a_n\})
 &=d\bigl(\min(F_{e_{\xi_1}(\Delta(a_{n-1},a_n))}(\xi_1)-\xi_0)\bigr)\\
 &=d\bigl(\min(F_m(\xi_1)\setminus \xi_0)\bigr)=d(N\cap\omega_1)=\tau.
\end{align*}

\end{proof}

When $n=2$, we obtain the conclusion of Theorem \ref{theorem: WKTdimensionN} under the additional hypothesis that the continuum is large. This hypothesis lets us exploit the difference between a Cantor tree and an Aronszajn tree.

\begin{theorem}\label{theorem: weakKurepaDimension2}
If there is a weak Kurepa tree with at least $\kappa$ branches and $\kappa\leq 2^\omega$, then $\kappa\not\to (\omega_1)^2_{\omega_1, <\omega_1}$ implies that $\kappa\not\to [\omega_1]^3_{\omega_1}$.
\end{theorem}

\begin{proof}
As before, we may assume $T\subseteq 2^{<\omega_1}$. Fix the following data.
	\begin{enumerate}
	\item Choose distinct reals $\langle r_\alpha:\alpha<\kappa\rangle$ in $2^\omega$ and distinct branches $\langle b_\alpha:\alpha<\kappa\rangle$ in $[T]$. Write $\Delta(\alpha,\beta)=\Delta(b_\alpha,b_\beta)$.
	\item Fix $d:\omega_1\to\omega_1$ such that each fiber of $d$ is stationary.
	\item For distinct reals $s,t$, let $\Delta_0(s,t)=\min\{k\in\omega:s(k)\neq t(k)\}$. Use the same notation for the first disagreement with an incompatible finite sequence. Write $\Delta_0(\alpha,\beta)=\Delta_0(r_\alpha,r_\beta)$.
	\item Fix $g: [\kappa]^2\to \omega_1$ witnessing that $\kappa\not\to (\omega_1)^2_{\omega_1, <\omega_1}$ and $h: [\omega_1]^2\to \omega_1$ witnessing that $\omega_1\not\to [\omega_1, \omega_1]^2_{\omega_1}$.
	\item For each $\nu<\omega_1$, fix an injection $e_\nu:\nu\to\omega$. For $m\in\omega$, let $F_m(\nu)=\{\xi<\nu:e_\nu(\xi)\leq m\}$.
	\item Let $\operatorname{index}_T: T\to \omega_1$ be a bijection.
	\end{enumerate}
The definition differs from that in Theorem \ref{theorem: WKTdimensionN} only in the Aronszajn tree case.

Given $\alpha<\beta<\gamma<\kappa$, define $c:[\kappa]^3\to\omega_1$ by the following cases.
\begin{enumerate}
	\item If $\max\{\Delta(\beta,\gamma), g(\alpha,\beta)\}<g(\alpha,\gamma)$, output $h(g(\alpha,\beta),g(\alpha,\gamma))$.
	\item Else if $\max\{g(\alpha,\beta),g(\alpha,\gamma)\}<\Delta(\beta,\gamma)$ and $\Delta(\alpha,\beta)<\Delta(\beta,\gamma)$, let $m=\Delta_0(\alpha,\gamma)$ and output
	\[
	d\bigl(\min(F_m(\Delta(\beta,\gamma))-\Delta(\alpha,\gamma))\bigr).
	\]
	\item Else let $\nu=g(\alpha,\gamma)$, $\xi_0=\operatorname{index}_T(b_\beta\restriction\nu)$ and $\xi_1=\operatorname{index}_T(b_\gamma\restriction\nu)$. For $m=e_{\xi_1}(\Delta(\beta,\gamma))$, output
	\[
	d\bigl(\min(F_m(\xi_1)\setminus \xi_0)\bigr).
	\]
	\item If the indicated expression is undefined, output $0$.
	\end{enumerate}
	
	\begin{figure}
    \centering 
   \includegraphics[scale=0.6]{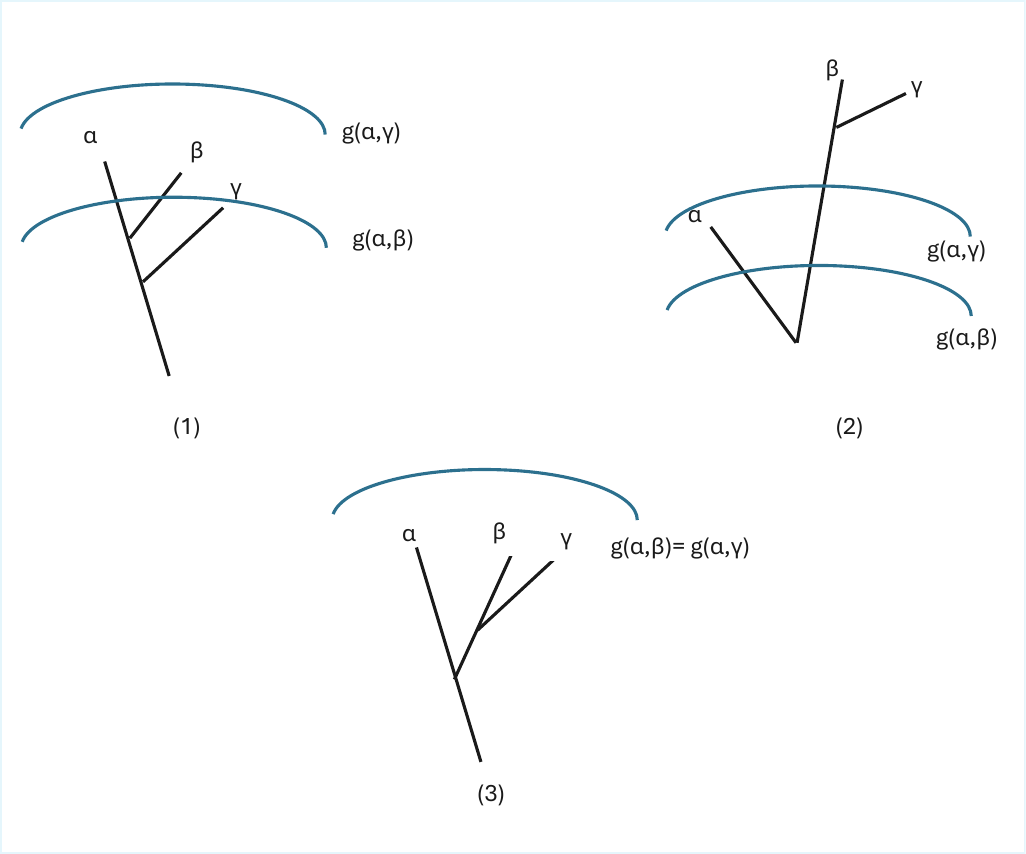}
    \caption{Tuple configurations based on a weak Kurepa tree} 
    \label{fig:KurepaContinuum} 
\end{figure}

Let us show that this coloring works. Fix $\Gamma\in [\kappa]^{\omega_1}$ and $\tau\in \omega_1$.

\textbf{(1):} If there is some $\alpha\in \Gamma$ such that $\{g(\alpha,\beta): \beta\in \Gamma,\ \beta>\alpha\}$ is uncountable, then we proceed exactly as in case \textbf{(1)} in the proof of Theorem \ref{theorem: WKTdimensionN}.

Suppose now that $A_\alpha=\{g(\alpha,\beta):\beta\in\Gamma,\ \beta>\alpha\}$ is countable for each $\alpha\in\Gamma$. By shrinking $\Gamma$ if necessary, we may assume that each $\nu\in A_\alpha$ is realized by uncountably many $\beta\in\Gamma$ above $\alpha$.

Consider the two trees generated by $\Gamma$, namely $S\subseteq 2^{<\omega}$ and $T_0\subseteq T$. More precisely, $\sigma\in S$ if and only if $\{\alpha\in \Gamma: \sigma\sqsubseteq r_\alpha\}$ is uncountable and $t\in T_0$ if and only if $\{\alpha\in \Gamma: t\sqsubseteq b_\alpha\}$ is uncountable.

\textbf{(2):} Suppose that $T_0$ either is an Aronszajn tree or contains an uncountable branch $b^*$.
Let $N\prec H(\theta)$ be a countable elementary submodel containing the relevant objects such that $N\cap \omega_1=\delta\in d^{-1}(\{\tau\})$. 

Call $\gamma\in\Gamma-N$ \emph{good} if $b_\gamma\restriction(\delta+1)\in T_0$ and for infinitely many $m\in\omega$ there is a $t\in S$ such that $\Delta_0(r_\gamma,t)=m$ and
\[
\{\Delta(\alpha,\gamma):\alpha\in\Gamma\cap N,\ t\sqsubseteq r_\alpha\}\cap\delta
\quad\text{is unbounded in }\delta.
\]
 
\begin{claim}
There exists a good $\gamma$.
\end{claim} 

\begin{proof}[Proof of the Claim]

Suppose first that $T_0$ contains an uncountable branch $b^*$. Let $\gamma\in \Gamma-N$. We may assume that $\langle \Delta(b_\gamma, b^*): \gamma\in \Gamma\rangle$ is an increasing sequence unbounded in $\omega_1$. For each $\sigma\in 2^{<\omega}$, consider the statement $\varphi(\sigma)$: there exist uncountably many $\gamma\in \Gamma$ such that $\sigma\sqsubseteq r_\gamma$. Since every finite initial segment $t$ of $r_\gamma$ satisfies $\varphi(t)$ and $r_\gamma\not \in N$, we know that there are unboundedly many $m\in \omega$ for which there is $t\in S$ such that $\Delta_0(t, r_\gamma)=m$ and $\varphi(t)$. For each such $t$, elementarity of $N$ gives, for every $\xi\in \delta$, some $\alpha\in N\cap \Gamma$ such that $t\sqsubseteq r_\alpha$ and $\Delta(b_\alpha, b^*)=\Delta(\alpha,\gamma)>\xi$, exactly as desired.

Suppose now that $T_0$ is an Aronszajn tree. We first show that any $\gamma\in\Gamma-N$ with $b_\gamma\restriction(\delta+1)\in T_0$ is \emph{pre-good} in the following sense. For every $k\in\omega$, there are unboundedly many $\xi<\delta$
such that some $\alpha\in\Gamma\cap N$ satisfies $\Delta(\alpha,\gamma)=\xi$ and $\Delta_0(\alpha,\gamma)>k$. Suppose otherwise and fix $k\in\omega$ and $\xi<\delta$ witnessing the failure.
Let $\sigma=r_\gamma\restriction(k+1)\in N$ and $\Gamma_0=\{a\in\Gamma:\sigma\sqsubseteq r_a\}\in N$. Since $\gamma\in\Gamma_0-N$, the set $\Gamma_0$ is uncountable. The tree generated by $\Gamma_0$ is a subtree of the Aronszajn tree $T_0$ and has uncountable height. Thus we can find some $\alpha\in\Gamma_0\cap N$ with $\Delta(\alpha,\gamma)>\xi$, contradicting our choice of $k$ and $\xi$.

Take an $\in$-chain $N=N_0\in N_1\in\cdots$ of countable elementary submodels. Let $M=\bigcup_{i\in\omega}N_i$, $\delta_i=N_i\cap\omega_1$ and $\delta_\omega=M\cap\omega_1$. Choose $\gamma\in\Gamma-M$ with $b_\gamma\restriction(\delta_\omega+1)\in T_0$. For each $i\in\omega$, choose $\gamma_i\in(N_{i+1}-N_i)\cap\Gamma$ such that $\Delta_0(\gamma_i,\gamma)>i$ and $\delta_i\leq\Delta(\gamma_i,\gamma)<\delta_{i+1}$. This follows by reflecting $\gamma$ and using the fact that $T_0$ has no uncountable branch. By the argument above, each $\gamma_i$ is pre-good. To see that $\gamma$ is good, fix $m\in\omega$ and let $t=r_{\gamma_m}\restriction(\Delta_0(\gamma_m,\gamma)+1)$.

Fix $\xi<\delta$. Since $\gamma_m$ is pre-good, we can find $\alpha\in\Gamma\cap N$ such that $\xi\leq\Delta(\alpha,\gamma_m)<\delta$ and $\Delta_0(\alpha,\gamma_m)>\Delta_0(\gamma_m,\gamma)+1$. Then $t\sqsubseteq r_\alpha$ and $\xi\leq\Delta(\alpha,\gamma)=\Delta(\alpha,\gamma_m)<\delta$. Thus the values $\Delta(\alpha,\gamma)$ for such $\alpha$ are unbounded in $\delta$.

Since $\Delta_0(t,r_\gamma)=\Delta_0(\gamma_m,\gamma)>m$ and $m$ was arbitrary, $\gamma$ is good.
\end{proof}

Fix a good $\gamma$ as given by the claim, and choose $\beta\in\Gamma$ with $\nu=\Delta(\beta,\gamma)>\delta$. Let $m_0=e_\nu(\delta)$. By goodness, we can find $t\in S$ such that $m=\Delta_0(r_\gamma,t)\geq m_0$ and the values $\Delta(\alpha,\gamma)<\delta$, for $\alpha\in\Gamma\cap N$ with $t\sqsubseteq r_\alpha$, are unbounded in $\delta$. Choose such an $\alpha$ with $\xi=\Delta(\alpha,\gamma)>\sup(F_m(\nu)\cap\delta)$.

We claim that $c(\alpha,\beta,\gamma)=\tau$. Since
\[
\max\{g(\alpha,\beta),g(\alpha,\gamma)\}<\delta\leq\Delta(\beta,\gamma)
\quad\text{and}\quad \Delta(\alpha,\beta)<\Delta(\beta,\gamma),
\]
case \textbf{(2)} of the definition applies, and
\begin{align*}
c(\alpha,\beta,\gamma)
 &=d\bigl(\min(F_{\Delta_0(\alpha,\gamma)}(\Delta(\beta,\gamma))-\Delta(\alpha,\gamma))\bigr)\\
 &=d\bigl(\min(F_m(\nu)\setminus \xi)\bigr)=d(\delta)=\tau.
\end{align*}

\textbf{(3):} Shrinking $\Gamma$ to an uncountable subset if necessary, we may assume that $T_0$ is a countable perfect normal tree satisfying \ref{cantor1} and \ref{cantor2}. Let $\eta=\operatorname{ht}(T_0)$.

Let $N\prec H(\theta)$ be a countable elementary submodel containing the relevant objects, including $T_0$, $\langle e_\beta:\beta<\omega_1\rangle$, $d$ and $g$, such that $\delta=N\cap\omega_1$ satisfies $d(\delta)=\tau$. By the unboundedness of $g$, we can find $\alpha\in\Gamma\cap N$ and $\nu\in A_\alpha$ with $\nu\geq\eta$. Choose $\gamma\in\Gamma-N$ large enough that $g(\alpha,\gamma)=\nu$ and $\xi_1=\operatorname{index}_T(b_\gamma\restriction\nu)>\delta$. Let $m_0=e_{\xi_1}(\delta)$. Since $b_\gamma\restriction\eta\notin N$ and all its proper restrictions belong to $T_0\subseteq N$, there are unboundedly many $\xi<\eta$ for which there is some $t=t_\xi\in T_0$ such that
\begin{itemize}
\item $\Delta(t, b_\gamma)=\xi>\Delta(\alpha,\gamma)$ and 
\item there are uncountably many $\gamma'\in \Gamma$ such that $t\sqsubseteq b_{\gamma'}$ and $g(\alpha,\gamma')=\nu$.
\end{itemize}
Choose such a $\xi<\eta$ and $t=t_\xi$ with $m=e_{\xi_1}(\xi)\geq m_0$. By the property of $t$, find $\beta\in\Gamma\cap N$ above $\alpha$ such that $t\sqsubseteq b_\beta$, $g(\alpha,\beta)=\nu$ and $\xi_0=\operatorname{index}_T(b_\beta\restriction\nu)>\sup(F_m(\xi_1)\cap\delta)$. Since $g(\alpha,\gamma)=g(\alpha,\beta)=\nu\geq\eta>\Delta(\beta,\gamma)$, case \textbf{(3)} of the definition applies, and
\begin{align*}
c(\alpha,\beta,\gamma)
 &=d\bigl(\min(F_{e_{\xi_1}(\Delta(\beta,\gamma))}(\xi_1)-\xi_0)\bigr)\\
 &=d\bigl(\min(F_m(\xi_1)\setminus \xi_0)\bigr)=d(N\cap\omega_1)=\tau.
\end{align*}
\end{proof}

\begin{corollary}
If there exists a Kurepa tree with at least $\kappa$ branches and $\kappa\leq 2^\omega$, then $\kappa\not\to [\omega_1]^3_{\omega_1}$.
\end{corollary}

\begin{remark}
If $\kappa=\omega_2$, then we do not need the extra assumption that $\kappa\leq 2^\omega$. In fact, \cite{Todorcevic3D} shows that $\omega_2\not\to [\omega_1]^3_{\omega_1}$ is equivalent to the failure of Chang's conjecture $(\omega_2,\omega_1) \twoheadrightarrow (\omega_1, \omega)$. The reason why $\kappa=\omega_2$ is special in this case is that a coloring witnessing $\omega_2\not\to (\omega_1)^2_{\omega_1, <\omega_1}$ can be upgraded to another witness that is in addition subadditive (see \cite[Theorem 9.2.3]{walksbook}). Such a function will be \emph{strongly unbounded}. Namely, the coloring restricted to any $\Gamma \otimes\Sigma$ has an uncountable range, where $\Gamma, \Sigma \in [\kappa]^{\omega_1}$ with $\sup \Gamma =\sup \Sigma$. 

In general, if $\kappa>\omega_2$, we cannot upgrade an unbounded function witnessing $\kappa\not\to (\omega_1)^2_{\omega_1, <\omega_1}$ to a \emph{strongly unbounded} one. To see this, suppose there exists a Kurepa tree with $\kappa > (2^\omega)^+$ many branches. We show that for any $c: [\kappa]^2\to \omega_1$, there are $\Sigma,\Gamma\in [\kappa]^{\omega_1}$ having the same supremum such that $c \restriction \Sigma\otimes \Gamma$ is constant. Let $M\prec H(\theta)$ be such that $M$ is countably closed and $M\cap \kappa=\delta\in \kappa\cap \cof(\omega_2)$. Find some unbounded $A\subseteq \delta$ of size $\omega_2$ and $\xi\in \omega_1$ such that $c$ is constantly $\xi$ on $A\times\{\delta\}$. Recursively build $\Sigma=\{\sigma_i\in A: i\in \omega_1\}$ and $\Gamma=\{\gamma_i\in \delta: i\in \omega_1\}$ with the same supremum such that for any $i<j\in \omega_1$, $c(\sigma_i, \gamma_j)=\xi$. At stage $i$, let us find $\sigma_i\in A - (\sup \{\sigma_j, \gamma_j: j<i\}+1)$. Let $\gamma_i>\sigma_i$ be such that for all $k\leq i$, $c(\sigma_k, \gamma_i)=\xi$. The existence of $\gamma_i$ is due to the fact that $\{\sigma_k\in A: k\leq i\}\in M$ and the elementarity of $M$.
\end{remark}

\begin{remark}
Theorem \ref{theorem: weakKurepaDimension2} complements previous results on constructing $3$-dimensional strong colorings from a Kurepa tree under additional assumptions. For example, 
\begin{enumerate}
\item Shore \cite{Shore} showed that if $\mathrm{CH}$ holds and there exists a Kurepa tree with at least $\kappa$ branches, then $\kappa\not\to [\omega_1]^3_{\omega_1}$, 
\item Todorcevic \cite{TodorcevicThesis, TodorcevicWalks87} showed that if there exists a Kurepa tree without an Aronszajn subtree and with at least $\kappa$ branches, then $\kappa\not\to [\omega_1]^3_{\omega_1}$.
\end{enumerate}

\end{remark}

%
%
%

\begin{corollary}\label{cor: CC}
Suppose $n\geq 2$ and there exists a weak Kurepa tree $T$ with at least $\omega_n$ branches. Then the following are equivalent: 
	\begin{enumerate}
	\item $(\omega_n, \omega_1)\twoheadrightarrow (\omega_1, \omega)$,
	\item $\omega_n \to [\omega_1]^{n+1}_{\omega_1}$. 
	\end{enumerate}
\end{corollary}

\begin{proof}
For $n=2$, this is \cite[Theorem 2.1]{Todorcevic3D}. If $n\geq 3$, then the conclusion follows from Theorem \ref{theorem: WKTdimensionN} and the well-known result that $(\omega_n, \omega_1)\twoheadrightarrow (\omega_1, \omega)$ is equivalent to $\omega_n \to (\omega_1)^{n}_{\omega_1, \omega}$ (see \cite[Proposition 8.2]{Kanamori}).
\end{proof}

\begin{remark}
By a theorem of Baumgartner \cite{BaumgartnerAlmostDisjoint}, the hypothesis of Corollary \ref{cor: CC} holds if $2^\omega<\min\{\aleph_{\omega_1},2^{\omega_1}\}$ and $\omega_n\leq 2^{\omega_1}$.
\end{remark}

\begin{remark}
Corollary \ref{cor: CC} is a higher-dimensional extension of an unpublished theorem of Laver, who proved that under $\mathrm{CH}$, $(\omega_2, \omega_1)\twoheadrightarrow (\omega_1, \omega_0)$ is equivalent to $\omega_2 \to [\omega_1]^{3}_{\omega_1}$. As we mentioned above, Todorcevic \cite{Todorcevic3D} removed the use of $\mathrm{CH}$ and derived the conclusion in $\mathrm{ZFC}$. It remains open whether the conclusion of Corollary \ref{cor: CC} can be derived in $\mathrm{ZFC}$ alone.
\end{remark}

\begin{remark}
While the existence of a Kurepa tree with $\kappa$ branches implies the failure of $(\kappa,\omega_1)\twoheadrightarrow (\omega_1, \omega)$, the existence of a \emph{weak} Kurepa tree with $\kappa$ branches does not. Hence we cannot strengthen Theorem \ref{theorem: WKTdimensionN} to conclude that $\kappa\not\to [\omega_1]^{n+1}_{\omega_1}$ follows from the existence of a weak Kurepa tree with at least $\kappa$ branches.

To see this, start with a measurable cardinal $\kappa$ and $n\geq 2$. Consider the forcing $\operatorname{Coll}(\omega_n, <\kappa)\times \operatorname{Add}(\omega_1, \kappa)$. In the forcing extension, $\kappa=\omega_{n+1}$, $(\kappa, \omega_1) \twoheadrightarrow (\omega_1, \omega)$, $2^\omega=\omega_1$ and $2^{<\omega_1}$ is a weak Kurepa tree with at least $\kappa$ branches. It is also possible to arrange configurations where these versions of Chang's conjecture hold simultaneously at all $\omega_n$ for $n\geq 2$ and $2^{\omega_1}>\aleph_\omega$, starting from infinitely many measurable cardinals, for example.
\end{remark}

\section{Concluding remarks and questions}\label{section: questions}

The following questions remain open in the theory of higher-dimensional strong colorings (see also \cite{higherwalks}). For $n\geq 2$,

\begin{itemize}
\item Does $\omega_n\not\to[\omega_1]^{n+1}_{\omega}$ hold?
\item Is $(\omega_n,\omega_1)\not\twoheadrightarrow(\omega_1,\omega)$ equivalent to $\omega_n\not\to[\omega_1]^{n+1}_{\omega_1}$?
\end{itemize}

In this article, we resolved the first question when $n=3$ and made progress on the other problems. We believe investigating these test questions can greatly improve our understanding of the higher-dimensional theory of the cardinals $\omega_k$ and expand the pool of methods and potential applications that we have.


\begin{question}[{\cite[Question 10.3.9]{walksbook}}]\label{question: nonch}
Does $\neg \mathrm{CH}$ imply that $\omega_2\not\to [\omega_1]^2_{\omega}$? 
\end{question}

\begin{question}
For $n\geq 4$, suppose there exists a weak Kurepa tree $T$ with at least $\omega_n$ branches. Is it true that $\omega_n\not\to [\omega_1]^{n+1}_{\omega}$?
\end{question}


The proof of Theorem \ref{theorem: WKTdimensionN} suggests looking for a higher-dimensional $\rho$-function and dividing into cases according to whether it is bounded on a set of size $\aleph_1$. The challenge is to obtain oscillation in the bounded case.

%

\begin{question}
If there exists a Kurepa tree with at least $\kappa$ branches, is $\kappa\not\to [\omega_1]^3_{\omega_1}$ true?
\end{question}

\section{AI disclosure}

ChatGPT 6.2-Astra was used for editorial purposes, for example, transcribing handwritten notes, proofreading and drawing Figure \ref{fig:continuum-fixed-triple} from a hand-drawn sketch.

\bibliographystyle{alpha}
\bibliography{bib}

\begin{thebibliography}{EHMR84}

\bibitem[ALAT05]{ALAT2006}
S.~A. Argyros, J.~Lopez-Abad, and S.~Todorcevic.
\newblock A class of {B}anach spaces with few non-strictly singular operators.
\newblock {\em J. Funct. Anal.}, 222(2):306--384, 2005.

\bibitem[Bau76]{BaumgartnerAlmostDisjoint}
James~E. Baumgartner.
\newblock Almost-disjoint sets, the dense set problem and the partition
  calculus.
\newblock {\em Ann. Math. Logic}, 9(4):401--439, 1976.

\bibitem[Ber21]{BergfalkHigher}
Jeffrey Bergfalk.
\newblock The first omega alephs: from simplices to trees of trees to higher
  walks.
\newblock {\em Adv. Math.}, 393:Paper No. 108083, 74, 2021.

\bibitem[Ber25]{higherwalks}
Jeffrey Bergfalk.
\newblock An introduction to higher walks.
\newblock {\em Zb. Rad. (Beogr.)}, 22(30)(Selected topics in combinatorial
  analysis. Part II):7--71, 2025.

\bibitem[BLHZ26]{BergfalkLambieHansonZhang}
Jeffery Bergfalk, Chris Lambie-Hanson, and Jing Zhang.
\newblock The cohomology of the ordinals.
\newblock {\em in preparation}, 2026.

\bibitem[CFS13]{ImprovedBound}
David Conlon, Jacob Fox, and Benny Sudakov.
\newblock An improved bound for the stepping-up lemma.
\newblock {\em Discrete Appl. Math.}, 161(9):1191--1196, 2013.

\bibitem[COV24]{COV2024}
Samuel~M. Corson, Alexander Olshanskii, and Olga Varghese.
\newblock Steep uncountable groups.
\newblock {\em Bull. Lond. Math. Soc.}, 56(2):701--710, 2024.

\bibitem[Cum05]{CummingsNotes}
James Cummings.
\newblock Notes on singular cardinal combinatorics.
\newblock {\em Notre Dame J. Formal Logic}, 46(3):251--282, 2005.

\bibitem[EH74]{ErdosHajnalProblem}
P.~Erd\H{o}s and A.~Hajnal.
\newblock Unsolved and solved problems in set theory.
\newblock In {\em Proceedings of the {T}arski {S}ymposium ({P}roc. {S}ympos.
  {P}ure {M}ath., {V}ol. {XXV}, {U}niv. {C}alifornia, {B}erkeley, {C}alif.,
  1971)}, Proc. Sympos. Pure Math., Vol. XXV, pages 269--287. Published for the
  Association for Symbolic Logic by the American Mathematical Society,
  Providence, RI, 1974.

\bibitem[EHMR84]{partitionrelationsbook}
Paul Erd\H{o}s, Andr\'{a}s Hajnal, Attila M\'{a}t\'{e}, and Richard Rado.
\newblock {\em Combinatorial set theory: partition relations for cardinals},
  volume 106 of {\em Studies in Logic and the Foundations of Mathematics}.
\newblock North-Holland Publishing Co., Amsterdam, 1984.

\bibitem[EHR65]{ErdosHajnalRado}
P.~Erd\H{o}s, A.~Hajnal, and R.~Rado.
\newblock Partition relations for cardinal numbers.
\newblock {\em Acta Math. Acad. Sci. Hungar.}, 16:93--196, 1965.

\bibitem[Erd47]{erdosprob}
P.~Erd{\H{o}}s.
\newblock Some remarks on the theory of graphs.
\newblock {\em Bull. Amer. Math. Soc.}, 53:292--294, 1947.

\bibitem[FR23]{FeldmanRinot}
Ido Feldman and Assaf Rinot.
\newblock Sums of triples in abelian groups.
\newblock {\em Mathematika}, 69(3):622--664, 2023.

\bibitem[GS73]{GalvinShelah}
Fred Galvin and Saharon Shelah.
\newblock Some counterexamples in the partition calculus.
\newblock {\em J. Combinatorial Theory Ser. A}, 15:167--174, 1973.

\bibitem[GZ21]{GartiZhang}
Shimon Garti and Jing Zhang.
\newblock Stationary and closed rainbow subsets.
\newblock {\em Ann. Pure Appl. Logic}, 172(2):Paper No. 102887, 16, 2021.

\bibitem[HRST21]{HRST2021}
Petr H\'{a}jek, Tommaso Russo, Jacopo Somaglia, and Stevo Todor\v{c}evi\'{c}.
\newblock An {A}splund space with norming {M}arku\v{s}evi\v{c} basis that is
  not weakly compactly generated.
\newblock {\em Adv. Math.}, 392:Paper No. 108041, 22, 2021.

\bibitem[Kan09]{Kanamori}
Akihiro Kanamori.
\newblock {\em The higher infinite}.
\newblock Springer Monographs in Mathematics. Springer-Verlag, Berlin, second
  edition, 2009.
\newblock Large cardinals in set theory from their beginnings, Paperback
  reprint of the 2003 edition.

\bibitem[Kom87]{Komjath}
P.~Komj\'{a}th.
\newblock Morasses and the {L}\'{e}vy-collapse.
\newblock {\em J. Symbolic Logic}, 52(1):111--115, 1987.

\bibitem[KPT05]{KPT}
A.~S. Kechris, V.~G. Pestov, and S.~Todorcevic.
\newblock Fra\"{\i}ss\'{e} limits, {R}amsey theory, and topological dynamics of
  automorphism groups.
\newblock {\em Geom. Funct. Anal.}, 15(1):106--189, 2005.

\bibitem[LAT13]{LAT2011}
Jordi Lopez-Abad and Stevo Todorcevic.
\newblock Positional graphs and conditional structure of weakly null sequences.
\newblock {\em Adv. Math.}, 242:163--186, 2013.

\bibitem[Lav75]{LaverParition}
R.~Laver.
\newblock Partition relations for uncountable cardinals {$\leq 2\sp{\aleph
  \sb{0}}$}.
\newblock In {\em Infinite and finite sets ({C}olloq., {K}eszthely, 1973;
  dedicated to {P}. {E}rd\H{o}s on his 60th birthday), {V}ols. {I}, {II},
  {III}}, Colloq. Math. Soc. J\'{a}nos Bolyai, Vol. 10, pages 1029--1042.
  North-Holland, Amsterdam-London, 1975.

\bibitem[LRT24]{HLRandST}
Hossein Lamei~Ramandi and Stevo Todorcevic.
\newblock Can you take {K}omjath's inaccessible away?
\newblock {\em Ann. Pure Appl. Logic}, 175(7):Paper No. 103452, 24, 2024.

\bibitem[Moo06]{MooreLspace}
Justin~Tatch Moore.
\newblock A solution to the {$L$} space problem.
\newblock {\em J. Amer. Math. Soc.}, 19(3):717--736, 2006.

\bibitem[PR25]{PR2025}
M\'{a}rk Po\'{o}r and Assaf Rinot.
\newblock A {S}helah group in {ZFC}.
\newblock {\em Forum Math. Pi}, 13:Paper No. e13, 28, 2025.

\bibitem[Ram29]{Ra30}
F.~P. Ramsey.
\newblock On a {P}roblem of {F}ormal {L}ogic.
\newblock {\em Proc. London Math. Soc. (2)}, 30(4):264--286, 1929.

\bibitem[RZ21]{transformation}
Assaf Rinot and Jing Zhang.
\newblock Transformations of the transfinite plane.
\newblock {\em Forum Math. Sigma}, 9:Paper No. e16, 25, 2021.

\bibitem[She26]{shelahPositive}
Saharon Shelah.
\newblock Consistency of square bracket partition relation, 2026.

\bibitem[Sho74]{Shore}
Richard~A. Shore.
\newblock Square bracket partition relations in {$L$}.
\newblock {\em Fund. Math.}, 84(2):101--106, 1974.

\bibitem[Sie33]{Sierpinski}
Waclaw Sierpi\'{n}ski.
\newblock Sur un probl\`eme de la th\'{e}orie des relations.
\newblock {\em Ann. Scuola Norm. Super. Pisa Cl. Sci. (2)}, 2(3):285--287,
  1933.

\bibitem[Tod81]{TodorcevicThesis}
Stevo~B. Todor\v{c}evi\'{c}.
\newblock Trees, subtrees and order types.
\newblock {\em Ann. Math. Logic}, 20(3):233--268, 1981.

\bibitem[Tod85]{TodorcevicSpecial}
Stevo Todor\v{c}evi\'{c}.
\newblock Aronszajn trees and partitions.
\newblock {\em Israel J. Math.}, 52(1-2):53--58, 1985.

\bibitem[Tod87]{TodorcevicWalks87}
Stevo Todor\v{c}evi\'{c}.
\newblock Partitioning pairs of countable ordinals.
\newblock {\em Acta Math.}, 159(3-4):261--294, 1987.

\bibitem[Tod94]{Todorcevic3D}
Stevo Todor\v{c}evi\'{c}.
\newblock Some partitions of three-dimensional combinatorial cubes.
\newblock {\em J. Combin. Theory Ser. A}, 68(2):410--437, 1994.

\bibitem[Tod07]{walksbook}
Stevo Todorcevic.
\newblock {\em Walks on ordinals and their characteristics}, volume 263 of {\em
  Progress in Mathematics}.
\newblock Birkh\"{a}user Verlag, Basel, 2007.

\bibitem[TZ26]{TodorcevicZhangOsc}
Stevo Todorcevic and Jing Zhang.
\newblock Strong colorings based on oscillations, 2026.

\end{thebibliography}

\end{document}